\documentclass[reqno]{amsart}

\usepackage{amssymb,amscd}

\usepackage{pdfsync}

       \numberwithin{equation}{section}

\theoremstyle{plain}

\newtheorem{Thm}{Theorem}[section]
\newtheorem{Lem}[Thm]{Lemma}
\newtheorem{Cor}[Thm]{Corollary}
\newtheorem{Prop}[Thm]{Proposition}
\newtheorem{MainThm}{Theorem}
\newtheorem{MainCor}[MainThm]{Corollary}

\theoremstyle{definition}

\newtheorem{Def}[Thm]{Definition}

\newtheorem{exmp}[Thm]{Example}
\newtheorem{prob}[Thm]{Problem}

\theoremstyle{remark}

\newtheorem{Rem}{Remark}
\newtheorem{Rems}{Remarks}
\newtheorem{Cla}{Claim}

\newtheorem{notation}{Notation}
\newtheorem{ack}{Acknowledgment}

\renewcommand{\sec}[2]{\section{#2}\label{s:#1}}
\newcommand{\ssec}[2]{\subsection{#2}\label{ss:#1}}

\newcommand{\refs}[1]{Section~\ref{s:#1}}
\newcommand{\refss}[1]{Section~\ref{ss:#1}}
\newcommand{\reft}[1]{Theorem~\ref{t:#1}}
\newcommand{\refl}[1]{Lemma~\ref{l:#1}}
\newcommand{\refp}[1]{Proposition~\ref{p:#1}}
\newcommand{\refc}[1]{Corollary~\ref{c:#1}}
\newcommand{\refd}[1]{Definition~\ref{d:#1}}
\newcommand{\refr}[1]{Remark~\ref{r:#1}}
\newcommand{\refcl}[1]{Claim~\ref{cl:#1}}
\newcommand{\refe}[1]{\eqref{e:#1}}

\newcommand{\refts}[2]{Theorems~\ref{t:#1} and~\ref{t:#2}}
\newcommand{\refmt}[1]{Theorem~\ref{mt:#1}}

\newcommand{\refls}[2]{Lemmas~\ref{l:#1} and~\ref{l:#2}}
\newcommand{\refcs}[2]{Corollaries~\ref{c:#1} and~\ref{c:#2}}
\newcommand{\refmc}[1]{Corollary~\ref{mc:#1}}
\newcommand{\refmcs}[2]{Corollaries~\ref{mc:#1} and~\ref{mc:#2}}

\newcommand{\refcls}[2]{Claims~\ref{cl:#1} and~\ref{cl:#2}}

\newenvironment{thm}[1]{\begin{Thm} \label{t:#1}}{\end{Thm}}
\renewcommand{\th}[1]{\begin{thm}{#1}}
\renewcommand{\eth}{\end{thm}}

\newenvironment{mainthm}[1]{\begin{MainThm}\label{mt:#1}}{\end{MainThm}}
\newcommand{\mainth}[1]{\begin{mainthm}{#1}}
\newcommand{\emainth}{\end{mainthm}}

\newenvironment{lemma}[1]{\begin{Lem}\label{l:#1}}{\end{Lem}}
\newcommand{\lem}[1]{\begin{lemma}{#1}}
\newcommand{\elem}{\end{lemma}}

\newenvironment{propos}[1]{\begin{Prop}\label{p:#1}}{\end{Prop}}
\newcommand{\prop}[1]{\begin{propos}{#1}}
\newcommand{\eprop}{\end{propos}}

\newenvironment{corol}[1]{\begin{Cor}\label{c:#1}}{\end{Cor}}
\newcommand{\cor}[1]{\begin{corol}{#1}}
\newcommand{\ecor}{\end{corol}}

\newenvironment{maincorol}[1]{\begin{MainCor}\label{mc:#1}}{\end{MainCor}}
\newcommand{\maincor}[1]{\begin{maincorol}{#1}}
\newcommand{\emaincor}{\end{maincorol}}

\newenvironment{defini}[1]{\begin{Def}\label{d:#1}}{\end{Def}}
\newcommand{\defi}[1]{\begin{defini}{#1}}
\newcommand{\edefi}{\end{defini}}

\newenvironment{remark}[1]{\begin{Rem}\label{r:#1}}{\end{Rem}}
\newcommand{\rem}[1]{\begin{remark}{#1}}
\newcommand{\erem}{\end{remark}}

\newenvironment{remarks}[1]{\begin{Rems}\label{r:#1}}{\end{Rems}}
\newcommand{\rems}[1]{\begin{remarks}{#1}}
\newcommand{\erems}{\end{remarks}}

\newenvironment{claim}[1]{\begin{Cla}\label{cl:#1}}{\end{Cla}}
\newcommand{\cla}[1]{\begin{claim}{#1}}
\newcommand{\ecla}{\end{claim}}

\newcommand{\prf}{\begin{proof}}
\newcommand{\eprf}{\end{proof}}

\newcommand{\EE}{\text{$\mathcal{E}$}}
\newcommand\FF{\text{$\mathcal{F}$}}
\newcommand{\GG}{\text{$\mathcal{G}$}}
\newcommand{\HH}{\text{$\mathcal{H}$}}

\newcommand{\KK}{\text{$\mathcal{K}$}}

\newcommand{\MM}{\text{$\mathcal{M}$}}
\newcommand{\NN}{\text{$\mathcal{N}$}}
\newcommand{\OO}{\text{$\mathcal{O}$}}
\newcommand{\PP}{\text{$\mathcal{P}$}}
\newcommand{\QQ}{\text{$\mathcal{Q}$}}

\newcommand{\TT}{\text{$\mathcal{T}$}}
\newcommand{\UU}{\text{$\mathcal{U}$}}
\newcommand{\VV}{\text{$\mathcal{V}$}}
\newcommand{\WW}{\text{$\mathcal{W}$}}

\newcommand{\YY}{\text{$\mathcal{Y}$}}

\newcommand{\frakg}{\text{$\mathfrak{g}$}}

\newcommand{\frakk}{\text{$\mathfrak{k}$}}

\newcommand{\frakr}{\text{$\mathfrak{r}$}}

\newcommand{\frakG}{\text{$\mathfrak{G}$}}

\newcommand{\frakX}{\text{$\mathfrak{X}$}}

        \newcommand{\field}[1]{\text{$\mathbb{#1}$}}
        \newcommand{\N}{\field{N}}
        \newcommand{\Z}{\field{Z}}
        
        \newcommand{\R}{\field{R}}

\newcommand{\supp}{\operatorname{supp}}
\newcommand{\codim}{\operatorname{codim}}
\newcommand{\im}{\operatorname{im}}
\newcommand{\dom}{\operatorname{dom}}

\newcommand{\id}{\operatorname{id}}
\newcommand{\diam}{\operatorname{diam}}
\newcommand{\Pen}{\operatorname{Pen}}
\newcommand{\Top}{\operatorname{Top}}
\newcommand{\PsGr}{\operatorname{PsGr}}
\newcommand{\Fol}{\operatorname{Fol}}
\newcommand{\Hol}{\operatorname{Hol}}

\newcommand{\cat}{\operatorname{cat}}

\newcommand{\Man}{\operatorname{Man}}
\newcommand{\length}{\operatorname{length}}
\newcommand{\pr}{\operatorname{pr}}
\newcommand{\Exp}{\operatorname{Exp}}

\newcommand{\Or}{\operatorname{O}}

\newcommand{\Proper}{\operatorname{Prop}}

\newcommand{\cinf}{\text{$C^\infty$}}
\newcommand{\coinf}{\text{$C^{0,\infty}$}}
\newcommand{\cinfo}{\text{$C^{\infty,0}$}}
\newcommand{\lar}{\text{$\;\longrightarrow\;$}}
\newcommand{\sm}{\smallsetminus}

\newcommand{\germ}[2]{\gamma({#1},{#2})}
\newcommand{\germs}[1]{\gamma({#1})}
\newcommand{\jet}[1]{j^1(#1)}
\newcommand{\jetgerm}[2]{\jet{\germ{#1}{#2}}}
\newcommand{\jetgerms}[1]{\jet{\germs{#1}}}

\begin{document}

\bibliographystyle{plain}

%-------------------------------------------------------------------
%-------------------------------------------------------------------

\title[Morphisms between pseudogroups]{Morphisms between complete Riemannian pseudogroups}

\author[J.A. \'Alvarez L\'opez]{Jes\'us A. \'Alvarez L\'opez}

\author[X.M. Masa]{Xos\'e M. Masa}

\thanks{Partially supported by MEC (Spain), grant~MTM2004-08214}

\address{Departamento de Xeometr\'{\i}a e Topolox\'{\i}a\\
         Facultade de Matem\'aticas\\
         Universidade de Santiago de Compostela\\
         15706 Santiago de Compostela\\ Spain}

\email[J.A. \'Alvarez L\'opez]{jalvarez@usc.es}

\email[X.M. Masa]{xosemasa@usc.es}

\begin{abstract}
We introduce the concept of morphism of pseudogroups generalizing the \'etal\'e morphisms of Haefliger. With our definition, any continuous foliated map induces a morphism between the corresponding holonomy pseudogroups. The main theorem states that any morphism between complete Riemannian pseudogroups is complete, has a closure and its maps are \cinf\ along the orbit closures. Here, completeness and closure are versions for morphisms of concepts introduced by Haefliger for pseudogroups. This result is applied to approximate foliated maps by smooth ones in the case of transversely complete Riemannian foliations, yielding the foliated homotopy invariance of their spectral sequence. This generalizes the topological invariance of their basic cohomology, shown by El~Kacimi-Nicolau. A different proof of the spectral sequence invariance was also given by the second author.
\end{abstract}



\maketitle

\tableofcontents

%--------------------------------------------------------------
%--------------------------------------------------------------

\sec{intro}{Introduction and main results}

We introduce the concept of {\em morphism\/} of pseudogroups generalizing the {\em \'etal\'e morphisms\/} of Haefliger \cite{Haefliger:88Leaf}: it is a slight change of Haefliger's definition so that they may contain, not only homeomorphisms, but also arbitrary continuous maps. With our definition, any continuous map between foliated spaces mapping leaves to leaves (a {\em foliated map\/}) induces a morphism between the corresponding holonomy pseudogroups; indeed, this defines a functor between the categories of foliated maps and morphisms of pseudogroups, called the {\em holonomy functor\/}. A morphism of pseudogroups can be interpreted as a {\em morphism\/} of {\em $S$-atlases\/}, defined by Van~Est \cite{vanEst:84}, and is more general than a homomorphism of \'etal\'e groupoids: only transverse foliated maps induce homomorphisms of holonomy groupoids. Any pseudogroup morphism induces a continuous map between the corresponding orbit spaces.

Haefliger has also introduced other morphisms between topological groupoids \cite{Haefliger:84classifiants}. Any Haefliger morphism between \'etal\'e groupoids induces a morphism of our type (see Example~\ref{ex:Haefliger}). Somehow, our morphisms extract the information about local maps from Haefliger morphisms, and forget about their additional structure. Since our main results will be about maps, our morphisms seems to be more appropriate for our study. 

For each topological space, we can consider the pseudogroup generated by its identity map. Then, each continuous map between topological spaces {\em generates\/} (in an obvious sense) a morphism between the corresponding pseudogroups. This assignment defines a canonical injective functor between the category of continuous maps and the category of morphisms. In this sense, many topological and geometric concepts can be generalized to pseudogroups. For instance, the concepts of {\em homotopy\/} and {\em homotopy equivalence\/} will be used for pseudogroups. The definition of the {\em fundamental group\/} of a pseudogroup with a distinguished orbit was given by Haefliger \cite{Haefliger:85Pseudogroups} and Van~Est \cite{vanEst:84}.  For a pseudogroup consisting of \cinf\ transformations, its {\em de~Rham cohomology\/} is isomorphic to its {\em invariant cohomology\/} (the cohomology of the complex of differential forms that are invariant by the local transformations of the pseudogroup).

We will mainly consider morphisms between pseudogroups of local isometries of Riemannian manifolds, called {\em Riemannian pseudogroups\/}. In \cite{Haefliger:88Leaf}, Haefliger has introduced the condition of {\em completeness\/} for pseudogroups, and the {\em closure\/} of a complete Riemannian pseudogroup. Now, we give versions of these concepts for morphisms.

Recall that the orbit closures of a complete Riemannian pseudogroup are the leaves of a \cinf\ singular foliation \cite{Salem:88}. Then a morphism between such pseudogroups is said to be of {\em class \coinf\/} when it consists of maps that are \cinf\ along the orbit closures, and moreover the corresponding leafwise derivatives of arbitrary order are continuous on the ambient manifold; this makes sense even for singular foliations! Our main result is the following.

\mainth{main}
Any morphism between complete Riemannian pseudogroups is complete, has a closure and is of class \coinf.
\emainth

The last part of Theorem~\ref{mt:main} is a generalization of the well known fact that any continuous homomorphism between Lie groups is \cinf. Thus a version of last part of Theorem~\ref{mt:main} for ``measurable morphisms'' seems to be possible (Problem~\ref{pb:measurable morphisms}).

The proof of \refmt{main} only involves basic techniques, but it is rather complicated and has a lot of steps. A simplified version of the proof, without details and only for the case of dense orbits, was given in \cite{AlvMasa:Lodz}; it may be useful to understand the complete proof of this paper.

A {\em Riemannian foliation\/} is a \cinf\ foliation whose holonomy pseudogroup is Riemannian. These foliations can be locally described by Riemannian submersions for some Riemannian metric (a {\em bundle-like metric\/}) \cite{Reinhart:59bundle-like}. A Riemannian foliation is said to be {\em transversely complete\/} when, for some bundle-like metric, the geodesics orthogonal to the leaves are complete; this condition implies that its holonomy pseudogroup is complete \cite[Example~1.2.1]{Haefliger:88Leaf}. 

We will give examples of non-complete morphisms between complete pseudo\-groups (Examples~\ref{ex:product} and~\ref{ex:R2}), showing that the Riemannian condition cannot be removed in \refmt{main}. However, a weaker equicontinuity condition could be enough to get completeness and the existence of a closure (see \cite[Appendix~E]{Molino:88Riemannian} and \cite{AlvCandel:equicont}), although it is not enough to get smoothness along the orbit closures (Example~\ref{ex:Arnold}).

Certain topology, called the {\em strong adapted topology\/}, is introduced on the set of foliated maps between transversely complete Riemannian foliations. Its definition takes into account the special structure of these foliations \cite{Molino:88Riemannian}, but it equals the strong compact-open topology when the leaf closures are compact. 

A homotopy consisting of foliated maps is called a  {\em foliated homotopy\/}, and the corresponding concept of {\em foliated homotopy equivalence\/} can be defined too. The strong adapted topology has the following nice behavior with respect to foliated homotopies.

\mainth{foliated homotopy}
If two foliated maps between transversely complete Riemannian foliations are close enough with respect to the strong adapted topology, then there is a foliated homotopy between them. Moreover the foliated homotopy can be chosen to be proper if the foliated maps are proper.
\emainth

The following result gives smooth approximations of continuous foliated maps in the case of transversely complete Riemannian foliations. It is easy to find counterexamples showing that it cannot be generalized to arbitrary \cinf\ foliations.

\mainth{cinf}
In the space of continuous foliated maps between two transversely complete Riemannian foliations with the strong adapted topology, the subset of \cinf\ foliated maps is dense.
\emainth

Indeed, a relative version of \refmt{cinf} is proved, which, together with \refmt{foliated homotopy}, has the following consequences.

\maincor{cinf foliated homotopy}
If two continuous foliated maps between transversely complete Riemannian foliations are foliatedly homotopic, then there is a \cinf\ foliated homotopy between them.
\emaincor

\maincor{foliated homotopy}
Any continuous foliated map between transversely complete Riemannian foliations is foliatedly homotopic to a \cinf\ foliated map.
\emaincor

In the proof of \refmt{cinf}, the transverse smoothness inside the leaf closures is granted by \refmt{main} for any foliated map. Then standard arguments are used to improve differentiability, firstly, along the leaves and, secondly, along the geodesics orthogonal to the leaf closures. In particular, to improve leafwise differentiability, certain topology is introduced on the set of foliation maps, called the {\em strong plaquewise topology\/}, whose study is continued in \cite{AlvKordyukov:transversality}.

The cohomological information of a \cinf\ manifold $M$ can be enlarged by the presence of a \cinf\ foliation \FF. Precisely, \FF\ gives rise to the {\em spectral sequence\/} $(E_i=E_i(\FF),d_i)$ \cite{Masa:85Lie}, \cite{Sarkaria:78finite}, \cite{KamberTondeur:83Fm}, \cite{Alv:89finiteness}, which is the obvious generalization of the de~Rham version of the Leray spectral sequence of a \cinf\ fiber bundle: it is induced by a filtration of the de~Rham complex of $M$ locally defined like for fiber bundles. Several cohomologies usually associated to \cinf\ foliations are contained in this spectral sequence. For instance, the terms $E_1^{0,\bullet}$ and $E_2^{\bullet,0}$ are respectively called {\em leafwise cohomology\/} and {\em basic cohomology\/}.

Another version of the spectral sequence, $(E_{c,i},d_i)$, can be defined by considering compactly supported differential forms. The term $E_{c,2}^{\bullet,p}$ ($p=\dim\FF$) is isomorphic to the {\em transverse cohomology\/} \cite{Haefliger:80minimal}, also called {\em Haefliger cohomology\/}.

In general, $(E_i,d_i)$ is not a topological invariant (a counterexample is given in \cite{KacimiNicolau:93topological}). Nevertheless, by \refmcs{cinf foliated homotopy}{foliated homotopy}, any continuous foliated map between transversely complete Riemannian foliations induces a homomorphism between their spectral sequences, yielding the following spectral sequence invariance:

\maincor{invariance}
Any foliated homotopy equivalence between transversely complete Riemannian foliations induces an isomorphism on $E_i$ for $i\ge2$.
\emaincor

In the case of dense leaves, there is a version of \refmc{invariance} for $E_1$.

\refmc{invariance} generalizes the topological invariance of the basic cohomology for Riemannian foliations on closed manifolds, obtained by El~Kacimi-Nicolau \cite{KacimiNicolau:93topological}. 

The second author has given a different proof of \refmc{invariance} in \cite{Masa:02A-S}. He introduces an Alexander-Spanier version of the spectral sequence, which is a topological invariant by definition, and indeed is shown to be invariant by foliated homotopy equivalences. Then, for transversely complete Riemannian foliations, he shows that both spectral sequences are isomorphic from the second term on.

Versions of the above approximation and invariance results are also proved for proper foliated maps: 

\mainth{proper cinf}
In the space of continuous foliated maps between two transversely complete Riemannian foliations with the strong adapted topology, the subset of proper continuous foliated maps is open, and therefore the subset of proper \cinf\ foliated maps is dense.
\emainth

\maincor{proper cinf foliated homotopy}
For transversely complete Riemannian foliations, if two proper continuous foliated maps are foliatedly homotopic, then there is a proper \cinf\ foliated homotopy between them.
\emaincor

\maincor{proper foliated homotopy}
For transversely complete Riemannian foliations, there is a proper foliated homotopy between any proper continuous foliated map and some proper \cinf\ foliated map.
\emaincor

\maincor{proper invariance}
Any proper foliated homotopy equivalence between transversely complete Riemannian foliations induces an isomorphism on $E_{c,i}$ for $i\ge2$.
\emaincor

We also introduce a version of the strong compact-open topology on the set of morphisms between two pseudogroups. Then the above approximation and invariance results have the following versions for morphisms:

\mainth{psgr cinf}
In the space of morphisms between two complete Riemannian pseudogroups with the strong compact-open topology, the subset of \cinf\ morphisms is dense.
\emainth

\maincor{psgr cinf homotopy}
If two morphisms between complete Riemannian pseudogroups are homotopic, then there is a \cinf\ homotopy between them.
\emaincor

\maincor{psgr homotopy}
Any morphism between complete Riemannian pseudogroups is homotopic to a \cinf\ morphism.
\emaincor

\maincor{psgr invariance}
Any homotopy equivalence between complete Riemannian pseudo\-groups induces an isomorphism between their invariant cohomologies.
\emaincor

In Sections~\ref{s:morphisms}--\ref{s:existence}, besides introducing some new concepts, we recall the needed preliminaries about pseudogroups and foliated spaces. Special emphasis is put on the case of complete Riemannian pseudogroups and transversely complete Riemannian foliations. Some of the preliminaries are elaborated to fit our needs. \refmt{main} is proved in Sections~\ref{s:main} and~\ref{s:proof}. Sections~\ref{s:SP}--\ref{s:cinf approx} are devoted to the approximation results for foliated maps. The preliminaries on the spectral sequence are given in \refs{spectral seq}, and its invariance is shown in \refs{invariance}. Sections~\ref{s:S for morphisms}-\ref{s:invariance for psgrs} are devoted to the pseudogroup versions of the approximation and invariance results. Finally, examples and open problems are given in the last two sections.

\begin{ack}
We thank A.~Haefliger and the referee for interesting remarks.
\end{ack}

\sec{morphisms}{Morphisms of pseudogroups}

Recall the following definitions from \cite[Section~1.4]{Haefliger:88Leaf}. A {\em
pseudogroup\/} of local transformations of a topological space $T$ (or {\em
acting\/} on $T$) is a collection $\HH$ of homeomorphisms between open subsets of
$T$ which contains the identity map $\id_T$ of $T$, and is closed under
the operations of composition (wherever defined), inversion, restriction to open sets
and combination (defining homeomorphisms). Here, the composition $h_2\circ h_1$ of two local transformations $h_i:U_i\to V_i$, $i=1,2$, refers to the composite of the maps
$$
\begin{CD}
h_1^{-1}(V_1\cap U_2) @>{h_1}>> V_1\cap U_2 @>{h_2}>> h_2(V_1\cap U_2)\;.
\end{CD}
$$
On the other hand, being closed by combination means that, if a family of maps in \HH\ can be combined to define a homeomorphism, then this combination is in \HH. Sometimes it will be convenient to consider the elements of \HH\ as open embeddings with target space $T$.
The pseudogroup
$\HH$ is said to be {\em generated\/} by a subset
$S\subset\HH$ if any map in $\HH$ can be obtained from maps in $S$
involving the above operations.
The {\em restriction\/} of $\HH$ to a subspace $T_0\subset
T$ is the pseudogroup
$\HH|_{T_0}$ consisting of the homeomorphisms between open subsets of $T_0$ that can be locally extended to maps in \HH. If $T_0$ is open in $T$, then $\HH|_{T_0}$
consists of the maps in $\HH$ whose domain and image is contained in $T_0$. The
{\em orbit\/} of a point $x\in T$ is the set $\HH(x)$ of the images of $x$ by all
maps in $\HH$ whose domains contain $x$. The union of orbits that meet some subset
$T_0\subset T$ is called the {\em saturation\/} of $T_0$ and denoted by
$\HH(T_0)$; the set $T_0$ is said to be {\em invariant\/} or {\em saturated\/} if
$T_0=\HH(T_0)$. The quotient space of orbits is denoted by
$\HH\backslash T$.  The pseudogroup \HH\ is said to be of {\em class \cinf\/} if $T$ is a \cinf\ manifold and \HH\ consists of \cinf\ maps. 

Suppose that $T=\bigcup_iT_i$. Given a pseudogroup $\HH_i$ acting on each $T_i$, there is a pseudogroup \HH\ acting on $T$ whose elements are the homeomorphisms between open subsets of $T$, $h:U\to V$, such that $h(U\cap T_i)=V\cap T_i$ and the restriction $h:U\cap T_i\to V\cap T_i$ is in $\HH_i$ for all $i$. If the maps in each $\HH_i$ can be locally extended to maps in \HH, then $\HH|_{T_i}=\HH_i$ for all $i$, and \HH\ is called the {\em combination\/} of the pseudogroups $\HH_i$.  When every $T_i$ is open in $T$, then the combination of the pseudogroups $\HH_i$ is defined just when $\HH_i|_{T_i\cap T_j}=\HH_j|_{T_i\cap T_j}$ for all $i$ and $j$.

Let \HH\ and $\HH'$ be pseudogroups of local transformations of spaces $T$ and $T'$, respectively.
According to \cite[Section~1.4]{Haefliger:88Leaf}, an {\em \'etal\'e morphism\/} $\Phi:\HH\to\HH'$ is a
maximal collection of homeomorphisms of open subsets of $T$ to open subsets of $T'$
satisfying the following conditions:
\begin{itemize}

\item[(i)] If $\phi\in\Phi$, $h\in\HH$ and $h'\in\HH'$, then $h'\circ\phi\circ
h\in\Phi$.

\item[(ii)] The domains of elements of $\Phi$ cover $T$.

\item[(iii)] If $\phi,\psi\in\Phi$, then $\psi\circ\phi^{-1}\in\HH'$.

\end{itemize}
If moreover $\Phi^{-1}=\{\phi^{-1}\ |\ \phi\in\Phi\}$ is an \'etal\'e morphism, then $\Phi$ is called an {\em equivalence\/} of pseudogroups, and $\HH$ is said to be {\em equivalent\/} to $\HH'$. If $\Phi:\HH\to\HH'$ is an equivalence, then $\HH''=\HH\cup\HH'\cup\Phi\cup\Phi^{-1}$ is a pseudogroup on the topological sum $T''=T\sqcup T'$ whose orbits cut $T$ and $T'$, and so that $\HH''|_T=\HH$ and $\HH''|_{T'}=\HH'$. Reciprocally, for any pseudogroup $\HH''$ on $T''$ satisfying these conditions, the family
$$
\{\phi\in\HH''\ |\ \dom\phi\subset T,\ \im\phi\subset T'\}
$$
is an equivalence $\HH\to\HH'$. 

We want to generalize \'etal\'e morphisms by involving arbitrary local maps, and thus the above condition~(iii) must be weakened accordingly. Precisely, we introduce the following concept. 

\defi{morphism}
A {\em morphism\/} $\Phi:\HH\to\HH'$ is a maximal collection of continuous maps of open subsets of $T$ to $T'$ satisfying the following properties:
\begin{itemize}

\item[(i)] If $\phi\in\Phi$, $h\in\HH$ and $h'\in\HH'$, then $h'\circ\phi\circ h\in\Phi$.

\item[(ii)] The domains of elements of $\Phi$ cover $T$.

\item[(iii)] If $\phi,\psi\in\Phi$ and $x\in\dom\phi\cap\dom\psi$, then there is some $h'\in\HH'$ with
$\phi(x)\in\dom h'$ and so that $h'\circ\phi=\psi$ on some neighborhood of $x$.

\end{itemize}
\edefi

\rems{morphism} 
In \refd{morphism}, observe the following:
\begin{itemize}

\item[(a)] Property~(i) implies that $\Phi$ is closed by restrictions to open sets. 

\item[(b)] Properties~(i) and~(iii) imply that, if $\phi,\psi\in\Phi$,
$h\in\HH$ and $x\in\dom\phi\cap h^{-1}(\dom\psi)$, then there is some
$h'\in\HH'$ with $\phi(x)\in\dom h'$ and so that $h'\circ\phi=\psi\circ h$ on some neighborhood of $x$.

\item[(c)] $\Phi$ induces a continuous map $\Phi_{\text{\rm orb}}:\HH\backslash T\to\HH'\backslash T'$ defined by 
$$
\Phi_{\text{\rm orb}}(\OO)=\bigcup_{\phi\in\Phi}\phi(\OO\cap\dom\phi)\;.
$$

\end{itemize}
\erems

The set of morphisms $\HH\to\HH'$ will be denoted by $C(\HH,\HH')$. A morphism $\Phi:\HH\to\HH'$ is said to be of {\em class \cinf\/} if $\HH$ and $\HH'$ are \cinf\ pseudogroups and $\Phi$ consists of \cinf\ maps; the set of \cinf\ morphisms $\HH\to\HH'$ will be denoted by $\cinf(\HH,\HH')$. According to \refr{morphism}-(c), the mapping $\Phi\mapsto\Phi_{\text{\rm orb}}$ defines a canonical map $C(\HH,\HH')\to C(\HH\backslash T,\HH'\backslash T')$ called the {\em orbit map\/}.

\lem{combination}
Let $\Phi$ be collection of continuous maps of open subsets of $T$
to $T'$ satisfying the properties~\upn{(}i\upn{)},~\upn{(}ii\upn{)} and~\upn{(}iii\upn{)} of \refd{morphism}.
Then 
$\Phi$ is a morphism $\HH\to\HH'$ if and only if $\Phi$ is closed under combination of maps.
\elem

\prf
Suppose that $\Phi$ is a morphism $\HH\to\HH'$. The family $\Psi$ of all possible combinations of
maps in $\Phi$ contains $\Phi$ and also satisfies properties~(i),~(ii) and~(iii) of \refd{morphism}. Hence
$\Phi=\Psi$ by the maximality of $\Phi$.

Now, assume that $\Phi$ is closed under combination of maps and let us show the maximality of
\refd{morphism}. Suppose $\Phi$ is contained in another collection $\Psi$ of continuous maps of
open subsets of $T$ to $T'$ satisfying the properties~(i),~(ii) and~(iii) of \refd{morphism}. Take any
$\psi\in\Psi$. By property~(ii) for $\Phi$ and property~(iii) for
$\Psi$, for every $x\in\dom\psi$ there is some $\phi\in\Phi$ with $x\in\dom\phi$ and there is some
$h'\in\HH'$ with $\phi(x)\in\dom h'$ and such that
$h'\circ\phi=\psi$ on some neighborhood of $x$. But $h'\circ\phi\in\Phi$ because $\Phi$ satisfies property~(i).
Therefore every germ of
$\psi$ is the germ of some element of
$\Phi$, yielding that $\psi$ is a combination of elements of $\Phi$. Thus $\psi\in\Phi$ as desired because
$\Phi$ is closed under combination of maps.
\eprf

\lem{morphism}
Let $\Phi_0$ be a family of continuous maps of open subsets of $T$ to $T'$ satisfying the following
properties:
\begin{itemize}

\item[(ii')] The \HH-saturation of the domains of maps in $\Phi_0$ cover $T$.

\item[(iii')] There is a subset $S$ of generators of $\HH$ such that, if $\phi,\psi\in\Phi_0$,
$h\in S$ and $x\in\dom\phi\cap h^{-1}(\dom\psi)$, then there is some
$h'\in\HH'$ with $\phi(x)\in\dom h'$ and so that $h'\circ\phi=\psi\circ h$ on some neighborhood of $x$.

\end{itemize}
Then there is a unique morphism $\Phi:\HH\to\HH'$ containing $\Phi_0$.
\elem

\prf Let $\Phi$ be the family of the following types of maps of open subsets of $T$ to $T'$:
\begin{itemize}

\item All composites $h'\circ\phi\circ h$ with $\phi\in\Phi_0$, $h\in\HH$ and $h'\in\HH'$, wherever defined.

\item All possible combinations of composites of the above type.

\end{itemize}
This collection $\Phi$ clearly satisfies properties~(i),~(ii) and~(iii) of \refd{morphism}.  Moreover $\Phi$
is closed under combination of maps, and thus it is a morphism by \refl{combination}.

Finally, the uniqueness of $\Phi$ follows because, if $\Psi$ is another morphism $\HH\to\HH'$ containing $\Phi_0$, then it also contains $\Phi$ by property~(i) for $\Psi$ and its maximality, and thus $\Phi=\Psi$ by the maximality of $\Phi$.
\eprf

In \refl{morphism}, it will be said that $\Phi_0$ {\em generates\/} $\Phi$.  Observe that morphisms consisting of local homeomorphisms are precisely those generated by \'etal\'e morphisms (by taking the appropriate target space).

The {\em composition\/} of two consecutive morphisms,
$$
\begin{CD} 
\HH @>{\Phi}>> \HH' @>{\Psi}>> \HH''\;,
\end{CD}
$$
is the morphism $\Psi\circ\Phi:\HH\to\HH''$ generated by all composites of maps
in $\Phi$ with maps in $\Psi$ (wherever defined). With this operation, the morphisms
of pseudogroups form a category $\PsGr$, whose isomorphisms are the morphisms
generated by equivalences of pseudogroups. For a pseudogroup \HH\ acting on a
space $T$, the {\em identity morphism\/} $\id_{\HH}$ of $\PsGr$ at \HH\ is the
morphism generated by $\id_T$; observe that $\HH\subset\id_{\HH}$. 

The {\em restriction\/} of a morphism $\Phi:\HH\to\HH'$ to a subspace $T_0\subset T$ is the morphism $\Phi|_{T_0}:\HH|_{T_0}\to\HH'$ consisting of all maps of open subsets of $T_0$ to $T'$ that can be locally extended to maps in $\Phi$. If $T_0$ is open in $T$, then $\Phi|_{T_0}$ consists of all maps in $\Phi$ whose domain is contained in $T_0$. The inclusion map $T_0\hookrightarrow
T$ generates a morphism $\HH|_{T_0}\to\HH$, whose composition with $\Phi$ is $\Phi|_{T_0}$. 

Suppose that $T=\bigcup_iT_i$. Given a morphism $\Phi_i:\HH|_{T_i}\to\HH'$ for each $i$, let $\Phi$ be the family of continuous maps $\phi:U\to T'$, where $U$ is an open subset of $T$, such that $\phi|_{U\cap T_i}\in\Phi_i$ for all $i$. If $\Phi$ is a morphism and the maps in each $\Phi_i$ can be locally extended to maps in $\Phi$, then $\Phi|_{T_i}=\Phi_i$ for all $i$, and $\Phi$ is called the {\em combination\/} of the morphisms $\Phi_i$.  When every $T_i$ is open in $T$, then the combination of the morphisms $\Phi_i$ is defined just when $\Phi_i|_{T_i\cap T_j}=\Phi_j|_{T_i\cap T_j}$ for all $i$ and $j$.

The {\em image\/} of a morphism $\Phi:\HH\to\HH'$ is the $\HH'$-saturated set $\im\Phi=\bigcup_{\phi\in\Phi}\im\phi$. If \HH\ acts on a space $T$ and $T_0\subset T$, define the {\em direct image\/} 
$$
\Phi(T_0)=\im\Phi|_{T_0}=\bigcup_{\phi\in\Phi}\phi(T_0\cap\dom\phi)\;,
$$
which is $\HH'$-saturated; thus $\Phi_{\text{\rm orb}}(\OO)=\Phi(\OO)$ for any orbit \OO\ of \HH.  It is said that $\Phi$ is {\em constant\/} if $\im\Phi$ is one orbit ($\Phi_{\text{\rm orb}}$ is constant).  If $\HH'$ acts on a space $T'$ and $T'_0\subset T'$, define the {\em inverse image\/} 
$$
\Phi^{-1}(T'_0)=\bigcup_{\phi\in\Phi}\phi^{-1}(T'_0\cap\im\phi)\;,
$$
which is \HH-saturated. If $\im\Phi\subset T'_0$, then the restrictions $\phi:\dom\phi\to T'_0$, for $\phi\in\Phi$, form a morphism that may be denoted by $\Phi:\HH\to\HH'|_{T'_0}$ and called the {\em restriction\/} of $\Phi$ too.

The {\em product\/} of two pseudogroups, $\HH_1$ and $\HH_2$ acting on spaces $T_1$ and $T_2$, respectively, is the pseudogroup $\HH_1\times\HH_2$ acting on $T_1\times T_2$ generated by the products of maps in $\HH_1$ and $\HH_2$. The {\em product\/} of two morphisms $\Phi_i:\HH_i\to\HH'_i$, $i=1,2$, is the morphism $\Phi_1\times\Phi_2:\HH_1\times\HH_2\to\HH'_1\times\HH'_2$ generated by the products of maps in $\Phi_1$ and $\Phi_2$. The {\em pair\/} of two morphisms  $\Phi_i:\HH\to\HH'_i$, $i=1,2$, is the morphism $(\Phi_1,\Phi_2):\HH\to\HH'_1\times\HH'_2$ generated by the pairs $(\phi_1,\phi_2)$, where $\phi_1\in\Phi_1$ and $\phi_2\in\Phi_2$ have the same domain.

\sec{generalization}{Extending concepts to pseudogroups}

Let $\Top$ denote the category of continuous maps between topological spaces. There is a canonical injective covariant functor $\Top\to\PsGr$ which assigns the pseudogroup generated by $\id_T$ to each space $T$, and assigns the morphism generated by $f$ to each continuous map $f$. We will consider $\Top$ as a subcategory of $\PsGr$ in this way; i.e., topological spaces and continuous maps will be also considered as pseudogroups and morphisms.

Let $X$ be a space and \HH\ a pseudogroup acting on another space $T$. Any continuous map $X\to T$ generates a morphism $X\to\HH$. Nevertheless, a continuous map $f:T\to X$ generates a morphism $\HH\to X$ if and only if $f$ is constant on the \HH-orbits; in this case, the notation $f:\HH\to T'$ will be used for such a morphism, which consists of the restrictions of $f$ to all open subsets of $T$. 

Some topological concepts can be extended to pseudogroups by using orbits and morphisms instead of points and continuous maps. For instance, a pseudogroup \HH\ is said to be {\em path connected\/} if, for any pair of orbits \OO\ and $\OO'$ of \HH, there is a morphism $\Phi:I\to\HH$ with $\Phi(0)=\OO$ and $\Phi(1)=\OO'$, where $I=[0,1]$ is considered as a pseudogroup in the above sense. If \HH\ is path connected, then its space of orbits is path connected. A {\em homotopy\/} between two morphisms $\Phi_0,\Phi_1:\HH\to\HH'$ is a morphism $\Psi:\HH\times I\to\HH'$ such that $\Phi_0$ and $\Phi_1$ can be identified to the restrictions of $\Psi$ to $T\times\{0\}\equiv T$ and $T\times\{1\}\equiv T$. Since the restriction of $\HH\times I$ to each slice $T\times\{t\}\equiv T$ can be identified with \HH, a homotopy $\Psi:\HH\times I\to\HH'$ can be considered as the family of its restrictions $\Psi_t=\Psi|_{T\times\{t\}}:\HH\to\HH'$. A morphism $\Phi:\HH\to\HH'$ is called a {\em homotopy equivalence\/} if there is a morphism $\Phi':\HH'\to\HH$ such that $\Phi'\circ\Phi$ and $\Phi\circ\Phi'$ are homotopic to the identity morphisms $\id_{\HH}$ and $\id_{\HH'}$, respectively. If $\id_{\HH}$ is homotopic to a constant morphism, then \HH\ is said to be {\em contractible\/}. 

Other typical topological concepts that can be obviously generalized to pseudogroups in this way are {\em coverings\/}, {\em fiber bundles\/}, the most usual types of {\em $($co$)$ho\-mo\-lo\-gies\/}, {\em homotopy\/} groups, {\em LS category\/}, {\em $K$-theory\/} and {\em cobordism\/}. In particular, {\em coverings\/} and {\em fundamental groups\/} of pseudogroups were defined in \cite{Haefliger:85Pseudogroups} and \cite{vanEst:84}.

The differential calculus can be also extended to \cinf\ pseudogroups by using morphisms. For instance, we can define their {\em de~Rham cohomology\/} (see \refs{cohoms of psgrs}) and study {\em transversality\/} of \cinf\ morphisms.

Finally, morphisms also allow the generalization of many typical concepts of differential geometry to pseudogroups. For example, for pseudogroups of local isometries of Riemannian manifolds (called {\em Riemannian pseudogroups\/}), the {\em geodesics\/} and {\em geodesic segments\/} can be defined as morphisms, and, with more generality, we can consider {\em harmonic morphisms\/} between Riemannian pseudogroups. The {\em length\/} of a geodesic segment of a Riemannian pseudogroup makes sense in this way. Moreover {\em geodesic completeness\/} can be considered for Riemannian pseudogroups.

\sec{foliated sp}{Holonomy pseudogroups of foliated spaces}

With the terminology of \cite{Haefliger:88Leaf} and  \cite{Haefliger:02compactly}, a {\em foliated
structure\/} \FF\ of {\em dimension\/} $n\in\N$ on a space $X$ can be described by
a {\em foliated cocycle\/}, which is a collection $\{U_i,p_i\}$, where $\{U_i\}$ is
an open cover of $X$ and each $p_i$ is a topological submersion of $U_i$ onto some
space $T_i$ whose fibers are connected open subsets of
$\R^n$ and such that the following {\em compatibility condition\/} is satisfied: for every
$x\in U_i\cap U_j$, there is an open neighborhood $U_{i,j}^x$ of $x$ in $U_i\cap
U_j$ and a homeomorphism $h_{i,j}^x:p_i(U_{i,j}^x)\to p_j(U_{i,j}^x)$ such that
$p_j=h_{i,j}^x\circ p_i$ on $U_{i,j}^x$; the notation $T_{i,j}^x=p_i(U_{i,j}^x)$ will be
used for the sake of simplicity.  Two foliated cocycles determine the same foliated structure when their union is a foliated cocycle. Thus \FF\ is given as an equivalence class of foliated cocycles---the logical problems of this definition can be easily avoided, either by using the theory of universes of Grothendieck's school, or by using only foliated cocycles with values in a given set of spaces (see e.g. \cite{MacLane69:universe}). The space $X$ endowed with \FF\ is called a {\em foliated space\/}. 
 
Many interesting examples of foliated spaces are given e.g. in
\cite{CandelConlonI}. Foliated structures and foliated spaces with {\em boundary\/} or {\em corners\/}
can be defined similarly. Usually, foliated spaces are supposed to be Polish (completely metrizable and separable), and local compactness is also assumed often; recall that a space is locally compact and Polish if and only if it is locally compact and second countable \cite[Theorem~5.3]{Kechris:95}. But these conditions are not needed for the purposes of this section and the next one. Also, it would be enough that the fibers of the submersions $p_i$ be arbitrary connected locally path connected spaces. When a foliated space is a manifold, then the common term {\em foliation\/} is used for its foliated structure.

Let us recall the generalization to foliated spaces of some
common terminology for foliations. The domains of the submersions of the foliated cocycles are called {\em simple open sets\/}. An open covering of $X$ is called {\em simple\/} when it consists of simple open sets. The simple open sets form a base of the topology of $X$. The {\em local quotient\/} of a simple open set is its quotient space whose points are the plaques, which can be identified with the target space of the corresponding submersion. The fibers of those submersions of any foliated cocycle are called {\em plaques\/}. The plaques of all foliated cocycles form a base of a finer topology on $X$, called the {\em leaf topology\/}, whose connected components are called {\em leaves}. The leaf through each point $x$ is usually denoted by $L_x$. The quotient space of leaves will be denoted by $X/\FF$. The images of the local sections of the submersions of foliated cocycles are called {\em local transversals\/}.

The condition on each $p_i$ to be a topological submersion means that, for each
$x\in U_i$, there is a homeomorphism $\theta_i^x$ of an open neighborhood $U_i^x$ of
$x$ to $p_i(U_i^x)\times B$, for some fixed ball $B$ of $\R^n$, so that $p_i$
corresponds to the first factor projection by $\theta_i^x$. Then
$(U_i^x,\theta_i^x)$ is called a {\em foliated chart\/}. A collection of foliated charts whose domains cover $X$ is called a {\em foliated atlas\/}.

The foliated structure \FF\ can be also identified with its {\em canonical foliated cocycle\/}, which  consists of the canonical projections of all simple open subsets of $X$ onto their local quotients. Other well known descriptions of \FF\ can be given by using foliated atlases, or by using the partition of $X$ into the leaves, satisfying appropriate conditions.

A simple open set $V$ is said to be {\em uniform\/} in another simple open set $U$
when $V\subset U$ and every plaque of $V$ meets only one plaque in $U$. A simple
open covering $\{U_i\}$ of $X$ is called {\em regular\/} if, whenever $U_i$ meets $U_j$, there is some simple open set of $X$ where both $U_i$ and $U_j$ are uniform;
in particular, each plaque of $U_i$ meets at most one plaque of $U_j$, and thus we
can take $U_{i,j}^x=U_i\cap U_j$ in the above compatibility condition. If $X$ is
Polish and locally compact, then it has arbitrarily fine locally finite regular simple open coverings (use the arguments of \cite{HectorHirsch:81}, \cite{Godbillon} and \cite{AlvCandel:ggol}).

For a foliated cocycle $\{U_i,p_i\}$ of $\FF$ with $p_i:U_i\to T_i$, the homeomorphisms $h_{i,j}^x$,
given by its compatibility condition, generate a pseudogroup \HH\ acting on the
topological sum $T=\bigsqcup_i T_i$, and the maps $p_i$ generate
a morphism $\PP:X\to\HH$; these \HH\ and \PP\ are said to be {\em
induced\/} by $\{U_i,p_i\}$. Let $\{U'_a,p'_a\}$ be another
foliated cocycle of \FF\ with
$p'_a:U'_a\to T'_a$, which induces a pseudogroup $\HH'$ acting on 
$T'=\bigsqcup_a T'_a$ and a morphism $\PP':X\to\HH'$. Then the foliated cocycle
$\{U_i,p_i\}\cup\{U'_a,p'_a\}$ induces a
pseudogroup $\HH''$ acting on $T''=T\sqcup T'$ whose orbits meet $T$ and $T'$, and
so that $\HH''|_T=\HH$ and $\HH''|_{T'}=\HH$. Therefore $\HH''$
defines a canonical equivalence $\Phi_0:\HH\to\HH'$, which generates a canonical isomorphism $\Phi:\HH\to\HH'$ satisfying
$\Phi\circ\PP=\PP'$. The equivalence class of the pseudogroup induced by any
foliated cocycle of \FF\ is usually called the {\em holonomy pseudogroup\/}. But the canonical foliated cocycle induces a canonical representative of this class, which will be also called {\em holonomy pseudogroup\/} and denoted by $\Hol(\FF)$. By identifying the local quotient of each simple open set with a local transversal, we see that the holonomy pseudogroup can be given by sliding local transversals along the leaves; thus it represents the ``transverse dynamics'' of the foliated space.

Let us introduce the following terminology. Fix $m,n\in\N$, and spaces $T$ and $T'$. Consider the foliated structures on $T\times\R^m$ and $T'\times\R^n$ with leaves $\{x\}\times\R^m$ and $\{x'\}\times\R^n$ for $x\in T$ and $x'\in T'$, as well as their restrictions to open subsets $U\subset T\times\R^m$ and $V\subset T'\times\R^n$. Then a foliated map $f:U\to V$ is said to of {\em class \coinf\/} if, for each $y\in Y$, the mapping $z\mapsto\pr_2\circ f(y,z)$ is \cinf\ with partial derivatives of arbitrary order depending continuously on $y$, where $\pr_2:T'\times\R^n\to\R^n$ is the second factor projection. A {\em \coinf\ structure\/} on \FF\ is a maximal foliated atlas $\{U_i,\theta_i\}$ so that each composite $\theta_j\circ\theta_i^{-1}$ is \coinf. When $X$ is endowed with a \coinf\ structure, it is called a {\em \coinf\ foliated space\/}.

The foliated structure \FF\ is said to be of {\em class \cinfo\/} when $\Hol(\FF)$ is \cinf; in particular, \FF\ has to be a foliation. 

If $X$ is a manifold, the topological submersions of any foliated cocycle of \FF\ have values in manifolds. If $X$ is a \cinf\ manifold and there is a foliated cocycle of \FF\ consisting of \cinf\ submersions, then \FF\ is called a  {\em\cinf foliation\/}.

\sec{foliated map}{Holonomy morphisms of foliated maps}

Let $X$ and $Y$ be foliated spaces with respective foliated structures \FF\ and
\GG.  A {\em foliated map\/}
$f:(X,\FF)\to(Y,\GG)$, or simply $f:\FF\to\GG$, is a map $f:X\to Y$ which
maps leaves of \FF\ to leaves of \GG; thus $f$ induces a map $\bar f:X/\FF\to Y/\GG$, which is continuous if $f$ is continuous. The identity map $\id_X$, considered as a foliated map $\FF\to\FF$, will be denoted by $\id_{\FF}$. The set of continuous foliated maps
$\FF\to\GG$ will be denoted by $C(X,\FF;Y,\GG)$, or simply
$C(\FF,\GG)$. Continuous foliated maps between foliated spaces form a category with the
operation of composition, which will be denoted by $\Fol$. The mapping $f\mapsto\bar f$ defines a functor $\Fol\to\Top$.

When \FF\ and \GG\ are of class \coinf, a foliated map $f:\FF\to\GG$ is said to be of {\em class \coinf\/} when, for all foliated charts $(U,\theta)$ and $(U',\theta')$ of the \coinf\ structures of \FF\ and \GG\ such that $f(U)\subset U'$, the composite $\theta'\circ f\circ\theta^{-1}$ is \coinf. The set of \coinf\ foliated maps $\FF\to\GG$ will be denoted by $\coinf(\FF,\GG)$.

When \FF\ and \GG\ are \cinf\ foliations, the set of \cinf\ foliated maps $\FF\to\GG$ will be denoted by $\cinf(\FF,\GG)$.

Let $\{U_i,p_i\}$ and $\{V_a,p'_a\}$ be foliated cocycles of \FF\ and \GG\ with
$p_i:U_i\to T_i$ and $p'_a:V_a\to T'_a$. For $x\in U_i\cap U_j$ and $y\in V_a\cap V_b$, the compatibility condition is satisfied with open sets $U_{i,j}^x$ and $V_{a,b}^y$, and homeomorphisms $h_{i,j}^x:T_{i,j}^x\to T_{j,i}^x$ and $h_{a,b}^{\prime y}:T_{a,b}^y\to T_{b,a}^{\prime y}$.
Let $\HH$ and $\HH'$ be the
pseudogroups induced by $\{U_i,p_i\}$ and $\{V_a,p'_a\}$, acting on $T=\bigsqcup_iT_i$ and $T'=\bigsqcup_aT'_a$, and let
$\PP:X\to\HH$ and $\PP':Y\to\HH'$ be the corresponding morphisms.  

For any fixed $f\in C(\FF,\GG)$, we can choose the open sets $U_{i,j}^x$ and $V_{a,b}^y$
such that $f$ maps each fiber of $p_i$ in $U_{i,j}^x$ to a fiber of $p'_{a_i}$ on $V_{a_i,a_j}^{f(x)}$ for some mapping $i\mapsto a_i$. So there are continuous maps $\phi_{i,j}^x:T_{i,j}^x\to T^{\prime f(x)}_{a_i,a_j}$ satisfying 
\begin{equation}\label{e:phi i,j x}
\phi_{i,j}^x\circ p_i=p'_{a_i}\circ f 
\end{equation}
on $U_{i,j}^x$. Let $\Phi_0$ be the family of such maps $\phi_{i,j}^x$.

\lem{Phi}
$\Phi_0$ generates a morphism $\Phi:\HH\to\HH'$ such that $\PP'\circ f=\Phi\circ\PP$.
\elem

\prf
This $\Phi_0$ obviously satisfies  hypothesis~(ii') of \refl{morphism}, and the hypothesis~(iii') is given by the commutativity of the diagram
$$
\begin{CD}
T_{i,j}^x @>{\phi_{i,j}^x}>> T^{\prime f(x)}_{a_i,a_j}\\
@V{h_{i,j}^x}VV@VV{h^{\prime f(x)}_{a_i,a_j}}V\\
T_{j,i}^x @>{\phi_{j,i}^x}>> T^{\prime f(x)}_{a_j,a_i}
\end{CD}
$$
for $x\in U_i\cap U_j$, which holds because
$$
\phi_{j,i}^x\circ h_{i,j}^x\circ p_i=\phi_{j,i}^x\circ p_j=p'_{a_j}\circ f
=h^{\prime f(x)}_{a_i,a_j}\circ p'_{a_i}\circ f=h^{\prime f(x)}_{a_i,a_j}\circ\phi_{i,j}^x\circ p_i
$$
on $U_{i,j}^x$ by~\eqref{e:phi i,j x}. So $\Phi_0$ generates a morphism $\Phi:\HH\to\HH'$, and the equality $\PP'\circ f=\Phi\circ\PP$ also follows from~\eqref{e:phi i,j x}.
\eprf

It will be said that $\Phi$ is {\em induced\/} by $f$, or that $f$ is a {\em lift\/} of $\Phi$. In particular, for the canonical representatives, we get a morphism $\Hol(f):\Hol(\FF)\to\Hol(\GG)$, which is called the {\em holonomy morphism\/} induced by $f$; we may also say that $\Phi$ is a {\em representative\/} of $\Hol(f)$. In this way, we get a covariant functor $\Hol:\Fol\to\PsGr$, which is called the {\em holonomy functor\/}. 

When \FF\ and \GG\ are \cinfo\ foliations, a foliated map $f:\FF\to\GG$ is said to be of {\em class \cinfo\/} if $\Hol(f):\Hol(\FF)\to\Hol(\GG)$ is a \cinf\ morphism. The set of \cinfo\ foliated maps $\FF\to\GG$ will be denoted by $\cinfo(\FF,\GG)$.

Any topological space $X$ can be considered as a foliated space whose leaves are its
points; the notation $X_{\text{\rm pt}}$ will be used in this case. Moreover any map
$X\to Y$, between topological spaces, is a foliated map $X_{\text{\rm pt}}\to
Y_{\text{\rm pt}}$. This defines an injective functor $\Top\to\Fol$ whose
composition with the holonomy functor is the canonical injective functor
$\Top\to\PsGr$; this functor $\Top\to\Fol$ will be considered as an inclusion. 

Let $X$ be a space and \FF\ a foliated structure on another space $Y$. Any continuous map $X\to Y$ is a foliated map $X_{\text{\rm pt}}\to\FF$. However, a continuous map $Y\to X$ is a foliated map $\FF\to X$ if and only if if is constant on the leaves of \FF.

Any manifold $M$ (possibly with boundary) can be also
considered as a foliated space (possibly with boundary) whose leaves are
its connected components; this foliated space will be also denoted by $M$. Furthermore any map
between manifolds is a foliated map in this sense. This defines an injective functor
$\Man\to\Fol$, which will be also considered as an inclusion, where $\Man$ denotes
the category of continuous maps between manifolds. The composition of $\Man\to\Fol$ with the holonomy functor is the functor of projection to the discrete space of connected components.

Let $Y$ be a foliated space with foliated structure \GG, and let \FF\ be a foliated structure on some subspace $X\subset Y$. If each leaf of \FF\ is a submanifold of some leaf of \GG, then \FF\ is called a {\em subfoliated structure\/} of \GG; in this case, the inclusion map $X\hookrightarrow Y$ is a foliated map of \FF\ to \GG\ denoted by $\FF\hookrightarrow\GG$.

For foliated spaces $X_1$ and $X_2$ with respective foliated structures $\FF_1$ and $\FF_2$, the {\em product\/} $\FF_1\times\FF_2$ is the foliated structure on $X_1\times X_2$
whose leaves are the products of leaves of $\FF_1$ and leaves of $\FF_2$. If $(U^1_i,p^1_i)$ and $(U^2_j,p^2_j)$ are foliated cocycles of $\FF_1$ and $\FF_2$, respectively, then $(U^1_i\times U^2_j,p^1_i\times p^2_j)$ is a foliated cocycle of $\FF_1\times\FF_2$, called the {\em product\/} of $(U^1_i,p^1_i)$ and $(U^2_j,p^2_j)$. The factor projections $\pr_k:X_1\times X_2\to X_k$ are foliated maps $\pr_k:\FF_1\times\FF_2\to\FF_k$, $k=1,2$. Then 
$$
(\Hol(\pr_1),\Hol(\pr_2)):\Hol(\FF_1\times\FF_2)\to\Hol(\FF_1)\times\Hol(\FF_2)
$$ 
is an isomorphism. Indeed, $\Hol(\FF_1)\times\Hol(\FF_2)$ is induced by the product of the canonical  foliated cocycles of $\FF_1$ and $\FF_2$, and $(\Hol(\pr_1),\Hol(\pr_2))$ is the canonical isomorphism. Observe that $\Hol(\FF_1)\times\Hol(\FF_2)$ is the restriction of $\Hol(\FF_1\times\FF_2)$ to an open set that cuts every orbit, and the corresponding inclusion map generates $(\Hol(\pr_1),\Hol(\pr_2))^{-1}$.

The product $f_1\times f_2$ of  foliated maps $f_i:\FF_i\to\GG_i$, $i=1,2$, is a
foliated map $\FF_1\times\FF_2\to\GG_1\times\GG_2$. When $f_1$ and $f_2$ are
continuous, the morphism $\Hol(f_1)\times\Hol(f_2)$ corresponds to $\Hol(f_1\times f_2)$ by the canonical isomorphisms. The pair $(f_1,f_2)$ of  foliated maps $f_i:\FF\to\GG_i$, $i=1,2$, is a
foliated map $\FF\to\GG_1\times\GG_2$. When $f_1$ and $f_2$ are
continuous, the morphism $(\Hol(f_1),\Hol(f_2))$ corresponds to $\Hol(f_1,f_2)$ by the canonical equivalence.

If $M$ is a connected manifold (possibly with boundary or corners) considered as a foliation with one leaf, then $\Hol(M)$ has only one orbit, and it is thus isomorphic to a singleton space $\{*\}$. More precisely, the holonomy morphism of the map $M\to\{*\}$ is an isomorphism $\Hol(M)\to\Hol(\{*\})\equiv\{*\}$. Let \FF\ be a foliated structure, and let $\pr_1:\FF\times M\to\FF$ and $\pr_2:\FF\times M\to M$ be the factor projections. Then the equivalence
$$
(\Hol(\pr_1),\Hol(\pr_2)):\Hol(\FF\times M)\to\Hol(\FF)\times\Hol(M)
$$
corresponds to $\Hol(\pr_1):\Hol(\FF\times M)\to\Hol(\FF)$ by the isomorphism 
$$
\Hol(\FF)\times\Hol(M)\cong\Hol(\FF)\times\{*\}\equiv\Hol(\FF)\;.
$$
Given any point $y\in M$, let $\iota_y:\FF\to\FF\times M$ be the foliated map defined by $\iota_y(x)=(x,y)$. Since $\pr_1\circ\iota_y=\id_{\FF}$, we get $\Hol(\iota_y)=\Hol(\pr_1)^{-1}$, which is independent of $y$.

\sec{homotop}{Integrable and foliated homotopies}

Let $X$ and $Y$ be foliated spaces with foliated structures \FF\ and \GG. An {\em integrable homotopy\/}
between continuous foliated maps $f,g:\FF\to\GG$ is a homotopy $H:X\times I\to Y$ between $f$ and $g$ which is a foliated map $\FF\times I\to\GG$; {\em i.e.}, each homotopy curve $t\mapsto H(x,t)$ lies in a leaf of $\GG$. If there is an integrable homotopy
between $f$ and $g$, then these maps are said to be {\em integrably homotopic\/}. A continuous
foliated map $f:\FF\to\GG$ is called a {\em integrable homotopy equivalence\/} if
there is a continuous foliated map $g:\GG\to\FF$ such that $g\circ f$ and
$f\circ g$ are integrably homotopic to $\id_{\FF}$ and $\id_{\GG}$. 

\prop{integrable homotopy}
Integrably homotopic foliated maps define the same holonomy morphism.
\eprop

\prf
Let $H:\FF\times I\to\GG$ be an integrable homotopy  between foliated maps $f,g:\FF\to\GG$, and let $\iota_i:\FF\to\FF\times I$, $i=0,1$, be the foliated maps defined by $\iota_i(x)=(x,i)$. Since $I$ is connected, we have $\Hol(\iota_0)=\Hol(\iota_1)$ (\refs{foliated map}). Therefore
\begin{multline*}
\Hol(f)=\Hol(H\circ\iota_0)=\Hol(H)\circ\Hol(\iota_0)\\
=\Hol(H)\circ\Hol(\iota_1)=\Hol(H\circ\iota_1)=\Hol(g)\;.\qed
\end{multline*}
\renewcommand{\qed}{}
\eprf

\cor{integrable homotopy}
With the above notation, if $f:\FF\to\GG$ is an integrable homotopy equivalence, then $\Hol(f)$ is an isomorphism.
\ecor

A {\em foliated homotopy\/} between continuous foliated maps $f,g:\FF\to\GG$ is a homotopy $H$ between $f$ and $g$ which is a foliated map $\FF\times I_{\text{\rm pt}}\to\GG$; i.e., the homotopy $H$ consists of foliated maps $H_t=H(\cdot,t):\FF\to\GG$. Any integrable homotopy is a foliated homotopy. Two foliated maps are said to be {\em foliatedly homotopic\/} if there is a foliated homotopy between them. A continuous foliated map $f:\FF\to\GG$ is called a {\em foliated homotopy equivalence\/} if there is a continuous foliated map $g:\GG\to\FF$ such that $g\circ f$ and $f\circ g$ are foliatedly homotopic to $\id_{\FF}$ and $\id_{\GG}$, respectively. Since $\Hol(\FF\times I_{\text{\rm pt}})\cong\Hol(\FF)\times I$ canonically, if $H:\FF\times I_{\text{\rm pt}}\to\GG$ if a foliated homotopy  between $f$ and $g$, then $\Hol(H):\Hol(\FF\times I_{\text{\rm pt}})\to\Hol(\GG)$ canonically defines a homotopy between $\Hol(f)$ and $\Hol(g)$. Therefore we get the following.

\prop{foliated homotopy}
If $f$ is a foliated homotopy equivalence, then $\Hol(f)$ is a homotopy equivalence.
\eprop

A continuous foliated map $f:\FF\to\GG$ is called a {\em proper integrable homotopy equivalence\/} if it is proper and there is a proper continuous foliated map $g:\GG\to\FF$ such that there are proper integrable homotopies between $g\circ f$ and $\id_{\FF}$, and between $f\circ g$ and $\id_{\GG}$. A {\em proper foliated homotopy equivalence\/} can be defined similarly by using proper foliated homotopies instead of proper integrable homotopies.

\sec{complete}{Complete pseudogroups and complete morphisms}

 For any map $h:T\to T'$ between topological spaces, let $\germ{h}{x}$ denote its germ at any
$x\in T$. If \HH\ is a family of maps of open subsets
of $T$ to open subsets of $T'$, let $\germs{\HH}$ denote the space of all germs of maps in \HH\ with the \'etal\'e topology. If $T$ is a manifold, then $\germs{\HH}$ is a manifold of the same dimension. If \HH\ is a pseudogroup acting on $T$, then $\germs{\HH}$ is a topological groupoid with the operation induced by composition. 

Recall from \cite{Haefliger:88Leaf} that a pseudogroup \HH\ acting on a space $T$ is said to be {\em complete\/} if, for
all $x,y\in T$, there are open neighborhoods $U$ and $V$ of $x$ and $y$ such that,
for any $h\in\HH$ and any
$z\in U\cap \dom h$ with $h(z)\in V$, there is some $\tilde h\in\HH$ so that $U\subset\dom\tilde
h$ and $\gamma(\tilde h,z)=\gamma(h,z)$; in this case, $(U,V)$ is called a {\em completeness pair\/}. 

\defi{complete morphism}
Let \HH\ and $\HH'$ be pseudogroups acting on topological spaces $T$ and $T'$, respectively. A morphism $\Phi:\HH\to\HH'$ is said to be {\em complete\/} when, given
$\phi,\psi\in\Phi$, $x\in\dom\phi$ and $y\in\dom\psi$, there are open neighborhoods $U$ and $V$ of $x$ and $y$ in $\dom\phi$ and $\dom\psi$, respectively, such that, for all $h\in\HH$ and every $z\in U\cap \dom h$ with $h(z)\in V$, there is some $\tilde
h\in\HH$ and some $h'\in\HH'$ so that
$U\subset\dom\tilde h$, $\tilde h(U)\subset\dom\psi$, $\gamma(\tilde h,z)=\gamma(h,z)$, $\phi(U)\subset\dom h'$, and
$h'\circ\phi=\psi\circ\tilde h$ on $U$. In this case, $(\phi,U;\psi,V)$ is called a {\em completeness quadruple\/}.
\edefi

\rem{complete morphism}
In \refd{complete morphism}, if \HH\ is complete, then $(U,V)$ can be chosen to be a completeness pair of \HH. Thus, in this case, to prove that $(\phi,U;\psi,V)$ is a completeness quadruple of $\Phi$, it is
enough to take elements $h\in\HH$ with $\dom h=U$.
\erem

Observe that a pseudogroup \HH\ is complete if and only if the identity morphism $\id_{\HH}$ is complete.

\sec{singular}{Singular foliated spaces}

By adapting a definition by \cite{Stefan:74}, a {\em singular foliated structure\/} \FF\ on a space $X$ can be described by a {\em singular foliated cocycle\/}, which is a collection $\{U_i,p_i\}$, where $\{U_i\}$ is an open cover of $X$ and each $p_i$ is a topological submersion of $U_i$ onto some
space $T_i$ whose fibers are connected open subsets of $\R^{n_i}$, with $n_i\in\N$ depending on $i$, and such that the following properties are satisfied:
\begin{itemize}

\item {\em Compatibility condition\/}: any $x\in X$ is contained in some $U_i$ such that, if
$x\in U_j$, then there is some open neighborhood $U_{j,i}^x$ of $x$ in $U_j\cap U_i$ and a topological submersion $h_{j,i}^x:p_j(U_{j,i}^x)\to p_i(U_{j,i}^x)$ such that $p_i=h_{j,i}^x\circ p_j$ on $U_{j,i}^x$; in this case, the fiber $p_i^{-1}(p_i(x))$ is called a {\em plaque\/} of $\{U_i,p_i\}$ (or of $p_i$ or $U_i$).

\item Some fiber of any $p_i$ is a plaque.

\end{itemize}
Two singular foliated cocycles determine the same singular foliated structure when their union is a singular foliated cocycle. The space $X$ endowed with \FF\ is called a {\em singular foliated space\/}. 

The plaques of all singular foliated cocycles defining \FF\ form a base of a topology, called the {\em leaf topology\/}. The connected components of $X$ with the leaf topology are called the {\em leaves\/} of \FF. The quotient space of leaves will be denoted by $X/\FF$. 

Given a foliated cocycle $\{U_i,p_i\}$ and a leaf $L$ of \FF, each intersection $L\cap U_i$ is a union of fibers of $p_i$; such fibers are open in $L$ just when they are plaques.

For $x$ is in a plaque of some $p_i$, then the image $\Sigma$ of any local section of $p_i$ containing $x$ is called a {\em local transversal\/} of \FF\ through $x$. There is a unique singular foliated structure $\FF_\Sigma$ on $\Sigma$ such that, for $V=p_i^{-1}(p_i(\Sigma))$, the submersion $p_i$ defines a foliated map $p_i:\FF|_{V}\to\FF_\Sigma$, which restricts to submersions of the leaves of $\FF|_V$ to the leaves of $\FF_\Sigma$. Observe that $\FF_T$ is the foliated structure by points in the regular case. Notice also that, in the singular case, an open subset of a local transversal may not be a local transversal through any point.

Other concepts and notations of foliated structures can be directly generalized to the singular case: {\em saturations\/}, {\em restrictions\/}, {\em products\/} of singular foliated structures,  {\em foliated maps\/}, {\em integrable homotopies\/}, {\em foliated homotopies\/}, etc. If $X$ is a manifold, then \FF\ is called a {\em singular foliation\/}; in this case, the submersions of a foliated cocycle have values in manifolds.  If $X$ is a \cinf\ manifold and there is some foliated cocycle consisting of \cinf\ submersions, then \FF\ is called a {\em \cinf\ singular foliation\/}. 

Singular foliated structures whose leaves are manifolds with boundary, or arbitrary connected  locally path connected spaces, can be defined similarly.

Any foliated structure is a singular foliated structure; it may be said that these foliated structures are {\em regular\/} for emphasis. A (\cinf) singular foliated structure is a (\cinf) regular foliated structure if and only if all of its leaves have the same dimension.

By \cite{Ramsay:91}, for any continuous local action of a local Lie group on a separable metric space, the connected components of the orbits are the leaves of a singular foliated structure, which is regular just when the local action is locally free. In the case of \cinf\ local actions of local Lie groups on \cinf\ manifolds, we get \cinf\ singular foliations.

Let \FF\ and \GG\ be singular foliated structures. The set of continuous foliated maps $\FF\to\GG$ will be denoted by $C(\FF,\GG)$. When \FF\ and \GG\ are \cinf\ singular foliations, the set of \cinf\ foliated maps $\FF\to\GG$ will be denoted by $\cinf(\FF,\GG)$.

The following notation will be used.

\begin{notation}
The tangent vector bundle of \cinf\ manifold $M$ will be denoted by $\TT M$, and the tangent space of $M$ at some point $x\in M$ will be denoted by $\TT_xM$. As usual, the Lie algebra of \cinf\ tangent vector fields on $M$ is denoted by $\frakX(M)$. For any \cinf\ map between \cinf\ manifolds, $f:M\to N$, its tangent homomorphism $\TT M\to\TT N$ is denoted by $\TT f$ or $f_*$.
\end{notation}

By the main result of \cite{Stefan:74}, a partition \FF\ of $M$ into \cinf\ immersed connected submanifolds is a \cinf\ singular foliation if and only if the evaluation map $\frakX(\FF)\to\TT_x\FF$ is surjective for each $x\in M$, where $\frakX(\FF)\subset\frakX(M)$ is the Lie subalgebra of \cinf\ vector fields on $M$ tangent to the leaves, and $\TT_x\FF\subset\TT_x M$ is the linear subspace of
vectors tangent to the leaf through $x$.

Let \FF\ and \GG\ be \cinf\ singular foliations on \cinf\ manifolds $M$ and $N$, and let $f\in C(\FF,\GG)$ with \cinf\ restrictions to the leaves of \FF. Let also
$\phi:U\to\R^m$ and $\psi:V\to\R^n$ be charts of $M$ and $N$ such that $f(U)\subset V$. For any positive integer
$r$, let $X_1,\dots,X_r\in\frakX(\FF)$. Since
the restriction of $f$ to the leaves of \FF\ is \cinf, the order $r$ derivative
\begin{equation}\label{e:Dr}
D^r\left(\psi\circ f\circ\phi^{-1}\right)(x)(\phi_\ast X_1(x),\dots,\phi_\ast X_r(x))
\end{equation}
is well defined for all $x\in\phi(U)$. 

\defi{coinf}
The map $f$ will be said to be of {\em class $C^{0,k}$\/} if, for all $r\leq k$ and all possible
charts $\phi$ and $\psi$ as above, the derivatives of the type~\refe{Dr} depend continuously on $x$. The map $f$
will be said to be of {\em class \coinf\/} if it is of class $C^{0,k}$ for all $k$. The set of \coinf~maps
$\FF\to\GG$ will be denoted by $\coinf(\FF,\GG)$.
\edefi

This definition generalizes to the \cinf\ singular foliations the concept of \coinf\ foliated map between \coinf\ foliated spaces (\refs{foliated map}).

\lem{coinf}
Consider foliated maps
$$
\begin{CD}
\FF@>f>>\FF'@>g>>\FF''
\end{CD}
$$
between singular \cinf\ foliations.
Suppose the composite $g\circ f$ is \coinf, and 
$f$ is a \cinf\ surjective submersion such that, for all $Y\in\frakX(\FF')$, there is some
$X\in\frakX(\FF)$ with $f_\ast(X)=Y$. Then $g$ is \coinf.
\elem

\prf
Observe that the hypothesis on $f$ imply
that it is also a surjective submersion as a map from leaves to leaves. So the restriction of $g$ to the
leaves is of class \cinf\ because the restriction of $g\circ f$ to the leaves is of class \cinf.

Let $M$, $M'$ and $M''$ be the ambient manifolds of \FF, $\FF'$ and $\FF''$, and let $\psi:V\to\R^n$ and
$\chi:W\to\R^p$ be charts of $M'$ and $M''$ such that $g(V)\subset W$. Since $f$ is a \cinf\ surjective
submersion, we can assume $V=f(U)$ for some chart $\phi:U\to\R^m$ of $M$. For any positive integer $r$,
given $Y_1,\dots,Y_r\in\frakX(\FF')$, we know the existence of some $X_1,\dots,X_r\in\frakX(\FF)$ such
that $f_\ast(X_i)=Y_i$, $i=1,\dots,r$. Then, if
$y=\psi\circ f\circ \phi^{-1}(x)$ for $x\in\phi(U)$ and $y\in\psi(V)$, we have
\begin{multline}\label{e:coinf}
D^r(\chi\circ g\circ\psi^{-1})(y)(\psi_\ast Y_1(y),\dots,\psi_\ast Y_r(y))\\
=D^r(\chi\circ g\circ f\phi^{-1})(x)(\phi_\ast X_1(x),\dots,\phi_\ast X_r(x))\;.
\end{multline}
But the right hand side of~\refe{coinf} depends continuously on $x\in\phi(U)$ because $g\circ f$ is of class \coinf. Then the left hand side of~\refe{coinf} depends continuously on $y\in\psi(V)$ because $\psi\circ f\circ \phi^{-1}:\phi(U)\to\psi(V)$ is a surjective submersion, and thus a quotient map.
\eprf

\rem{coinf}
Since being \coinf\ is a local property, in \refl{coinf}, it is enough to require that the lifting property of $f$ is
satisfied locally.
\erem

\sec{Riem psgrs}{Riemannian pseudogroups}

A pseudogroup \HH\ of local isometries of an $n$-dimensional Riemannian manifold $T$ will be called a {\em Riemannian pseudogroup\/} on $T$. The corresponding groupoid of germs $\germs{\HH}$ is Hausdorff because a local isometry with connected domain is the identity if it is the identity on some nontrivial open subset ({\em quasi-analyticity\/}). Let $J^1(T)$ denote the topological groupoid of $1$-jets of local diffeomorphisms of $T$, and
$j^1:\germs{\HH}\to J^1(T)$ be defined by mapping each germ to its $1$-jet. $J^1(T)$ is a manifold of
dimension $n^2+2n$ and $j^1$ is a continuous homomorphism. Moreover $j^1$ is injective because germs of local isometries are determined by their $1$-jets. 

For any open subset $U\subset T$, let $\HH_U=\{h\in\HH\ |\ \dom h=U\}$ endowed with
the topology of uniform convergence.

\lem{jetgerm}
If $U$ is a connected relatively compact subset of $T$ and $x\in U$, then the map
$$
U\times\HH_U\lar J^1(T)\;,\quad (x,f)\mapsto j^1(\germ{f}{x})\;,
$$
is an embedding.
\elem

\prf 
This follows because, by the conditions on $U$, any $f\in\HH_U$ is well known to be continuously
determined by its $1$-jet at any fixed point. 
\eprf

From now on in this section, assume that \HH\ is complete.

\begin{Thm}[{Haefliger \cite[Proposition~3.1]{Haefliger:88Leaf}}]\label{t:Haefliger}
With the above notation and conditions, there is a unique Riemannian pseudogroup $\overline{\HH}$ on $T$ such that $\jetgerms{\overline{\HH}}=\overline{j^1(\germs{\HH})}$.
Moreover $\overline{\HH}$ is complete, its orbits are the closures of the \HH-orbits, and
$\overline{\HH}\backslash T$ is Hausdorff.
\end{Thm}

The pseudogroup $\overline{\HH}$ of \reft{Haefliger} is called the {\em closure\/} of \HH, and \HH\ is
said to be {\em closed\/} if $\HH=\overline{\HH}$.

From \reft{Haefliger}, we get $\overline{\HH(x)}=\overline{\HH(y)}$ for all  all $x\in T$ and all $y\in\overline{\HH(x)}$. It follows that any \HH-saturated open subset $U\subset T$ is $\overline{\HH}$-saturated too; indeed, by the above observation, $U$ cuts all \HH-orbits in $\overline{\HH(x)}$ for any $x\in U$.

Let $\HH'$ be another complete Riemannian pseudogroup on a Riemannian manifold $T$. For any morphism $\Phi:\HH\to\HH'$, the orbit map $\Phi_{\text{\rm orb}}:\HH\backslash T\to\HH'\backslash T'$ induces a continuous map 
$\Phi_{\overline{\text{\rm orb}}}:\overline{\HH}\backslash T\to\overline{\HH'}\backslash T'$ called the {\em orbit closure map\/}.

The following is a generalization of the theorem of Myers-Steenrod.

\begin{Thm}[{Salem \cite{Salem:88}}]\label{t:Salem}
With the above notation and conditions, suppose that \HH\ is closed. For each point $x\in T$, there is an open neighborhood $U$ of $x$ and a finite dimensional Lie algebra
$\frakG(U)$ of Killing vector fields over $U$ such that, for each relatively compact
open set $V$ with $\overline{V}\subset U$, the elements of \HH\ with domain $V$ and close enough to the identity are the maps of the form $\exp\xi$ for $\xi\in\frakG(U)$ small enough.  
\end{Thm}

In \reft{Salem}, the notation $\exp t\xi$ is used for the local uniparametric group of diffeomorphisms
defined by a \cinf\ vector field $\xi$.

The elements of the Lie algebra $\frakG(U)$ of \reft{Salem} are the
sections on $U$ of a locally constant sheaf $\frakG $ of Lie algebras of germs of vector
fields over $T$, upon which \HH\ acts by automorphisms \cite[Section~3.4]{Haefliger:88Leaf}. Such a 
$\frakG$ is called the {\em sheaf of infinitesimal transformations of \HH\/}, and its
typical stack, denoted by \frakg, is called the {\em structural Lie algebra\/} of \HH. When \HH\ is not closed, the same terminology is used for the sheaf and Lie algebra associated to $\overline{\HH}$.

\sec{calH0}{The pseudogroup generated by the elements close to identity maps}

Let \HH\ be a closed complete Riemannian pseudogroup on a Riemannian manifold $T$, and let
$\frakG $ denote its sheaf of infinitesimal transformations. Then let $\HH_0$ denote the
Riemannian pseudogroup on $T$ generated by the maps $\exp\xi$, where $\xi$ is any local section of
$\frakG $. Using combination of maps and the formula $\exp\xi=(\exp\xi/N)^N$,
$N\in\N$, it follows that $\HH_0$ is also generated by the maps $\exp\xi$ with
$\xi\in\frakG(U)$ small enough, and $U$ small enough as well. So $\HH_0\subset\HH$; in fact,
according to \reft{Salem}, $\HH_0$ is also generated by the elements of \HH\ that are close enough to
the identity map in their domains. The main goal of this section is to prove the following.

\th{calH0}
With the above notation and conditions, the Riemannian
pseudogroup $\HH_0$ is complete and closed, and its orbits are the connected components
of the orbits of \HH.
\eth

The following notation will be used.

\begin{notation}
Let $G$ be a groupoid. For units $x$ and $y$ of $G$, let $G_x$ (respectively, $G^x$) denote the subset of elements of $G$ with source (respectively, target) $x$, and let
$G_x^y=G_x\cap G^y$. If $X,Y$ are subsets of the space of units of $G$, let $G_X=\bigcup_{x\in X}G_x$, $G^Y=\bigcup_{y\in Y}G^y$, and $G_X^Y=G_X\cap G^Y$. Sometimes, the unit $x$ will be also
denoted by $1_x$.
\end{notation}

For $x,y\in T$, the subspace $\jet{\germs{\HH}_x^y}=\jetgerms{\HH}_x^y$ is a closed
subspace of $J^1(T)_x^y$ because \HH\ is closed. Then, since $J^1(T)_x^y$ can be considered as the space of
linear maps $\TT_xT\to\TT_yT$, and because the linear maps in $\jetgerms{\HH}_x^y$ are
orthogonal, it follows that $\jetgerms{\HH}_x^y$ is compact. 

\lem{jetgermscalH0}
We have:
\begin{itemize}

\item[$($i$)$] For all $x,y\in T$, $\jetgerms{\HH_0}_x^y$ is closed in
$\jetgerms{\HH}_x^y$, and thus $\jetgerms{\HH_0}_x^y$ is also
compact.

\item[$($ii$)$] $\jetgerms{\HH_0}$ is open in $\jetgerms{\HH}$.

\end{itemize}
\elem

\prf To prove property~(i), let $h_m$ be a sequence of maps in $\HH_0$ with $x\in\dom h_m$ and
$h_m(x)=y$ for all $m$. Suppose the sequence $j^1(\germ{h_m}{x})$ is convergent in $J^1(T)$. Since
$j^1(\germs{\HH})_x^y$ is closed in $J^1(T)$, there is some $h\in\HH$ such that $x\in\dom h$, $h(x)=y$ and $j^1(\germ{h_m}{x})\to j^1(\germ{h}{x})$ as $m\to\infty$.

By the completeness of \HH, we can assume that there are open neighborhoods $U$ and $V$ of $x$ and $y$, and that there is a sequence
$g_m\in\HH$ such that $h$ and every $g_m$ are isometries $U\to V$, and $\germ{g_m}{x}=\germ{h_m}{x}$ for
all $m$. By \refl{jetgerm} and since $j^1(\germ{h_m}{x})$ converges to $j^1(\germ{h}{x})$, it follows that
$g_m\to h$ uniformly on $U$, and thus $g_m^{-1}\circ h\to\id_U$ uniformly. So $g_m^{-1}\circ h\in\HH_0$ for $m$ large
enough. Hence
$$
j^1(\germ{h}{x})=j^1(\germ{g_m}{x})\cdot j^1(\germ{g_m^{-1}\circ h}{x})
=j^1(\germ{h_m}{x})\cdot j^1(\germ{g_m^{-1}\circ h}{x})\;,
$$
which is in
$$
j^1(\germs{\HH_0})_x^y\cdot j^1(\germs{\HH_0})_x^x\subset j^1(\germs{\HH_0})_x^y
$$
for $m$ large enough, as desired.

To prove property~(ii), take any $h\in\HH_0$ and any $x\in\dom h$. We have to show that, for $g\in\HH$
and $y\in\dom g$, if $j^1(\germ{g}{y})$ is close enough to $j^1(\germ{h}{x})$, then $j^1(\germ{g}{y})\in
\jetgerms{\HH_0}$; i.e., the restriction of $g$ to some open neighborhood of $y$ is in
$\HH_0$.

Let $(U,V)$ be a completeness pair of \HH\ with $x\in U$ and $h(x)\in V$. Let
$W$ be a relatively compact connected open neighborhood of $x$ with $\overline{W} \subset\dom h\cap U$ and
$h(\overline{W})\subset V$. Since $j^1(\germ{g}{y})$ is as close as desired to $j^1(\germ{h}{x})$, we can suppose that
$y\in W$ and $\dom g=U$. Then $g$ is uniformly close to $h$ on $W$ by \refl{jetgerm}, and thus
$g\circ h^{-1}$ is close to the identity map on $h(W)$, yielding $g\circ h^{-1}\in\HH_0$ by \reft{Salem}. Now we have
$$
\jetgerm{g}{y}=\jetgerm{g\circ h^{-1}}{h(y)}\cdot\jetgerm{h}{y}\;,
$$
where both $\jetgerm{g\circ h^{-1}}{h(y)}$ and $\jetgerm{h}{y}$ are in
$\jetgerms{\HH_0}$. Thus $\jetgerm{g}{y}$ is also
in $\jetgerms{\HH_0}$ as desired.
\eprf

\cor{jetgermscalH0UV}
If $U$ and $V$ are relatively compact open subsets of $T$, then $\jetgerms{\HH_0}_U^V$ is a
relatively compact open subset of $\jetgerms{\HH_0}$.
\ecor

\prf We know that $\jetgerms{\HH_0}_U^V$ is open in
$\jetgerms{\HH}$ by \refl{jetgermscalH0}-(ii). 

When $P$ and $Q$ run in the family of open neighborhoods of points $x,y\in T$, the subsets
$\jetgerms{\HH_0}_P^Q$ form a base of open neighborhoods of
$\jetgerms{\HH_0}_x^y$ in $\jetgerms{\HH_0}$. Moreover
$\jetgerms{\HH_0}_x^y$ is compact by \refl{jetgermscalH0}-(i). Also,
$\jetgerms{\HH_0}$ is locally compact because it is open in the closed subset
$\jetgerms{\HH}$ of the Hausdorff manifold $J^1(T)$. So $\jetgerms{\HH_0}_P^Q$
is relatively compact in $\jetgerms{\HH_0}$ for $P$ and $Q$ small enough.

Now, for each $x\in\overline{U}$ and $y\in\overline{V}$, choose corresponding open neighborhoods $P_{x,y}$
and $Q_{x,y}$ such that $\jetgerms{\HH_0}_{P_{x,y}}^{Q_{x,y}}$ is relatively compact in
$\jetgerms{\HH_0}$. For each $x\in\overline{U}$, the compact subspace $\overline{V}$ can be
covered by a finite number of open sets $Q_{x,y}$, which are denoted by $Q_{x,j}=Q_{x,y_j}$,
$j=1,\dots,\ell_x$. For $P_{x,j}=P_{x,y_j}$, the open subsets
$$
P_x=P_{x,1}\cap\dots\cap P_{x,\ell_x}\;,\quad x\in\overline{U}\;,
$$
cover the compact subspace $\overline{U}$, and thus there is a finite subcovering consisting of sets $P_i=P_{x_i}$,
$i=1,\dots,k$. Let $Q_{i,j}=Q_{x_i,j}$, $i=1,\dots,k$, $j=1,\dots,\ell_i=\ell_{x_i}$. Then
$$
\jetgerms{\HH_0}_U^V\subset
\bigcup_{i=1}^k\bigcup_{j=1}^{\ell_i}\jetgerms{\HH_0}_{P_i}^{Q_{i,j}}\;,
$$
where each subset of the right hand side is relatively compact in $\jetgerms{\HH_0}$, and
thus $\jetgerms{\HH_0}_U^V$ is relatively compact in $\jetgerms{\HH_0}$.
\eprf

\cor{orbits of calH0}
The orbits of $\HH_0$ are the connected components of the orbits of \HH.
\ecor

\prf First we prove that the orbits of $\HH_0$ are connected. Let $h\in\HH_0$ and $x\in\dom h$. Then
$$
h=\exp\xi_k\circ\dots\circ\exp\xi_1
$$
over some neighborhood $W$ of $x$, where each $\xi_i$ is a local section of $\frakG $. Arguing by
induction on $k$, it is enough to prove that $h(x)$ is in the same connected component of $\HH_0(x)$ as 
$$
h_0(x)=\exp\xi_{k-1}\circ\dots\circ\exp\xi_1(x)\;.
$$
But the composite
$$
h_t=\exp t\xi_k\circ\exp\xi_{k-1}\circ\dots\circ\exp\xi_1\in\HH_0
$$
is also defined on $W$ for all $t\in I$. So $t\mapsto h_t(x)$, $t\in I$, is a curve in $\HH_0(x)$
joining
$h_0(x)$ and $h(x)$.

Second, we prove that $\HH_0(x)$ is open in $\HH(x)$ for each $x\in T$. By \refl{jetgermscalH0}-(ii), the subset $\jetgerms{\HH_0}_x$ is open in $\jetgerms{\HH}_x$. But it is easy to see that the target
projection $\beta:\jetgerms{\HH}_x\to\HH(x)$ is open. Therefore $\HH_0(x)=\beta\left(\jetgerms{\HH_0}_x\right)$ is open in $\HH(x)$. 

The result now follows because each orbit of \HH\ is a disjoint union of orbits of $\HH_0$ since
$\HH_0\subset\HH$.
\eprf

\cor{calH0 is closed}
The pseudogroup of local isometries $\HH_0$ is closed.
\ecor

\prf Since \HH\ is closed, it is enough to prove that $\jetgerms{\HH_0}$ is closed in
$\jetgerms{\HH}$. Let $\sigma_m$ be a sequence in $\jetgerms{\HH_0}$
converging to some $\sigma\in \jetgerms{\HH}$; we have to prove that $\sigma$ is in
$\jetgerms{\HH_0}$. Let $x$ and $y$ be the source and target of $\sigma$, and let $U$ and $V$ be
relatively compact open neighborhoods of $x$ and $y$. Then $\jetgerms{\HH}_U^V$ is an
open neighborhood of
$\sigma\in\jetgerms{\HH}$, and thus we can assume $\sigma_m\in
\jetgerms{\HH}_U^V$ for all $m$; so $\sigma_m$ in $\jetgerms{\HH_0}_U^V$. Hence $\sigma$ is in the closure 
of $\jetgerms{\HH_0}_U^V$ in $\jetgerms{\HH}$. But the closure of $\jetgerms{\HH_0}_U^V$ in
$\jetgerms{\HH}$ equals its closure in $\jetgerms{\HH_0}$ because $\jetgerms{\HH}$ is Hausdorff and $\jetgerms{\HH_0}_U^V$ is relatively compact in $\jetgerms{\HH_0}$ (\refc{jetgermscalH0UV}). Therefore
$\sigma\in\jetgerms{\HH_0}$.
\eprf

\lem{calH0 is complete}
Let $\sigma\in\jetgerms{\HH_0}_x$ for some $x\in T$. Then there is some open neighborhood $V$ of
$x$ such that any $\tau\in\jetgerms{\HH_0}$ close enough to $\sigma$ is the $1$-jet of the
germ of some element of $\HH_0$ with domain $V$.
\elem

\prf
Let $\sigma=\jetgerm{h}{x}$ for some $h\in\HH_0$ and let $P=\dom h$. We can suppose $\overline{P}\subset
U$ for some open set $U$ of those given by \reft{Salem}; i.e., for every relatively compact open set $W$ with
$\overline{W}\subset U$, the elements of \HH\ with domain $W$ and close enough to the identity are those
of the form $\exp\xi$ for $\xi\in\frakG(U)$ small enough. Then take as $V$ any open neighborhood of
$x$ with $\overline{V}\subset P$.

Let $\tau=\jetgerm{g}{y}$ for some $g\in\HH_0$ and $y\in\dom g=W$. Assuming that $W$ is connected, by
\refl{jetgerm}, if $\tau$ is close enough to $\sigma$, we can suppose that $W\subset V$, $g(W)\subset
h(V)$, and $g$ is as close to $h$ on $W$ as desired. Thus the composite $h^{-1}\circ g$ is defined in $W$,
and can be made close enough to the identity to be of the form $\exp\xi$ for some small
$\xi\in\frakG(U)$ by \reft{Salem}. But if $h^{-1}\circ g$ is close enough to the identity on $W$, then
$\xi$ is so small that $\exp\xi$ is defined on $V$ and $\exp\xi(V)\subset P$ because $\overline{V}\subset P$
and
$\overline{P}\subset U$. Thus $h\circ\exp\xi\in\HH_0$ is an extension of $g$ defined on the
whole of $V$.
\eprf

\cor{calH0 is complete}
$\HH_0$ is complete.
\ecor

\prf
Let $x,y\in T$. Since the orbits of $\HH_0$ are closed by \refc{orbits of calH0}, if $x$ and $y$ are not
in the same orbit of $\HH_0$, then there are open neighborhoods $U$ and $V$ of $x$ and $y$ so that no orbit of
$\HH_0$ intersects both $U$ and $V$. Such a pair $(U,V)$ obviously satisfies the completeness condition of
$\HH_0$. Therefore we can assume $x$ and $y$ are in the same orbit of $\HH_0$; so
$\jetgerms{\HH_0}_x^y\neq\emptyset$. Now, by \refl{calH0 is
complete}, each $\sigma\in\jetgerms{\HH_0}_x^y$ has an open neighborhood $\Theta_\sigma$ in
$\jetgerms{\HH_0}$ such that, for some neighborhood $U_\sigma$ of $x$, every element of
$\Theta_\sigma$ is the $1$-jet of some map in $\HH_0$ defined on the whole of $U_\sigma$. Because
$\jetgerms{\HH_0}_x^y$ is compact (\refl{jetgermscalH0}-(i)), it can be covered by a finite number of open
sets $\Theta_i=\Theta_{\sigma_i}$, $i=1,\dots,k$, and let $U_i=U_{\sigma_i}$. Since the open sets
$\jetgerms{\HH_0}_U^V$ of $\jetgerms{\HH_0}$ form a neighborhood base of the compact set
$\jetgerms{\HH_0}_x^y$ when $U$ and $V$ run over the open neighborhoods of $x$ and $y$, there is such a pair $(U,V)$ so
that
$$
\jetgerms{\HH_0}_U^V\subset\Theta_1\cup\dots\cup\Theta_k\;.
$$
Moreover we can assume $U\subset U_1\cap\dots\cap U_k$. Then it is easy to check that $(U,V)$
satisfies the completeness condition of $\HH_0$; indeed, if $h\in\HH_0$ and $z\in U\cap\dom h$, then
$\jetgerm{h}{z}\in\jetgerms{\HH_0}_U^V$, and thus $\jetgerm{h}{z}\in\Theta_i$ for some $i$,
yielding that $\germ{h}{z}$ is the germ at $z$ of some $g\in\HH_0$ with $\dom g=U_i\supset U$.
\eprf

\reft{calH0} is the combination of
Corollaries~\ref{c:orbits of calH0},~\ref{c:calH0 is closed} and~\ref{c:calH0 is complete}.

The following is a version of \reft{Salem} in terms of $\HH_0$ and local actions of local Lie groups.

\th{local action}
Let \HH\ be a closed complete Riemannian pseudogroup on a Riemannian manifold $T$. Then $\HH_0$
is generated by an effective isometric local action on
$T$ of a local Lie group $G$ whose Lie algebra is the structural Lie algebra \frakg\ of \HH.
\eth

\prf
According to \reft{Salem}, the sheaf $\frakG $ defines an infinitesimal action of \frakg\ on each
open set $U$ considered in its statement. By \cite[Corollary~1, page~184, Chapters~2 and~3]{Bourbaki:72Lie}, this
infinitesimal action on each $U$ is induced by a local action of a local Lie group with Lie algebra \frakg.
These local actions can be glued
together by \cite[Proposition~11, page~182, Chapters~2 and~3]{Bourbaki:72Lie}, defining the required local action
on $T$, which is effective and isometric by the definition of $\frakG $.
\eprf

Now suppose that \HH\ is a (possibly non-closed) complete Riemannian pseudogroup on a Riemannian
manifold $T$. By \reft{local action}, the pseudogroup $\overline{\HH}_0$ is generated by an effective
isometric local action of a local Lie group $G$ on $T$. Let $\Lambda$ be the subset of $g\in G$ whose local
action on $T$ defines an element of \HH. It easily follows that such a $\Lambda$ is a dense local 
subgroup of $G$, and the restriction of the local action of $G$ to $\Lambda$ induces the pseudogroup
$\HH_0=\HH\cap\overline{\HH}_0$. We get that $\HH_0$ is complete and
$\overline{\HH_0}=\overline{\HH}_0$.

By \refts{calH0}{local action}, or directly from \reft{Salem}, the connected components of the orbits closures of \HH\ are the leaves of a \cinf\ singular foliation. It follows that any \HH-saturated open subset of $T$ is $\overline{\HH}$-saturated too.

\sec{locally homogeneous orbits}{Locally homogeneous structure of the orbit closures}

Let \HH\ be a complete Riemannian pseudogroup on a Riemannian manifold $T$. To describe locally the orbit closures of \HH, we can assume that \HH\ is closed by \reft{Haefliger}. Even with this assumption, there may not be any pseudogroup in the equivalence class of \HH\ induced by the action of a Lie group \cite{Haefliger:88Leaf} (the orbits would be homogeneous spaces). But, since $\HH_0$ is
generated by the local action of a local Lie group (\reft{local action}), and since the orbits of $\HH_0$
are the connected components of the orbits of \HH\ (\reft{calH0}), the orbits of \HH\ have certain
``local homogeneous structure'' that will be described in this section.

Let $G$ be the local Lie group with an effective isometric left local action on $T$ that
induces $\HH_0$. The identity element of $G$ will be denoted by $e$. Such a local action will be denoted by $\mu:\Omega\to T$, where $\Omega$ is an open neighborhood of $\{e\}\times T$ in $G\times T$. If \frakg\ denotes the Lie algebra of right invariant tangent vector fields on $G$, we have a corresponding infinitesimal action of \frakg\ on $T$.

Fix $x\in T$. Let $\frakk\subset\frakg$ be the
Lie subalgebra of elements of \frakg\ whose infinitesimal transformations of $T$ vanish at $x$, and let \KK\ be the distribution on $G$ defined by the left translates of \frakk; i.e., $\KK_g=L_{g\ast}\frakk$ for each $g\in G$. Such \KK\ is completely integrable because \frakk\ is a subalgebra, and thus defines a foliation on $G$ that will be also denoted by \KK; this is the {\em left foliation\/} on $G$ induced by the Lie subalgebra $\frakk\subset\frakg$ according to the terminology
of \cite[Chapter~III, \S~4.1, pp.~166--167]{Bourbaki:72Lie}.  

Let $\Omega_x=\{g\in G\ |\ (g,x)\in\Omega\}$, which is open in $G$,  and let $\mu_x:\Omega_x\to T$ be defined by $\mu_x(g)=\mu(g,x)$. Since $\mu$ is locally transitive on $\HH_0(x)$, the restriction $\mu_x:\Omega_x\to\HH_0(x)$ is a \cinf\ submersion; in particular, it is open. Thus the connected components of the fibers of $\mu_x$ are the leaves of a foliation, which is easily seen to be equal to the restriction of \KK\ to $\Omega_x$. 

\lem{h n to id U}
If $x_n\to x$ in $\HH(x)$, then there is an open neighborhood $U$ of $x$ in $T$ and a sequence
$h_n\in\HH_U$ such that $h_n(x)=x_n$ and $h_n\to\id_U$ uniformly.
\elem

\prf
This holds because $\mu_x:\Omega_x\to\HH_0(x)$ is open.
\eprf

Let $V$ be an open neighborhood of $e$ in $\Omega_x$ which is simple with respect to \KK. Then  $\mu_x(V)$ is open in $\HH_0(x)$ because $\mu_x:\Omega_x\to\HH_0(x)$ is open. Let $Q$ be the corresponding local quotient of $V$, which is a manifold with a unique \cinf\ structure so that the canonical projection $V\to Q$ is a \cinf\ submersion. Since $\mu_x:\Omega_x\to\HH_0(x)$ is a \cinf\ submersion, it induces a local diffeomorphism $\bar\mu_x:Q\to\HH_0(x)$. We can choose $V$ cutting each fiber of $\mu_x$ at most in one plaque, which means that $\bar\mu_x:Q\to\mu_x(V)$ is a bijection, and thus a diffeomorphism.

\sec{Molino}{Version of Molino's theory for pseudogroups}

The following is the obvious adaptation to pseudogroups of Molino's description of the so called {\em Riemannian foliations\/} \cite{Molino:88Riemannian}. The proofs from \cite{Molino:88Riemannian} can be easily adapted to this setting too. 

A \cinf\ pseudogroup acting on a \cinf\ manifold $T$ is said to be {\em parallelizable\/} if its maps preserve same parallelism of $T$; such a pseudogroup becomes Riemannian by declaring this parallelism to be orthonormal. Suppose that moreover \HH\ is complete. Then the following properties hold:
\begin{itemize}

\item $\overline{\HH}\backslash T$ is a manifold with a unique \cinf\ structure so that the canonical projection $\pi:T\to\overline{\HH}\backslash T$ is a \cinf\ submersion.

\item  Each \HH-orbit closure $F$ has an open neighborhood $U$ in $\overline{\HH}\backslash T$ such that 
\begin{equation}\label{e:calH| pi -1(U)}
\HH|_{\pi^{-1}(U)}\cong\HH|_F\times U\;.
\end{equation}

\end{itemize}

Now, let \HH\ be a complete Riemannian pseudogroup on a Riemannian manifold $T$. The tangent homomorphisms $\TT h$ of maps $h\in\HH$ generate a pseudogroup $\TT\HH$ on the tangent bundle $\TT T$, and the bundle projection $\pi_\TT:\TT T\to T$ generates a morphism $\Pi_\TT:\TT\HH\to\HH$. The pseudogroup $\TT\HH$ is complete and Riemannian with respect to the Sasaki metric on $\TT T$.

Let $\Or(T)$ denote the $\Or(n)$-principal bundle of orthonormal frames of $\TT T$, $n=\dim T$. For any $h\in\HH$, let $\Or(h):\Or(\dom h)\to\Or(\im h)$ be the map defined by
$$
\Or(h)(f_1,\dots,f_n)=(\TT h(f_1),\dots,\TT h(f_n))\;.
$$
The maps $\Or(h)$ generate a pseudogroup $\Or(\HH)$ acting on $\Or(T)$, and the bundle projection $\pi_{\Or}:\Or(T)\to T$ generates a homomorphism $\Pi_{\Or}:\Or(\HH)\to\HH$. The pseudogroup $\Or(\HH)$ is parallelizable and complete, obtaining the following properties:
\begin{itemize}

\item $W=\overline{\Or(\HH)}\backslash\Or(T)$ is a manifold with a unique \cinf\ structure so that the canonical projection $\pi:\Or(T)\to\overline{\Or(\HH)}\backslash\Or(T)$ is a \cinf\ submersion.

\item  The action of $\Or(n)$ on $\Or(T)$ induces a \cinf\ right action of $\Or(n)$ on $W$. The canonical projection $W\to W/\Or(n)$ will be denoted by $\bar\pi$. 

\item The map 
\begin{equation}\label{e:overline calH backslash T}
\overline{\HH}\backslash T\to W/\Or(n)\;,\quad
F\mapsto\bar\pi\circ\pi(\pi_{\Or}^{-1}(F))\;,
\end{equation}
is a homeomorphism, whose inverse is given by $E\mapsto\pi_{\Or}((\bar\pi\circ\pi)^{-1}(E))$.

\item The homeomorphism~\refe{overline calH backslash T} induces a bijection $\cinf(W/\Or(n))\to\cinf(\overline{\HH}\backslash T)$, where $\cinf(W/\Or(n))$ and $\cinf(\overline{\HH}\backslash T)$ are the sets of real valued functions on $W/\Or(n)$ and $\overline{\HH}\backslash T$ inducing \cinf\ functions on $W$ and $T$ via $\bar\pi$ and the canonical projection $T\to\overline{\HH}\backslash T$, respectively.

\end{itemize}

As a first consequence of the above properties, it follows that each locally finite open covering of $\overline{\HH}\backslash T$ has a subordinated partition of unity consisting of functions in $\cinf(\overline{\HH}\backslash T)$.

Consider the nested sequence
$$
\emptyset=T_{-1}\subset T_0\subset\dots\subset T_n=T\;,
$$
where each $T_\ell$ is the union of all orbit closures of \HH\ with dimension $\le\ell$. Since $\overline{\HH}\backslash T$ is homeomorphic to $W/\Or(n)$, it follows that every $T_\ell$ is closed in $T$, each $T_\ell\sm T_{\ell-1}$ is open and dense in $T_\ell$, and $T_\ell\sm T_{\ell-1}$ is an $\overline{\HH}$-invariant \cinf\ submanifold of $T$. Thus the orbit closures define a \cinf\ regular foliation on $T_\ell\sm T_{\ell-1}$ because all of them have the same dimension.

Recall the following definitions from e.g. \cite{Pears:75}. The {\em order\/} of a family of subsets, not all empty, of some set is the largest integer $k$ for which there is a subfamily with $k+1$ elements whose intersection is non-empty, or is $\infty$ if there is no such largest integer. A family of empty sets is declared to have {\em order\/} $-1$. Then the {\em covering dimension\/} $\dim X$ of a space $X$ is the least integer $k$ such that every finite open covering of $X$ has an open refinement of order not exceeding $k$ or is $\infty$ if there is no such integer. By \cite[Theorem~4.3]{Pears:75}, if $X$ is normal, then we can use locally finite open coverings instead of finite ones in the definition of covering dimension; thus we can use arbitrary open coverings when $X$ is normal and paracompact.

For any \cinf\ action of a compact Lie group $G$ on a \cinf\ manifold $M$, consider the {\em isotropy type stratification\/} of $M/G$ \cite{Schwarz:80}: two orbits have the same {\em isotropy type\/} (and thus belong to the same stratum) when they have the same conjugacy classes of isotropy groups in $G$. Any orbit type stratum is a manifold, and is open in its closure, which consists of strata of lower dimension \cite{Schwarz:80}. Then, by applying Proposition~1.5, Theorem~2.5 and Corollary 5.8 of Chapter~3 of \cite{Pears:75} to the inclusion $S\subset\overline{S}$ for each isotropy type stratum $S$, we get by induction on $\dim S$ that $\dim M/G$ equals the top codimension of the orbits. Thus, using the above homeomorphism between $\overline{\HH}\backslash T$ and $W/\Or(n)$, we get that $\dim\overline{\HH}\backslash T$ equals the top codimension of the orbit closures; in particular, it is finite.

\sec{around the orbit closures}{Description around the orbit closures}

In this section, we describe complete Riemannian pseudogroups around the orbit closures. This is an adaptation of the description of transversely complete Riemannian foliations around the leaf closures given by Molino \cite{MolinoPierrot}. A more refined description around the orbit closures, with a complete set of invariants, was given by Haefliger \cite{Haefliger:88Leaf}.

\begin{notation}
Let M be a metric space with distance function $d$. For each $x\in M$ and $r>0$, we use the standard notation $B(x,r)$ for the open $r$-ball in $M$ centered at $x$. For any subset $S\subset M$, $\Pen(F,r)$ denotes the {\em $r$-penumbra\/} of $S$ in $M$, whose definition is
$$
\Pen(F,r)=\{y\in M\ |\ d(y,K)<r\}=\bigcup_{x\in M}B(x,r)\;.
$$
\end{notation}

Let \HH\ be a complete Riemannian pseudogroup on a Riemannian manifold $T$. With the notation of the above section, the exponential map is defined on some open neighborhood $\widetilde{\Omega}$ of the zero section of $\TT T$; let $\widetilde{\Omega}'=\TT\HH(\widetilde{\Omega})$. By elementary properties, $\exp:\widetilde{\Omega}\to T$ generates a morphism $\Exp:\TT\HH|_{\widetilde{\Omega}'}\to\HH$. Let $\Omega=(\pi,\exp)(\widetilde{\Omega})$ and $\Omega'=(\HH\times\HH)(\Omega)$. Observe that, if two orbit closures $F_1$ and $F_2$ are close enough in $\overline{\HH}\backslash T$, then $F_1\times F_2\subset\Omega'$. We can choose $\widetilde{\Omega}$ so that $(\pi,\exp):\widetilde{\Omega}\to\Omega$ is a diffeomorphism, and thus $(\Pi,\Exp):\TT\HH|_{\widetilde{\Omega}'}\to(\HH\times\HH)|_{\Omega'}$ is an isomorphism. Moreover we can assume that $t\cdot\widetilde{\Omega}\subset\widetilde{\Omega}$ for all $t\in I$. Hence the mapping $(v,t)\mapsto t\cdot v$ generates a homotopy $\widetilde{\Psi}:\TT\HH|_{\widetilde{\Omega}'}\times I\to\TT\HH|_{\widetilde{\Omega}'}$, which defines a  homotopy $\Psi:(\HH\times\HH)|_{\Omega'}\times I\to(\HH\times\HH)|_{\Omega'}$ via $(\Pi,\Exp):\TT\HH|_{\widetilde{\Omega}'}\to(\HH\times\HH)|_{\Omega'}$. Observe that $\Psi_0$ is generated by mapping $(x,y)\mapsto(x,x)$, and $\Psi_1$ is the identity morphism at $(\HH\times\HH)|_{\Omega'}$.

For any $x\in T$, let $F=\overline{\HH}(x)$, which is \HH-invariant. The normal subbundle $\TT F^\perp\subset\TT T$ is invariant by $\TT\HH$. For $\epsilon>0$, besides the $\epsilon$-penumbra $\Pen(F,\epsilon)$ of $F$ in $M$, consider the $\epsilon$-penumbra $\widetilde{\Pen}(F,\epsilon)$ of the zero section in $\TT F^\perp$. There is some $\epsilon>0$ so that $\widetilde{\Pen}(F,\epsilon)\subset\widetilde{\Omega}'$ and the restriction 
\begin{equation}\label{e:exp}
\exp:\widetilde{\Omega}\cap\widetilde{\Pen}(F,\epsilon)\to\exp(\widetilde{\Omega})\cap\Pen(F,\epsilon)
\end{equation}
is a diffeomorphism. Let $\HH_{F,\epsilon}$ and $\TT\HH_{F,\epsilon}^\perp$ denote the restrictions of \HH\ and $\TT\HH$ to $\Pen(F,\epsilon)$ and $\widetilde{\Pen}(F,\epsilon)$, respectively. Then the restriction $\Exp:\TT\HH_{F,\epsilon}^\perp\to\HH_{F,\epsilon}$ is an isomorphism generated by~\refe{exp}. We have $\widetilde{\Psi}(\widetilde{\Pen}(F,\epsilon)\times I)\subset\widetilde{\Pen}(F,\epsilon)$, and thus $\widetilde{\Psi}$ induces a homotopy $\Psi_{F,\epsilon}:\HH_{F,\epsilon}\times I\to\HH_{F,\epsilon}$ via $\Exp:\TT\HH_{F,\epsilon}^\perp\to\HH_{F,\epsilon}$. Notice that $\im\Psi_{F,\epsilon,0}\subset F$, $\Psi_{F,\epsilon,1}$ is the identity morphism at $\HH_{F,\epsilon}$, and, for each $x\in\Pen(F,\epsilon)$, the restriction of $\Psi_{F,\epsilon}$ to $\{x\}\times I\equiv I$ is a geodesic segment of \HH\ whose length is $<d(x,F)$ (in the sense of \refs{generalization}). 

Let $f=(f_1,\dots,f_n)\in\Or(T)|_F$. Suppose that $f$ is {\em adapted\/} to $F$ in the sense that $f_1,\dots,f_r\in\TT F$ and $f_{r+1},\dots,f_n\in\TT F^\perp$, where $r=\dim F$. The orbit closure $P=\overline{\Or(\HH)(f)}$ is a principal bundle over $F$ with structural Lie group $H\subset\Or(n)$. Since $P$ consists of frames adapted to $F$, it follows that $H\subset\Or(r)\times\Or(n')$, where $n'=n-r$. Let $H'$ denote the projection of $H$ to $\Or(n')$. The restriction of $\Or(\HH)$ to $P$ will be denoted by $\HH_P$, and the bundle projection $P\to F$ generates a morphism $\HH_P\to\HH_F$. 

On the other hand, let $\Or(\TT F^\perp)$ denote the $\Or(n')$-principal bundle over $F$ of orthonormal frames of $\TT F^\perp$. For each $h\in\HH$, let 
$$
\Or_F^\perp(h):\Or(\TT(F\cap\dom h)^\perp)\to\Or(\TT(F\cap\im h)^\perp)
$$ 
be the map defined by
$$
\Or_F^\perp(h)(e_1,\dots,e_{n'})=(\TT h(e_1),\dots,\TT h(e_{n'}))\;.
$$
The maps $\Or_F^\perp(h)$ generate a pseudogroup $\Or_F^\perp(\HH)$ on $\Or(\TT F^\perp)$, and the bundle projection $\Or(\TT F^\perp)\to F$ generates a morphism $\Or_F^\perp(\HH)\to\HH_F$. There is a canonical map $P\to\Or(\TT F^\perp)$ defined by forgetting the first $r$ components of each frame. Its image is an $H'$-principal bundle over $F$ denoted by $P'$, which is invariant by $\Or_F^\perp(\HH)$. The restriction of $\Or_F^\perp(\HH)$ to $P'$ will be denoted by $\HH_{P'}$; observe that $\HH_{P'}$ has dense orbits. Moreover the canonical projection $P\to P'$ generates a morphism $\HH_P\to\HH_{P'}$.

Let $B_\epsilon$ denote the open ball in $\R^{n'}$ of radius $\epsilon$ and centered at the origin. Consider the diagonal action of $H'$ on $P'\times B_\epsilon$, and let $P'\times_{H'}B_\epsilon$ denote the corresponding quotient space. As usual, we get a canonical identity 
\begin{equation}\label{e:widetilde Pen(F,epsilon)}
\widetilde{\Pen}(F,\epsilon)\equiv P'\times_{H'} B_\epsilon\;.
\end{equation}
Furthermore $\HH_{P'}\times B_\epsilon$ induces  a pseudogroup $\HH_{P'}\times_{H'} B_\epsilon$ acting on the quotient $P'\times_{H'} B_\epsilon$. Then~\refe{widetilde Pen(F,epsilon)} generates an identity 
$$
\TT\HH_{F,\epsilon}^\perp\equiv\HH_{P'}\times_{H'} B_\epsilon\;,
$$
yielding that
\begin{equation}\label{e:calH F,epsilon}
\HH_{F,\epsilon}\cong\HH_{P'}\times_{H'} B_\epsilon
\end{equation}
via the isomorphism $\Exp:\TT\HH_{F,\epsilon}^\perp\to\HH_{F,\epsilon}$, which generalizes~\refe{calH| pi -1(U)} when $U$ is diffeomorphic to an open ball in $\R^{n'}$.

\sec{Riem folns}{Riemannian foliations}

A $C^\infty$ foliation \FF\ on a manifold $M$ is said to be {\em
Riemannian\/} when its holonomy pseudogroup is Riemannian for some
metric. In this case, there is a foliated cocycle of \FF\ consisting of Riemannian submersions for some Riemannian metric on $M$, which is called a {\em bundle-like metric\/} \cite{Reinhart:59bundle-like}, \cite{Molino:88Riemannian}; thus the geodesics orthogonal to the leaves are projected to geodesics by the maps of such foliated cocycle. A characteristic property of bundle-like metrics is that, if any geodesic is orthogonal to the leaves at some point, then it remains orthogonal to the leaves at every point \cite{Reinhart:59bundle-like}, \cite{Molino:88Riemannian}; these geodesics are called {\em horizontal\/}. This condition can be considered for a \cinf\ singular foliation too, obtaining the definition of {\em singular Riemannian foliation\/} \cite{Molino:88Riemannian}. A Riemannian foliation \FF\ on a manifold $M$ is said to be {\em transversely complete\/}  when the horizontal geodesics are complete for some bundle-like metric (a {\em transversely complete\/} bundle-like metric); in this case, $\Hol(\FF)$ is complete \cite{Molino:88Riemannian} (this is a direct consequence of \refl{s} bellow, which is included here for other purposes).

For instance, for any isometric local action of a local Lie group on a Riemannian manifold, the connected components of the orbits are the leaves of a singular Riemannian foliation. Therefore, by \refts{calH0}{local action}, or  by the version of Molino's theory for pseudogroups (\refs{Molino}), the connected components of the orbit closures of a complete Riemannian pseudogroup are the leaves of a singular Riemannian foliation \cite{Haefliger:88Leaf}. It follows that the leaf closures of any transversely complete Riemannian foliation \FF\ on a manifold $M$ are the leaves of a singular Riemannian foliation $\overline{\FF}$, which is also a consequence of Molino's theory \cite{Molino:88Riemannian}. Like for complete Riemannian pseudogroups, it follows that any \FF-saturated open set of $M$ is $\overline{\FF}$-saturated too. 

For any transversely complete Riemannian foliation \FF\  on a manifold $M$, the following properties follow from the corresponding ones for pseudogroups (Sections~\ref{s:Molino} and~\ref{s:around the orbit closures}), or from Molino's theory \cite{Molino:88Riemannian}:
\begin{itemize}

\item The space $M/\overline{\FF}$ is homeomorphic to a space of orbit closures of an action of a compact Lie group. Let $\cinf(M/\overline{\FF})$ be the set of functions $M/\overline{\FF}\to\R$ whose composite with the canonical projection $M\to M/\overline{\FF}$ is a \cinf\ function on $M$. For every locally finite open covering of $M/\overline{\FF}$, there is a subordinated partition of unity consisting of functions in $\cinf(M/\overline{\FF})$.

\item Let $p=\dim\FF$ and $q=\codim\FF$. Consider the nested sequence
$$
\emptyset=M_{-1}\subset M_0\subset\dots\subset M_q=M\;,
$$
where each $M_\ell$ is the union of all leaf closures of \FF\ with dimension $\le p+\ell$. Every $M_\ell$ is closed in $M$. Each $M_\ell\sm M_{\ell-1}$ is open and dense in $M_\ell$, and is a $\overline{\FF}$-saturated \cinf\ submanifold of $T$. The restriction of $\overline{\FF}$ to $M_\ell\sm M_{\ell-1}$ is a regular Riemannian foliation because all leaf closures have the same dimension on this submanifold.

\item The covering dimension of $M/\overline{\FF}$ equals the top codimension of the leaf closures, which is thus finite.

\item There is a description of \FF\ around each leaf closure  that corresponds to~\refe{calH F,epsilon}. A finer description around the leaf closures was also given by Haefliger \cite{Haefliger:88Leaf}.

\end{itemize}

Let $\FF_\TT$ be the foliation on $\TT M$ defined by the following condition. For each $x\in M$ and $v\in\TT_x\FF$, take some simple open neighborhood $U$ of $x$, and let $P$ be the plaque of $U$ that contains $x$. Then a neighborhood of $v$ in the corresponding leaf of $\FF_\TT$ is given by all tangent vectors of the form $X(y)$, where $y\in P$, and $X$ runs in the set of infinitesimal transformations of $\FF|_U$ so that $X(x)=v$. The projection of the leaves of $\FF_\TT$ to the normal bundle $\nu\FF=\TT M/\TT\FF$ are the leaves of a foliation $\FF_\nu$, which is the horizontal lift of \FF\ with respect to the partial Bott connection.  Observe that the bundle projections $\pi_\TT:\TT M\to M$ and $\pi_\nu:\nu\FF\to M$ are foliated maps $\FF_\TT\to\FF$ and $\FF_\nu\to\FF$, respectively. Moreover $\pi_\nu$ restricts to covering maps from the leaves of $\FF_\nu$ to the leaves of \FF.

Via the canonical identity $\nu\FF\equiv\TT\FF^\perp$, $\FF_\nu$ corresponds to a foliation $\widetilde{\FF}$ on $\TT\FF^\perp$. Then $\pi_\TT$ restricts to a foliated map $\widetilde{\FF}\to\FF$, whose restrictions to the leaves are covering maps.

By the above definition of $\FF_\TT$, the mapping $(v,t)\mapsto t\cdot v$ defines a \cinf\ foliated map $H_\TT: \FF_\TT\times\R_{\text{\rm pt}}\to\FF_\TT$, which induces a \cinf\ foliated map $H_\nu: \FF_\nu\times\R\to\FF_\nu$. Obviously, $H_\TT(\TT\FF^\perp\times\R)\subset\TT\FF^\perp$, and the restriction $H_\TT:\TT\FF^\perp\times\R\to\TT\FF^\perp$, which will be denoted by $\widetilde{H}$, corresponds to $H_\nu$ by the canonical identity $\nu\FF\equiv\TT\FF^\perp$.

Since the metric is bundle-like and transversely complete, the exponential map $\exp$ of $M$ is defined on the whole of $\TT\FF^\perp$ and is a foliated map $\widetilde{\FF}\to\FF$. Moreover there is a $\widetilde{\FF}$-saturated open neighborhood $\widetilde{\Omega}$ of the zero section of $\TT\FF^\perp$ such that the map $(\pi_\TT,\exp):\widetilde{\Omega}\to M\times M$ is a \cinf\ embedding, and $t\cdot\widetilde{\Omega}\subset\widetilde{\Omega}$ for all $t\in I$. Then $\Omega=(\pi,\exp)(\widetilde{\Omega})$ is a regular \cinf\ submanifold of $M\times M$, and $(\pi_\TT,\exp):\widetilde{\FF}|_{\widetilde{\Omega}}\to(\FF\times\FF)|_\Omega$ is a foliated diffeomorphism. So $\widetilde{H}$ induces via $(\pi,\exp)$ a \cinf\ foliated homotopy $H:(\FF\times\FF)|_\Omega\times I_{\text{\rm pt}}\to(\FF\times\FF)|_\Omega$. Observe that $H_0(x,y)=(x,x)$ for all $(x,y)\in\Omega$, and $H_1=\id_\Omega$.

Let $F$ be a leaf closure of \FF. For $\epsilon>0$, consider the $\epsilon$-penumbra $\Pen(F,\epsilon)$ in $M$, and the $\epsilon$-penumbra $\widetilde{\Pen}(F,\epsilon)$ of the zero section in $\TT F^\perp$. Let $\FF_{F,\epsilon}=\FF|_{\Pen(F,\epsilon)}$ and $\widetilde{\FF}_{F,\epsilon}=\widetilde{\FF}|_{\widetilde{\Pen}(F,\epsilon)}$. Then $\widetilde{H}$ restricts to a foliated map $\widetilde{\FF}\times\R\to\widetilde{\FF}$.
If $\epsilon$ is small enough, then $\widetilde{\Pen}(F,\epsilon)\subset\Omega$ and $\exp:\widetilde{\FF}_{F,\epsilon}\to\FF_{F,\epsilon}$ is a foliated \cinf\ diffeomorphism. On the other hand, we have $\widetilde{H}(\widetilde{\Pen}(F,\epsilon)\times I)\subset\widetilde{\Pen}(F,\epsilon)$, and thus $\widetilde{H}$ induces a \cinf\ foliated homotopy $H_{F,\epsilon}:\FF_{F,\epsilon}\times I_{\text{\rm pt}}\to\FF_{F,\epsilon}$ via $\exp:\widetilde{\FF}_{F,\epsilon}\to\FF_{F,\epsilon}$. Notice that $\im H_{F,\epsilon,0}\subset F$, $H_{F,\epsilon,1}=\id_{\FF_{F,\epsilon}}$, and the mapping $t\mapsto H_{F,\epsilon}(x,t)$ is a horizontal geodesic segment of length $<d(x,F)$ for each $x\in\Pen(F,\epsilon)$.

\sec{existence}{Existence of complete bundle-like metrics}

Let \FF\ be a transversely complete Riemannian foliation on a manifold $M$, and let $g$ be a transversely complete bundle like metric on M. The goal of this section is to prove the following result, which is an adaptation of a theorem of \cite{NomizuOzeki:existence} to Riemaniann foliations.

\prop{existence} There exists a complete bundle-like metric $g'$ on $M$ such that:
\begin{itemize}

\item[(i)] $g$ and $g'$ define the same orthogonal complement $\TT\FF^\perp$ of $\TT\FF$;

\item[(ii)] $g$ and $g'$ have the same restriction to $\TT\FF^\perp$; and

\item[(iii)] $g'$ is a conformal change of $g$ on the leaves.

\end{itemize}
\eprop

Let $d_\FF$ denote the distance function of the leaves. For each $x\in M$ and $r>0$, let $B_\FF(x,r)$ denote the open $r$-ball in the leaf $L_x$. Also, let $\widetilde B(x,r)$ be the open $r$-ball in $\TT_x^\perp$ centered at the origin, and let $B_\perp(x,r)=\exp(\widetilde B(x,r))$. Observe that, if $r$ is small enough, then $B_\perp(x,r)$ is a \cinf\ local transversal of \FF\ through $x$. For $S\subset M$, let 
\begin{gather*}
\widetilde{\Pen}(S,r)=\bigcup_{x\in S}\widetilde B(x,r)\;,\\
\Pen_\perp(S,r)=\bigcup_{x\in S}B_\perp(x,r)=\exp(\widetilde{\Pen}(S,r))\;.
\end{gather*}

\lem{s}
For each compact subset $K\subset M$, there is some $s>0$ such that:
\begin{itemize}

  \item[(i)] $\exp:\widetilde B(x,s)\to M$ is a \cinf\ embedding transverse to \FF\ for all $x$ in the \FF-saturation $K'$ of $K$; and
  
  \item[(ii)] for any leaf $L$ that meets $K$ and any curve $\gamma:I\to L$, there is a continuous map $h_\gamma:I\times B_\perp(\gamma(0),s)\to M$, which is $($piecewise$)$ \cinf\ if $\gamma$ is $($piecewise$)$ \cinf, and such that  
$$
h_\gamma(\{t\}\times B_\perp(x,s))=B_\perp(\gamma(t),s)
$$
for all $t\in I$, and 
$$
h_\gamma(0,z)=z\;,\quad h_\gamma(I\times\{z\})\subset L_z
$$
for every $z\in B_\perp(x,s)$.
  
\end{itemize} 
\elem

\prf
With the notation of \refs{Riem folns}, the map $\exp:\widetilde\Omega\cap\TT_x\FF^\perp\to M$ is a \cinf\ embedding transverse to \FF\ for each $x\in M$. On the other hand, by the compactness of $K$, there is some $s>0$ such that $\widetilde{\Pen}(K,s)\subset\widetilde{\Omega}$. Since $\widetilde{\Omega}$ is $\widetilde{\FF}$-saturated and $\widetilde{\Pen}(K',s)$ is the $\widetilde{\FF}$-saturation of $\widetilde{\Pen}(K,s)$, it follows that $\widetilde{\Pen}(K',s)\subset\widetilde{\Omega}$, yielding part~(i).

For each $v\in\TT_{\gamma(0)}\FF^\perp$, let $\widetilde L_v$ be the leaf of $\widetilde\FF$ though $v$. Since $\pi_\TT$ restricts to covering maps from the leaves of $\widetilde\FF$ to the leaves of \FF, there is a unique curve $\tilde\gamma_v:I\to\widetilde L_v$ such that $\pi_\TT\circ\tilde\gamma_v=\gamma$. This $\tilde\gamma_v$ is (piecewise) \cinf\ if $\gamma$ is (piecewise) \cinf, and has a \cinf\ dependence on $v$. Then part~(ii) is satisfied with the map $h_\gamma:I\times B_\perp(x,s)\to M$ defined by $h_\gamma(x,z)=\exp(\tilde\gamma_v(t))$, where $v$ is the unique point in $\widetilde B(\gamma(0),s)$ satisfying $\exp(v)=z$.
\eprf

\lem{Pen perp}
For any $S\subset M$ and $s>0$, $\overline{\Pen_\perp(S,s)}$ is compact if and only if $\overline S$ is compact.
\elem

\prf
The ``only if'' part is obvious because $S\subset\Pen_\perp(S,s)$.

If $\overline S$ is compact, then $\widetilde{\Pen}(S,s)$ has compact closure in $\TT\FF^\perp$. Since the domain of its exponential map contains $\TT\FF^\perp$ because $g$ is transversely complete, it follows that $\Pen_\perp(S,s)=\exp(\widetilde{\Pen}(S,s))$ has compact closure in $M$.
\eprf

For $x\in M$ and $r,s>0$, let 
$$
\Pi(x,r,s)=\Pen_\perp(B_\FF(x,r),s)\;.
$$
Observe that
\begin{equation}\label{e:Pi(x,r,s)}
\Pi(x,r,s)\subset B(x,r+s)\;.
\end{equation}

\lem{g complete}
If the leaves are complete Riemannian submanifolds, then $g$ is complete.
\elem

\prf
If the leaves are complete, then $\overline{\Pi(x,r,s)}$ is compact for all $x\in M$ and $r,s>0$ by \refl{Pen perp}. Hence $\overline{B(x,r)}$ is compact for all $x\in M$ and $r>0$ by~\refe{Pi(x,r,s)}, and thus $g$ is complete.
\eprf

Let $\frakr:M\to(0,\infty]$ be the function defined by
$$
\frakr(x)=\sup\{r>0\ |\ \text{$\overline{B_\FF(x,r)}$ is compact}\}\;.
$$
For any leaf $L$ of \FF\ and all $x,y\in L$, we easily get
\begin{equation}\label{e:frakr}
|\frakr(x)-\frakr(y)|<d_\FF(x,y)\;.
\end{equation}

\lem{frakr cont}
\frakr\ is continuous.
\elem

\prf
The continuity of \frakr\ on the leaves follows from~\refe{frakr}, but the statement asserts the continuity of \frakr\ on $M$. 

Fix some $x\in M$, and let $K$ be a compact neighborhood of $x$. Take some $s>0$ satisfying the statement of \refl{s} with this $K$. Let $x_n$ be a sequence in $K$ converging to $x$. By the above observation, we can assume that $x_n\in B_\perp(x,s)$ for all $n$. 

For $0<r'<r<\frakr(x)$ and $n$ large enough, we have $B_\FF(x_n,r')\subset\Pi(x,r,s)$ by \refl{s}, and thus $\overline{B_\FF(x_n,r')}$ is compact by \refl{Pen perp}. It follows that
$$
\frakr(x)\le\liminf_{n\to\infty}\frakr(x_n)\;.
$$

For $0<r<r'$ and $n$ large enough, we get $B_\FF(x,r)\subset\Pi(x_n,r',s)$ by \refl{s}. Hence $\overline{B_\FF(x,r)}$ is compact if $\overline{B_\FF(x_n,r')}$ is compact by \refl{Pen perp}. This yields
$$
\frakr(x)\ge\limsup_{n\to\infty}\frakr(x_n)\;,
$$
and the result follows.
\eprf

Observe that $\frakr^{-1}(\infty)$ is a saturated set: it is the union of complete leaves. Moreover it is closed by \refl{frakr cont}, but it may not be open. Thus there may be points $x$ and $y$ in the same connected component of $M$ such that $\frakr(x)=\infty$ and $\frakr(y)<\infty$. Therefore the argument of \cite{NomizuOzeki:existence} has to be slightly modified to finish the proof of \refp{existence}. We proceed as follows. 

\prf[Proof of \refp{existence}]
Setting $\frac{1}{\infty}=0$ as usual, the mapping $x\mapsto\max\{\frac{1}{\frakr(x)},1\}$ is continuous by \refl{frakr cont}. So there is a \cinf\ function $\omega:M\to\R$ such that $\omega(x)>\max\{\frac{1}{\frakr(x)},1\}$ for all $x\in M$. Let $g'$ be the riemannian tensor on $M$ satisfying the properties~(i)--(iii) of \refp{existence}, whose restriction to the leaves is equal to $\omega^2g$. For each $x\in M$ and every $r>0$, let $B'_\FF(x,r)$ be the open $g'$-ball in $L_x$ of center $x$ and radius $r$. 

\cla{existence}
For all $x\in M$, we have
$$
B'_\FF(x,\frac{1}{3})\subset
\begin{cases}
  B_\FF(x,\frac{\frakr(x)}{2}) & \text{if $\frakr(x)<\infty$}\\
  B_\FF(x,\frac{1}{3}) & \text{if $\frakr(x)=\infty$.}
\end{cases}
$$
\ecla

By \refcl{existence}, $\overline{B'_\FF(x,\frac{1}{3})}$ is compact for all $x\in M$. Hence $g'$ has complete restrictions to the leaves, and \refp{existence} follows by \refl{g complete}.

Let us prove \refcl{existence}. For any leaf $L$ and any \cinf\ curve $c:I\to L$ joining points $x$ and $y$, let $\ell$ and $\ell'$ its lengths defined by $g$ and $g'$, respectively. By the mean value theorem, we have
\begin{align}
\ell'&=\int_0^1\omega(c(t))\,\|\gamma'(t)\|\,dt\notag\\
&=\omega(c(t_0))\int_0^1\|\gamma'(t)\|\,dt\notag\\
&=\omega(c(t_0))\,\ell\label{e:ell'=}
\end{align}
for some $t_0\in I$. 

Assume first that $\frakr(x)<\infty$, and thus $\frakr(z)<\infty$ for all $z\in L$. Let $d'_\FF$ denote the $g'$-distance function on the leaves. Suppose that $d_\FF(x,y)\ge\frac{\frakr(x)}{2}$. Then 
\begin{equation}\label{e:ell'>}
\ell'>\frac{\ell}{\frakr(c(t_0))}
\end{equation}
by~\refe{ell'=}. On the other hand,
$$
|\frakr(c(t_0))-\frakr(x)|<d_\FF(x,c(t_0))\le\ell
$$
by~\refe{frakr}, yielding
$$
\frakr(c(t_0))<\frakr(x)+\ell\;.
$$
Therefore
$$
\ell'>\frac{\ell}{\frakr(x)+\ell}\ge\frac{1}{3}
$$
by~\refe{ell'>} and since $\ell\ge\frac{\frakr(x)}{2}$. Hence $d'_\FF(x,y)\ge\frac{1}{3}$, showing \refcl{existence} in this case.

Suppose now that $\frakr(x)=\infty$. So $\frakr(L)=\infty$ and $\omega(L)=1$. Therefore $\ell'=\ell$ by~\refe{ell'=}, yielding $d'_\FF(x,y)=d_\FF(x,y)$. We get $B'_\FF(x,\frac{1}{3})=B_\FF(x,\frac{1}{3})$, which completes the proof of \refcl{existence}.
\eprf

Let us use the term {\em horizontal metric\/} for the {\em Carnot-CarathŽodory metric\/} $d_H:M\times M\to[0,\infty]$ induced by the {\em polarization\/} $\TT\FF^\perp\subset\TT M$ (see e.g. \cite{Gromov:96C-C}), which is defined as follows. A {\em horizontal curve\/} in $M$ is a piecewise \cinf\ curve orthogonal to the leaves of \FF. Then $d_H(x,y)=\infty$ if there is no horizontal curve between $x$ and $y$, and, otherwise, $d_H(x,y)$ is the infimum of the lengths of horizontal curves between $x$ and $y$. For $x\in M$ and $r>0$, the {\em horizontal ball\/} of {\em radius\/} $r>0$ and {\em center\/} $x$ is the set
$$
B_H(x,r)=\{y\in N\ |\ d_H(x,y)<r\}\;.
$$ 
For $K\subset M$, the {\em horizontal penumbra\/} of {\em radius\/} $r$ around $K$ is the set 
$$
\Pen_H(K,r)=\bigcup_{x\in K}B_H(x,r)\;.
$$ 

\prop{Pen H}
If $K$ is compact, then $\overline{\Pen_H(K,r)}$ is compact.
\eprop

\prf
By \refp{existence}, there is a complete bundle-like metric on $M$ defining the same horizontal penumbras as $g$. So we can assume that $g$ is complete. Then $\overline{\Pen(K,r)}$ is compact, and the result follows because $\Pen_H(K,r)\subset\Pen(K,r)$.
\eprf

\sec{main}{Morphisms between complete Riemannian pseudogroups}

Let $\HH$ and $\HH'$ be complete Riemannian pseudogroups on Riemannian manifolds $T$ and $T'$, respectively, and let
$\FF$ and $\FF'$ be the corresponding possibly singular Riemannian foliations defined by their orbits closures. For any morphism $\Phi:\HH\to\HH'$, every $\phi:U\to T'$ in $\Phi$
is a foliated map $\FF|_U\to\FF'$. The morphism $\Phi$ is said to be of {\em class
\coinf\/} if all of its elements are \coinf\ as foliation maps in the above sense. The set of all
\coinf~morphisms $\HH\to\HH'$ will be denoted by $\coinf(\HH,\HH')$.

Now, our main result (\refmt{main}) can be stated as follows.

\th{main}
With the above notation and conditions, any morphism $\Phi:\HH\to\HH'$ satisfies the following properties:
\begin{itemize}

\item[$($i$)$] $\Phi$ is complete.

\item[$($ii$)$] $\Phi$ generates a morphism $\overline{\Phi}:\overline{\HH}\to\overline{\HH'}$.

\item[$($iii$)$] $\Phi$ is of class \coinf. 

\end{itemize}
\eth

The morphism $\overline{\Phi}$ of \reft{main}-(ii) will be called the {\em closure\/} of $\Phi$. 

\rems{main}
In \reft{main}, observe the following:
\begin{itemize}

\item[(a)] By property~(iii), $\phi$ is \cinf\ if $\HH'$ has dense orbits.

\item[(b)] With the notation of \refs{Riem psgrs} and \reft{main}, observe that $\Phi_{\overline{\text{\rm orb}}}=\overline{\Phi}_{\text{\rm orb}}:\overline{\HH}\backslash T\to\overline{\HH'}\backslash T'$.

\end{itemize}
\erems

\reft{main} will follow from the following proposition, whose large proof is given in
the following section.

\prop{main}
Let \HH, $\HH'$ and $\Phi$ be as in the statement of \reft{main}. For each $\phi\in\Phi$ and any
$x\in\dom\phi$, there is an open neighborhood $U$ of $x$ in $\dom\phi$ satisfying the following properties:
\begin{itemize}

\item[$($i$)$] $(U,U)$ is a completeness pair of $\overline{\HH}$.

\item[$($ii$)$] There exists a neighborhood \OO\ of $\id_U$ in $\overline{\HH}_U$ such that, if
$h\in\OO$, then $h(U)\subset\dom\phi$ and there is some $h'\in\overline{\HH'}$ with $\phi(U)\subset\dom
h'$ and so that $h'\circ \phi=\phi\circ h$ on the whole of $U$.

\item[$($iii$)$] If $h_1',h_2'\in\overline{\HH'}$ satisfy $\phi(U)\subset\dom h_1'\cap\dom h_2'$ and
$h_1'\circ \phi=h_2'\circ \phi$ on some neighborhood of $x$, then $h_1'\circ \phi=h_2'\circ \phi$ on the whole of $U$.

\item[$($iv$)$] The map $\phi$ is \coinf\ on $U$.

\end{itemize}
\eprop

The rest of this section will be devoted to prove that \reft{main} follows from \refp{main}. To begin with,
\reft{main}-(iii) is the same condition as \refp{main}-(iv).

To prove \reft{main}-(ii), by \refl{morphism} it is enough to prove that, given $\phi,\psi\in\Phi$,
$h\in\overline{\HH}$ and $x\in\dom\phi\cap\dom h$ with $h(x)\in\dom\psi$, there is some
$h'\in\overline{\HH'}$ with $\phi(x)\in\dom h'$ and such that $h'\circ\phi=\psi\circ h$ on some neighborhood of
$x$. To prove the existence of such an $h'$, we take the neighborhood $U$ of $y=h(x)$ and the neighborhood
\OO\ of $\id_U$ in $\overline{\HH}_U$ given by \refp{main} for $\psi$ and $y$; moreover we can assume
that $U$ is connected and relatively compact. Let
$$
\sigma=\jetgerm{h}{x}\in\jetgerms{\overline{\HH}}_x=\overline{\jetgerms{\HH}_x}\;.
$$
Then $\sigma$ can be approximated as much as desired by elements
$\tau\in\jetgerms{\HH}_x$, which are of the form $\tau=\jetgerm{f}{x}$ for $f\in\HH$ with $x\in\dom
f$. Thus $\sigma\cdot\tau^{-1}$ is as close as desired to $1_z\in\jet{\overline{\HH}}$, where $z=f(x)$ can be assumed to be in $U$. Because \OO\ is a neighborhood of $\id_U$ in $\overline{\HH}_U$, and since
$U$ is connected and relatively compact, it follows from \refl{jetgerm} that, for $\tau$ close enough to
$\sigma$, there is some $g\in\OO$ such that $\jetgerm{g}{z}=\sigma\cdot\tau^{-1}$. Therefore 
$g=h\circ f^{-1}$ on some neighborhood of $z$ because both maps are local isometries. Now, since $\Phi$ is a morphism
$\HH\to\HH'$, there is some $f'\in\HH'$ so that $\phi(x)\in\dom f'$ and $f'\circ \phi=\psi\circ f$ on some
neighborhood of $x$. On the other hand, from \refp{main}, we get $g(U)\subset\dom\psi$ and the
existence of some $g'\in\overline{\HH'}$ such that $\psi(U)\subset\dom g'$ and $g'\circ\psi=\psi\circ g$ on $U$. It follows that
$$
\psi\circ h=\psi\circ g\circ f=g'\circ\psi\circ f=g'\circ f'\circ\phi
$$
on some neighborhood of $x$, which is the desired property with $h'=g'\circ f'\in\overline{\HH'}$. Thus $\Phi$ generates a morphism $\overline{\Phi}:\overline{\HH}\to\overline{\HH'}$.

\lem{overlinePhi}
The morphism $\overline{\Phi}$ is complete.
\elem

Before proving \refl{overlinePhi}, we show how it implies property~(i) of \reft{main} (the completeness of
$\Phi$). Fix $\phi,\psi\in\Phi$, $x\in\dom\phi$ and $y\in\dom\psi$. Let $U$ and $V$
be neighborhoods of $x$ and $y$ such that $(\phi,U;\psi,V)$ is a completeness quadruple of $\overline{\Phi}$. We
can assume $U$ and $V$ are so small that $\phi(U)\subset U'$ and $\psi(V)\subset V'$ for some
completeness pair $(U',V')$ of $\overline{\HH'}$. We can also suppose that $U$ satisfies the properties of
\refp{main} with respect to $\phi$ and $x$. With these assumptions, we are going to show that
$(\phi,U;\psi,V)$ is also a completeness quadruple of $\Phi$. 

Take $h\in\HH_U$ with $h(U)\cap V\neq\emptyset$. Since $(\phi,U;\psi,V)$ is a completeness
quadruple of $\overline{\Phi}$, we have $h(U)\subset\dom\psi$, and moreover there is some $\bar
h'\in\overline{\HH'}$ satisfying $\phi(U)\subset\dom\bar h'$ and $\bar h'\circ\phi=\psi\circ h$ on $U$. On the other hand, because $h\in\HH$ and $\Phi$ is a morphism $\HH\to\HH'$, there is some $h'\in\HH'$ with $\phi(x)\in\dom h'$ and such that $h'\circ\phi=\psi\circ h$ on some neighborhood of $x$. We can suppose that $\dom h'=U'\supset\phi(U)$ because $(U',V')$ is a completeness pair of $\overline{\HH'}$. So $h'\circ\phi=\bar h'\circ\phi$ on some neighborhood of $x$, and thus also on $U$ by \refp{main}-(iii). Therefore $h'\circ\phi=\psi\circ h$ on $U$ as desired.

\begin{proof}[Proof of \refl{overlinePhi}]
Let $\phi,\psi\in\overline{\Phi}$, $x\in\dom\phi$ and $y\in\dom\psi$. Choose relatively compact open
neighborhoods $P$ and $Q$ of $x$ and $y$ in $T$ and $T'$, respectively, satisfying the following properties:
\begin{itemize}

\item[(A)] There are relatively compact connected open neighborhoods $P_1$ and $Q_1$ of
$\overline{P}$ and $\overline{Q}$ in $\dom\phi$ and $\dom\psi$, respectively, so that $(P_1,Q_1)$ is a completeness pair of $\overline{\HH}$.

\item[(B)] There is an open neighborhoods \OO\ of $\id_U$ in $\overline{\HH}_P$ such that, for all
$h\in\OO$, we have $h(P)\subset\dom\phi$, and moreover there is some $h'\in\overline{\HH'}$ with
$\phi(P)\subset\dom h'$ and $h'\circ\phi=\phi\circ h$ on $P$.

\item[(C)] The space $\jetgerms{\overline{\HH}}_{\overline{P}}^{\overline{Q}}$ is compact.

\end{itemize}
Here, it is obvious that $P$ and $Q$ can be chosen so that property~(A) holds. Property~(B) can be assumed by \refp{main}. Finally, property~(C) can be assumed because $\jetgerms{\overline{\HH}}$ is locally compact, $\jetgerms{\overline{\HH}}_x^y$ is compact, and the family of sets $\jetgerms{\overline{\HH}}_P^Q$ is a
base of open neighborhoods of $\jetgerms{\overline{\HH}}_x^y$ in $\jetgerms{\overline{\HH}}$ when $P$ and $Q$ run over the open neighborhoods of $x$ and $y$.

\cla{overlinePhi}
If some $h\in\overline{\HH}_{P_1}$ satisfies
$h(\overline{P})\cap\overline{Q}\neq\emptyset$, then there is some neighborhood $U_h$ of $x$ in $P$ and
some neighborhood $\GG_h$ of $h$ in $\overline{\HH}_{P_1}$ such that, for all $f\in\GG_h$, there
exists some $f'\in\overline{\HH'}$ with $\phi(U_h)\subset\dom f'$ and $f'\circ\phi=\psi\circ f$ on $U_h$.
\ecla

Fix $h\in\overline{\HH}_{P_1}$ with $h(\overline{P})\cap\overline{Q}\neq\emptyset$ to prove
\refcl{overlinePhi}. Because $\overline{\Phi}$ is a morphism, there is some
$h'\in\overline{\HH'}$ with $\phi(x)\in\dom h'$ and so that $\psi\circ h=h'\circ\phi$ on some neighborhood $W_h$ of
$x$ in $P$. Choose $U_h$ such that $\overline{U_h}\subset W_h$, and choose the neighborhood $\GG_h$ of
$h$ in $\overline{\HH}_{P_1}$ so small that $f(P)\subset h(P_1)$, $g=h^{-1}\circ f|_P\in\OO$ and
$g(U_h)\subset W_h$ for all $f\in\GG_h$. Then, by~(B), there is some $g'\in\overline{\HH'}$ with
$\phi(P)\subset\dom g'$ and $g'\circ\phi=\phi\circ g$ on $P$. Hence $f'=h'\circ g'\in\overline{\HH'}$ satisfies
$\phi(U_h)\subset\dom f'$ and we have
$$
f'\circ\phi=h'\circ g'\circ\phi=h'\circ\phi\circ g=\psi\circ h\circ g=\psi\circ f
$$
on $U_h$, which completes the proof of \refcl{overlinePhi}.

Now, to finish the proof of \refl{overlinePhi}, consider the subspace
$$
\FF=\left\{(h,z)\in\overline{\HH}_{P_1}\times\overline{P}\ \left|\ h(z)\in\overline{Q}\right.\right\}
\subset\overline{\HH}_{P_1}\times\overline{P}\;.
$$
Since $P_1$ is connected and relatively compact, the map
$$
\FF\lar\jetgerms{\overline{\HH}}_{\overline{P}}^{\overline{Q}}\;,\quad(h,z)\mapsto\jetgerm{h}{z}\;,
$$
is a homeomorphism by \refl{jetgerm}. Thus \FF\ is compact by~(C). When $h$ runs over the elements of
$\overline{\HH}_{P_1}$ satisfying $h\left(\overline{P}\right)\cap\overline{Q}\neq\emptyset$, the sets
$$
\FF_h=\FF\cap\left(\GG_h\times\overline{P}\right)\neq\emptyset\;,
$$
form an open covering of \FF, where the sets $\GG_h$ are given by \refcl{overlinePhi}. Then
there is a finite number of elements $h_1,\dots,h_n\in\overline{\HH}_{P_1}$ with
$h_i\left(\overline{P}\right)\cap\overline{Q}\neq\emptyset$ and so that the corresponding sets
$\FF_i=\FF_{h_i}$ cover \FF. Let also $\GG_i=\GG_{h_i}$ and $U_i=U_{h_i}$ according to
\refcl{overlinePhi}. Now let $U$ be any connected open neighborhood of $x$ in
$U_1\cap\dots\cap U_n$, and let $V$ be any open neighborhood of $y$ with $\overline{V}\subset Q$. With this
choice of $U$ and $V$, we are going to prove that $(\phi,U;\psi,V)$ is a completeness quadruple of
$\overline{\Phi}$. 

Take any $f\in\overline{\HH}_U$ with $f(U)\cap V\neq\emptyset$. Because $U$ is connected and $(P_1,Q_1)$
is a completeness pair of $\overline{\HH}$, there is some extension $\bar f\in\overline{\HH}_{P_1}$ of
$f$. We have $\bar f\left(\overline{P}\right)\cap\overline{Q}\supset f(U)\cap V\neq\emptyset$, and thus
there is some $z\in\overline{P}$ such that $\left(\bar f,z\right)\in\FF$. So $\left(\bar
f,z\right)\in\FF_i$ for some $i\in\{1,\dots,n\}$; in particular, $\bar f\in\GG_i$. By
\refcl{overlinePhi}, there is some $f'\in\overline{\HH'}$ with $\phi(U_i)\subset\dom f'$ and
$f'\circ\phi=\psi\circ\bar f$ on $U_i$. Therefore $f'\circ\phi=\psi\circ\bar f$ on $U$ as desired.
\eprf

\sec{proof}{Proof of \refp{main}}

Let \HH, $\HH'$ and $\Phi$ be as in the statement of \refp{main}, and fix $\phi\in\Phi$. To begin with, choose relatively compact connected open subsets, $U_0\subset T$ and
$U_0'\subset T'$, such that:
\begin{itemize}

\item[(A)] $\overline{U_0}\subset\dom\phi$ and $\phi\left(\overline{U_0}\right)\subset U_0'\subset\im\phi$.

\item[(B)] $(U_0,U_0)$ and $(U_0',U_0')$ are completeness pairs of $\overline{\HH}$ and
$\overline{\HH'}$, respectively.

\end{itemize}
Since $\overline{U_0}$ is compact and $\phi$ continuous, we have
$\phi\left(\overline{U_0}\right)=\overline{\phi(U_0)}$ and the restriction
$\phi:\overline{U_0}\to\phi\left(\overline{U_0}\right)$ is a quotient map.

We shall use the following notation. Let $d$ denote the distance function on both $T$ and $T'$. For $y\in T$, $y'\in T'$ and $R>0$, let $B(y,R)$ and $B(y',R)$ be the open balls in $T$ and $T'$ of radius $R$ centered at $y$ and $y'$, respectively. Let $\widetilde{B}(y',R)$ be the open ball in $\TT_{y'}T'$ of radius $R$ centered at the origin. Finally, let $\exp_{y'}$ denote the exponential map from some neighborhood of the origin in $\TT_{y'}T'$ to $T'$. 

Now, consider connected open subsets, $U_1\subset U_0$ and $U_1'\subset U_0'$, and choose
$R,R'>0$ such that the following properties hold:
\begin{itemize}

\item[(C)] $\phi(U_1)\subset U_1'$.

\item[(D)] $\diam(\overline{U_1})<R$, $d(\overline{U_1},T\sm U_0)>R$; thus
$\overline{U_1}\subset U_0$.

\item[(E)] $d(\overline{U_1'},T'\sm U_0')>R'$; thus $\overline{U_1'}\subset U_0'$.

\item[(F)] $\phi(\overline{B(y,R)})\subset B(\phi(y),R')$ for all $y\in\overline{U_1}$.

\item[(G)] The map $\exp_{y'}:\widetilde{B}(y',R')\to B(y',R')$ is defined and is a diffeomorphism for
all $y'\in\overline{U_1'}$.

\item[(H)] $U_1\cap\overline{\HH}(y)$ is connected for all $y\in U_1$; i.e., the orbits of
$\overline{\HH}$ have connected intersections with $U_1$.

\item[(I)] $U_1\cap h_1(U_1)\cap h_1\circ h_2(U_1)\cap\overline{\HH}(y)\neq\emptyset$ for all $y\in U_1$ and all
$h_1,h_2\in\overline{\HH}_{U_0}$ close enough to $\id_{U_0}$.

\end{itemize}
Observe that, for any $x\in\dom\phi$, we can choose neighborhoods $U_i$ and $U'_i$, $i=0,1$, of $x$ and $\phi(x)$
satisfying properties~(A)--(I). Indeed, properties~(A)--(E) can be
obviously assumed; property~(F) can be assumed because $\overline{U_0}$ is a compact subset of 
$\dom\phi$; property~(G) can be assumed because $\overline{U_1'}$ is compact; and finally, properties~(H)
and~(I) can be assumed by the description of a neighborhood of an orbit closure (\refs{around the orbit closures}). 

For $y'\in\overline{U_1'}$, let $\log_{y'}:B(y',R')\to\widetilde{B}(y',R')$ denote the inverse of
$\exp_{y'}:\widetilde{B}(y',R')\to B(y',R')$. Fix $0<r\leq R$ and $y\in\overline{U_1}$. Let $y'=\phi(y)$,
which is in $\overline{U_1'}$ by~(C). Then the subset $\log_{y'}(\phi(B(y,r)))\subset\TT_{y'}T'$ is well
defined by~(F) and~(G), and let
$E(y,r)$ denote the linear span of $\log_{y'}(\phi(B(y,r)))$ in $\TT_{y'}T'$. Finally let
$$
E(y)=\bigcap_{0<r\leq R}E(y,r)\;,
$$
which is a linear subspace of $\TT_{y'}T'$ too. 

\lem{E(y)=E(y,r)}
For all $y\in\overline{U_1}$ there is some $r\in(0,R]$ such that $E(y)=E(y,r)$.
\elem

\prf
This is an easy consequence of the finite dimension of $\TT_{y'}T'$.
\eprf

\lem{h'i}
Let $r\in(0,R]$, $y\in\overline{U_1}$ and $h'_1,h'_2\in\overline{\HH'}_{U_0'}$. Then
$h'_1\circ\phi=h'_2\circ\phi$ on $B(y,r)$ if and only if $h'_{1\ast}=h'_{2\ast}$ on $E(y,r)$.
\elem

\prf
Let $y'=\phi(y)$ and $y'_i=h'_i(y')$, $i=1,2$. Then the diagrams
\begin{equation}\label{e:h'i}
\begin{CD}
\widetilde{B}(y',R')@>{h'_{i\ast}}>>\widetilde{B}(y'_i,R')\\
@V{\exp_{y'}}VV@VV{\exp_{y'_i}}V\\
B(y',R')@>{h'_i}>>B(y'_i,R')
\end{CD}
\end{equation}
are well defined and commutative. We have $h'_1\circ\phi=h'_2\circ\phi$ on $B(y,r)$ if and only if $h'_1=h'_2$
on $\phi(B(y,r))$, which is equivalent to $h'_{1\ast}=h'_{2\ast}$ on $\log_{y'}\phi(B(y,r))$ by the
commutativity of \refe{h'i}, which in turn is obviously equivalent to $h'_{1\ast}=h'_{2\ast}$ on $E(y,r)$. 
\eprf

\cor{h'i}
Let $y\in\overline{U_1}$ and $h'_1,h'_2\in\overline{\HH'}_{U_0'}$. Then $h'_1\circ\phi=h'_2\circ\phi$ on some
neighborhood of $y$ if and only if $h'_{1\ast}=h'_{2\ast}$ on $E(y)$.
\ecor

\prf
This is a direct consequence of \refls{E(y)=E(y,r)}{h'i}.
\eprf

\lem{h'ast}
Let $h\in\HH'_{U_0}$, $h'\in\HH'_{U_0'}$, $y\in\overline{U_1}$ and $y'=\phi(y)$. If $\phi\circ h=h'\circ\phi$
around $y$, then the tangent map $h'_\ast:\TT_{y'}T'\to\TT_{h'(y')}T'$ restricts to an isomorphism
$h'_\ast:E(y)\stackrel{\cong}{\to}E(h(y))$.
\elem

\prf
We have $\phi\circ h=h'\circ\phi$ on $B(y,r)$ for some $r\in(0,R]$ by hypothesis. Then, because
$\phi:\overline{U_0}\to\phi\left(\overline{U_0}\right)$ is a quotient map and since $h:B(y,s)\to B(h(y),s)$
is a homeomorphism for all $s\in(0,R]$, it follows that
$$
h':\phi(B(y,s))\lar\phi(B(h(y),s))
$$
is a homeomorphism if $s\leq r$. So, for $s\leq r$,
$$
h'_\ast:\log_{y'}\phi(B(y,s))\lar\log_{h'(y')}\phi(B(h(y),s))
$$
is a homeomorphism as well because the diagram
$$
\begin{CD}
\widetilde{B}(y',R')@>{h'_\ast}>>\widetilde{B}(h'(y'),R')\\
@V{\exp_{y'}}VV@VV{\exp_{h'(y')}}V\\
B(y',R')@>{h'}>>B(h'(y'),R')
\end{CD}
$$
is well defined and commutative. Therefore $h'_\ast:\TT_{y'}T'\to\TT_{h'(y')}T'$ restricts to an isomorphism
$h'_\ast:E(y,s)\stackrel{\cong}{\to}E(h(y),s)$ for all $s\leq r$, and the result follows.
\eprf

For $X\subset\overline{U_1}$, let $E(X)=\bigcup_{y\in X}E(y)$. The following is some kind of semicontinuity for the spaces $E(y)$.

\lem{E(z) subset bigcap...}
Let $y\in U_1$ and let $V$ be an open subset of $U_1$. If $z\in\overline{V\cap\HH(y)}$, then
$$
E(z)\subset\bigcap_W\overline{E(V\cap W\cap\HH(y))}\;,
$$
where $W$ runs over the open neighborhoods of $z$ in $T$.
\elem

\prf
For $z'=\phi(z)$, let $\Sigma\subset\jetgerms{\overline{\HH'}}_{z'}^{z'}$ denote the subset whose elements
are limits in $\jetgerms{\overline{\HH'}}$ of sequences $\jetgerm{h'_n}{z'}$, where $h'_n\in\HH'_{U'_0}$
satisfies $h'_n\circ\phi=\phi\circ h_n$ around $z$ for some sequence $h_n\to\id_{U_0}$ in $\HH_{U_0}$ with $h_n(z)\in
V$ for all $n$. If $z\not\in V$, then we cannot take $h_n=\id_{U_0}$ and $h'_n=\id_{U'_0}$ for all $n$; thus, {\em a priori\/}, it is not clear that $1_{z'}\in\Sigma$.

\cla{Sigma neq emptyset}
$\Sigma\neq\emptyset$.
\ecla

Let us prove this assertion. Because $z\in\overline{V\cap\HH(y)}$, there is some sequence $h_n\in\HH$
with $z\in\dom h_n$ and $h_n(z)\to z$. Since $(U_0,U_0)$ is a completeness pair of
$\overline{\HH}$, we can assume $h_n\in\HH_{U_0}$ for all $n$. Furthermore we can suppose
$h_n\to\id_{U_0}$ in $\HH_{U_0}$ by \refl{h n to id U}. Since $\Phi$ is a morphism $\HH\to\HH'$ and
$(U'_0,U'_0)$ is a completeness pair of $\overline{\HH'}$, there is a sequence
$h'_n\in\HH'_{U'_0}$ satisfying $\phi\circ h_n=h'_n\circ\phi$ around $z$ for each $n$. The sequence
$\jetgerm{h'_n}{z'}\in\jetgerms{\overline{\HH'}}$ approaches the compact subspace
$\jetgerms{\overline{\HH'}}_{z'}^{z'}$, and thus some subsequence is convergent to some
$\sigma\in\jetgerms{\overline{\HH'}}_{z'}^{z'}$. It follows that $\sigma\in\Sigma$, concluding the proof
of \refcl{Sigma neq emptyset}.

\cla{sigmatau in Sigma}
If $\sigma,\tau\in\Sigma$, then $\sigma\cdot\tau\in\Sigma$.
\ecla

Let us prove \refcl{sigmatau in Sigma}. By the definition of $\Sigma$, there are sequences
$h_m,g_n\in\HH_{U_0}$ and $h'_m,g'_n\in\HH'_{U'_0}$ such that:
\begin{itemize}

\item $\jetgerm{h'_m}{z'}\to\sigma$, $\jetgerm{g'_n}{z'}\to\tau$;

\item $h'_m\circ\phi=\phi\circ h_m$ and $g'_n\circ\phi=\phi\circ g_n$ around $z$;

\item $h_m(z),g_n(z)\in V$; and

\item $h_m\to\id_{U_0}$ and $g_n\to\id_{U_0}$ in $\HH_{U_0}$.

\end{itemize}
For each $m\in\Z_+$, let $P_m$ be a neighborhood of $z$ where $h'_m\circ\phi=\phi\circ h_m$ and
$h_m(P_m)\subset V$. Since $g_n\to\id_{U_0}$, for each $m$ there is some $n_m\in\Z_+$ such
that $g_{n_m}(z)\in P_m$. Moreover we can take $n_m\uparrow\infty$ as $m\to\infty$. Then there is
some neighborhood $Q_m$ of $z$ so that $g_{n_m}(Q_m)\subset P_m$ and $\phi\circ g_{n_m}=g'_{n_m}\circ\phi$ on
$Q_m$. Because $(U_0,U_0)$ and $(U'_0,U'_0)$ are completeness pairs of $\overline{\HH}$ and
$\overline{\HH'}$, respectively, there are maps $f_m\in\HH_{U_0}$ and $f'_m\in\HH'_{U_0}$ that are
equal to $h_m\circ g_{n_m}$ and $h'_m\circ g'_{n_m}$ on some neighborhoods $W_m$ and $W'_m$ of $z$ and $z'$, respectively. We can
assume $\phi(W_m)\subset W'_m$. Hence
$$
\phi\circ f_m=\phi\circ h_m\circ g_{n_m}=h'_m\circ\phi\circ g_{n_m}=h'_m\circ g'_{n_m}\circ\phi=f'_m\circ\phi
$$
in the neighborhood $W_m\cap Q_m$ of $z$, and we have
$$
f_m(z)=h_m\circ g_{n_m}(z)\in h_m\circ g_{n_m}(Q_m)\subset h_m(P_m)\subset V\;.
$$
Moreover
\begin{align*}
\lim_m\jetgerm{f_m}{z}&=\lim_m\left(\jetgerm{h_m}{g_{n_m}(z)}\cdot\jetgerm{g_{n_m}}{z}\right)\\
&=\lim_m\jetgerm{h_m}{g_{n_m}(z)}\cdot\lim_m\jetgerm{g_{n_m}}{z}\\
&=1_z\;,\\
\lim_m\jetgerm{f'_m}{z'}&=\lim_m\left(\jetgerm{h'_m}{g'_{n_m}(z')}\cdot\jetgerm{g'_{n_m}}{z'}\right)\\
&=\lim_m\jetgerm{h'_m}{g'_{n_m}(z')}\cdot\lim_m\jetgerm{g'_{n_m}}{z'}\\
&=\sigma\cdot\tau
\end{align*}
by the properties of the maps $h_m$, $g_n$, $h'_m$ and $g'_n$, because $n_m\uparrow\infty$, and since
$$
g'_{n_m}(z')=\phi\circ g_{n_m}(z)\to\phi(z)=z'\;.
$$ 
So $f_m\to\id_{U_0}$ in $\HH_{U_0}$ by \refl{jetgerm}. It
follows that $\sigma\cdot\tau$ satisfies the conditions to be in $\Sigma$ as desired.

\cla{Sigma is closed}
$\Sigma$ is closed, and thus compact.
\ecla

To prove \refcl{Sigma is closed}, take any sequence $\sigma_m\in\Sigma$ converging to some element
$\tau\in\jetgerms{\overline{\HH'}}_{z'}^{z'}$. By the definition of $\Sigma$, for each $m$, there are
sequences $h_{m,n}\in\HH_{U_0}$ and $h'_{m,n}\in\HH'_{U'_0}$ such that:
\begin{itemize}

\item $\jetgerm{h'_{m,n}}{z'}\to\sigma_m$;

\item $h'_{m,n}\circ\phi=\phi\circ h_{m,n}$ around $z$;

\item $h_{m,n}(z)\in V$;

\item $h_{m,n}\to\id_{U_0}$ in $\HH_{U_0}$ as $n\to\infty$.

\end{itemize}
A typical diagonal convergence argument easily yields the existence of a sequence $n_m\uparrow\infty$ such
that $h_{m,n_m}\to\id_{U_0}$ in $\HH_{U_0}$ and $\jetgerm{h'_{m,n_m}}{z'}\to\tau$.
Therefore $\tau$ satisfies the conditions to be in $\Sigma$ as desired.

The following well known general assertion will be used.

\cla{K is subgroup}
Let $G$ be a first countable compact topological group. For each $g\in G$, the subset
$$
K=\overline{\{g^n\ |\ n\in\Z_+\}}\subset G
$$
is a subgroup.
\ecla

A proof of \refcl{K is subgroup} is included for completeness. First consider the closed subset
$$
L=\bigcap_{m=1}^\infty\overline{\{g^n\ |\ n\ge m\}}\subset K\;.
$$
We have $L\neq\emptyset$ because $G$ is compact. Furthermore $L$ is a subgroup; in fact, if $\sigma,\tau\in
L$, then there are sequences of positive integers, $m_k,n_k\uparrow\infty$, such that $\sigma=\lim_k
g^{m_k}$ and $\tau=\lim_k g^{n_k}$. We can also suppose $m_k-n_k\uparrow\infty$, yielding
$$
\sigma\cdot\tau^{-1}=\lim_kg^{m_k-n_k}\in L\;.
$$
On the other hand, we clearly have $gL\subset L$. Thus $g\in L$, yielding $K=L$
because $L$ is a closed subgroup of $G$. Hence $K$ is a subgroup as desired.

\cla{1z' in Sigma}
$1_{z'}\in\Sigma$.
\ecla

To prove this, we know the existence of some $\sigma\in\Sigma$ by \refcl{Sigma neq emptyset}. Since
$\jetgerms{\overline{\HH'}}_{z'}^{z'}$ is a compact Lie group, the subset
$$
K=\overline{\{\sigma^n\ |\ n\in\Z_+\}}\subset\jetgerms{\overline{\HH'}}_{z'}^{z'}
$$
is a subgroup by \refcl{K is subgroup}. Moreover $K\subset\Sigma$ by \refcls{sigmatau in Sigma}{Sigma is
closed}. Therefore $1_{z'}\in\Sigma$ as desired.

Now the proof of \refl{E(z) subset bigcap...} can be completed as follows. By \refcl{1z' in Sigma}, there are
sequences
$h_n$ and $h'_n$ in $\HH_{U_0}$ and $\HH'_{U'_0}$ such that:
\begin{itemize}

\item $\jetgerm{h'_n}{z'}\to1_{z'}$;

\item $h'_n\circ\phi=\phi\circ h_n$ around $z$;

\item $h_n(z)\in V$;

\item $h_n\to\id_{U_0}$ in $\HH_{U_0}$.

\end{itemize}
By \refl{h'ast}, every tangent homomorphism $h'_{n\ast}:\TT_{z'}T'\to\tau_{h'_n(z')}T'$ restricts to an isomorphism
$h'_{n\ast}:E(z)\stackrel{\cong}{\to}E(h_n(z))$. So, as $W$ runs over the neighborhoods of $z$ in $T$, we
get
\begin{equation}\label{e:bigcapW ...}
\bigcap_W\overline{E(V\cap W\cap\HH(y))}\supset\bigcap_m\overline{\bigcup_{n\ge m}E(h_n(z))}
=\bigcap_m\overline{\bigcup_{n\ge m}h'_{n\ast}(E(z))}\;.
\end{equation}
On the other hand, for each $u\in\TT_{z'}T'$, we have $h'_{n\ast}(u)\to u$ in $\TT T'$
because
$\jetgerm{h'_n}{z'}\to1_{z'}$, yielding
\begin{equation}\label{e:E(z) subset bigcapm ...}
E(z)\subset\bigcap_m\overline{\bigcup_{n\ge m}h'_{n\ast}(E(z))}\;.
\end{equation}
The result now follows from~\refe{bigcapW ...} and \refe{E(z) subset bigcapm ...}.
\eprf

From now on in this section, let $X_y=U_1\cap\overline{\HH}(y)$ and $X'_y=\phi(X_y)$ for each $y\in U_1$.

\cor{phi h=h' phi on some neighborhood of Xy}
Let $y\in U_1$, $h\in\HH_{U_0}$ and $h'\in\HH'_{U'_0}$ such that
$h\left(\overline{U_1}\right)\subset U_0$. If $\phi\circ h=h'\circ\phi$ on some neighborhood of $y$ in $T$, then 
$\phi\circ h=h'\circ\phi$ on some neighborhood of $X_y$ in $T$
\ecor

\prf
Let
$$
A=\left\{z\in U_1\ |\ \phi\circ h=h'\circ\phi\ \text{on some neighborhood of $z$ in $T$}\right\}\;,
$$
Clearly, $A$ is an open subset of $U_1$ and contains $y$. So, since
$X_y$ is connected (property~(H)), it is enough to prove that $A\cap\overline{\HH}(y)$
is closed in $X_y$. 

Let $z$ be a point in the closure of $A\cap\overline{\HH}(y)$ in $X_y$; i.e., 
$$
z\in U_1\cap\overline{A\cap\overline{\HH}(y)}=U_1\cap\overline{A\cap\HH(y)}\;.
$$
We have to prove that $z\in A$.

Because $\Phi$ is a morphism $\HH\to\HH'$
and $(U'_0,U'_0)$ a completeness pair of $\overline{\HH'}$, there is some $h''\in\HH'_{U'_0}$ such that
$\phi\circ h=h''\circ\phi$ on some open neighborhood $P$ of $z$. So $z\in\overline{V\cap\HH(y)}$, where $V=P\cap A$.
Furthermore
$h'\circ\phi=\phi\circ h=h''\phi$ in $V$, yielding $h'_\ast=h''_\ast$ on $E(V)$ by \refc{h'i}. Hence, by
the continuity of $h'_\ast$ and $h''_\ast$, we get $h'_\ast=h''_\ast$ on
$\bigcap_W\overline{E(W\cap V)}$, where $W$ runs over the neighborhoods of $z$. It follows
that $h'_\ast=h''_\ast$ on $E(z)$ by \refl{E(z) subset bigcap...}, and thus $h'\circ\phi=h''\circ\phi$ around $z$ by
\refc{h'i}. Therefore $h'\circ\phi=\phi\circ h$ around $z$, yielding $z\in A$ as desired.
\eprf

\cor{phi h=h' phi on Xy}
Let $y\in U_1$, $h\in\HH_{U_0}$ and $h'\in\HH'_{U'_0}$ such that
$h\left(\overline{U_1}\right)\subset U_0$. If $\phi\circ h=h'\circ\phi$ on some neighborhood of $y$ in
$X_y$, then  $\phi\circ h=h'\circ\phi$ on the whole of $X_y$.
\ecor

\prf
Let \GG\ denote the restriction of \HH\ to $\overline{\HH}(y)$, and let $\Psi:\GG\to\HH'$ denote the
restriction of $\Phi$. With respect to $\Psi$, the open sets $U_i\cap\overline{\HH}(y)$ and $U'_i$, $i=0,1$,
clearly satisfy properties~(A)--(H). Thus the result is a direct consequence of \refc{phi h=h' phi on some neighborhood of Xy}.
\eprf

\cor{phi h-1=h'-1 phi}
There are neighborhoods $\PP_0$ and $\QQ_0$ of $\id_{U_0}$ and $\id_{U'_0}$ in
$\overline{\HH}_{U_0}$ and $\overline{\HH'}_{U'_0}$, respectively, such that all
$h\in\PP_0$ and $h'\in\QQ_0$ satisfy 
$U_1\subset\dom h^{-1}$, $U'_1\subset\dom h^{\prime-1}$, and moreover, if $\phi\circ h=h'\circ\phi$ on $X_y$ for some
$y\in U_1$, then $\phi\circ h^{-1}=h^{\prime-1}\circ\phi$ on $X_y$.
\ecor

\prf
Let $\PP_0$ be the set of $h\in\overline{\HH}_{U_0}$ satisfying 
$h\left(\overline{U_1}\right)\subset U_0$, $\overline{U_1}\subset h(U_0)$ and
$h(X_y)\cap X_y\neq\emptyset$ for all $y\in U_1$. Let $\QQ_0$ be the set of
$h'\in\overline{\HH'}_{U'_0}$ satisfying $\overline{U'_1}\subset h'(U_0)$. By property~(I), the sets 
$\PP_0$ and $\QQ_0$ are neighborhoods of $\id_{U_0}$ and $\id_{U'_0}$, respectively.

Take $h\in\PP_0$ and $h'\in\QQ_0$. We clearly have $U_1\subset\dom h^{-1}$ and $U'_1\subset\dom
h^{\prime-1}$. Now, suppose that $\phi\circ h=h'\circ\phi$ on $X_y$ for some $y\in U_1$; in particular, the following
diagram is commutative:
$$
\begin{CD}
X_y\cap h^{-1}(X_y)@>h>>h(X_y)\cap X_y\\
@V{\phi}VV@VV{\phi}V\\
\phi\left(X_y\cap h^{-1}(X_y)\right)@>{h'}>>\phi\left(h(X_y)\cap X_y\right)\;.
\end{CD}
$$
So $\phi\circ h^{-1}=h^{\prime-1}\circ\phi$ on $h(X_y)\cap X_y$. Furthermore, since $(U_0,U_0)$ and $(U'_0,U'_0)$ are completeness pairs of
$\overline{\HH}'$ and $\overline{\HH'}$, and since $U_1$ and $U'_1$ are connected, there is some $f\in\overline{\HH}_{U_0}$ and some $f'\in\overline{\HH'}_{U'_0}$ that are equal to $h^{-1}$ and $h^{\prime-1}$ on $U_1$ and $U'_1$,
respectively. Thus $\phi\circ f=f'\circ\phi$ on $h(X_y)\cap X_y$. Moreover $h(X_y)\cap X_y\neq\emptyset$ and
$f(\overline{U_1})=h^{-1}(\overline{U_1})\subset U_0$ because $h\in\PP_0$.
Therefore 
$$
\phi\circ h^{-1}=\phi\circ f=f'\circ\phi=h^{\prime-1}\circ\phi
$$ 
on $X_y$ by \refc{phi h=h' phi on Xy}.
\eprf

\cor{phi h1 h2=h'1 h'2 phi}
There are neighborhoods $\PP_1$ and $\QQ_1$ of $\id_{U_0}$ and $\id_{U'_0}$ in
$\overline{\HH}_{U_0}$ and $\overline{\HH'}_{U'_0}$, respectively, such that all
$h_1,h_2\in\PP_1$ and all $h'_1,h'_2\in\QQ_1$ satisfy
$U_1\subset\dom(h_1\circ h_2)$, $U'_1\subset\dom(h'_1\circ h'_2)$, and moreover, if $\phi\circ h_i=h'_i\circ\phi$ on $X_y$
for some $y\in U_1$, $i=1,2$, then $\phi\circ h_1\circ h_2=h'_1\circ h'_2\circ\phi$ on $X_y$.
\ecor

\prf
By property~(I), there is a neighborhood $\PP_{01}$ of $\id_{U_0}$ in
$\overline{\HH}_{U_0}$ such that $X_y\cap h_1(X_y)\cap h_1\circ h_2(X_y)\neq\emptyset$ for all
$h_1,h_2\in\PP_{01}$ and all $y\in U_1$. Take an open set $V$ in $T$ such that $\overline{U_1}\subset V$
and $\overline{V}\subset U_0$. Then the family $\PP_1$ of maps $h\in\PP_{01}$ satisfying
$h(\overline{U_1})\subset V$ and $h(\overline{V})\subset U_0$ is
a neighborhood of $\id_{U_0}$ in $\overline{\HH}_{U_0}$. Let $\QQ_1\subset\overline{\HH'}_{U'_0}$ be
the open subset of maps $h'\in\overline{\HH'}_{U'_0}$ satisfying $h'(\overline{U'_1})\subset
U'_0$

Take $h_1,h_2\in\PP_1$ and $h'_1,h'_2\in\QQ_1$. We clearly have
$U_1\subset\dom(h_1\circ h_2)$ and
$U'_1\subset\dom(h'_1\circ h'_2)$. Now, suppose that $\phi\circ h_i=h'_i\circ\phi$ on $X_y$
for some $y\in U_1$, $i=1,2$. So, for
$$
Y=X_y\cap h_1(X_y)\cap h_1\circ h_2(X_y)\;,
$$ 
the following diagram is commutative:
$$
\begin{CD}
h_2^{-1}\circ h_1^{-1}(Y)@>{h_2}>>h_1^{-1}(Y)@>{h_1}>>Y\\ 
@V{\phi}VV@V{\phi}VV@VV{\phi}V\\
\phi(h_2^{-1}\circ h_1^{-1}(Y))@>{h_2}>>\phi(h_1^{-1}(Y))@>{h_1}>>\phi(Y)\;.
\end{CD}
$$
Hence $\phi\circ h_1\circ h_2=h'_1\circ h'_2\circ\phi$ on $Y$. Furthermore, since
$(U_0,U_0)$ and $(U'_0,U'_0)$ are completeness pairs of
$\overline{\HH}$ and $\overline{\HH'}$, and since $U_1$ and $U'_1$ are connected, there is some 
$f\in\overline{\HH}_{U_0}$ which equals $h_1\circ h_2$ on $U_1$, and there is some 
$f'\in\overline{\HH'}_{U'_0}$ which equals $h'_1\circ h'_2$ on $U'_1$.
Thus $\phi\circ f=f'\circ\phi$ on $h_2^{-1}\circ h_1^{-1}(Y)$. Moreover $Y\neq\emptyset$ and
$$
f(\overline{U_1})=h_1\circ h_2(\overline{U_1})\subset h_1(V)\subset U_0
$$ 
because $h_1,h_2\in\PP_1$. Therefore 
$$
\phi\circ h_1\circ h_2=\phi\circ f=f'\circ\phi=h'_1\circ h'_2\circ\phi
$$ 
on $X_y$ by \refc{phi h=h' phi on Xy}.
\eprf

\cor{there is calPy}
For any $y\in U_1$ and any neighborhood \QQ\ of $\id_{U'_0}$ in $\overline{\HH'}_{U'_0}$, there is some
neighborhood $\PP_y$ of $\id_{U_0}$ in $\overline{\HH}_{U_0}$ such that, for all $h\in\PP_y$, there 
exists some $h'\in\QQ$ such that $\phi\circ h=h'\circ\phi$ on $X_y$.
\ecor

\prf
Fix any $y\in U_1$ and let $y'=\phi(y)$.

\cla{there is some h'}
If $h\in\overline{\HH}_{U_0}$ satisfies $h(\overline{U_1})\subset U_0$, then
there is some $h'\in\overline{\HH'}_{U'_0}$ such that $\phi\circ h=h'\circ\phi$ on $X_y$.
\ecla

Let us prove this statement. By the definition of $\overline{\HH}$, by \refl{jetgerm}, and since
$(U_0,U_0)$ is a completeness pair of
$\overline{\HH}$, we get the existence of a sequence $h_n\in\HH_{U_0}$ that converges to $h$ in
$\overline{\HH}_{U_0}$. We can suppose that
$h_n(y)\in U_0$ for all $n$. Thus there is another sequence $h'_n\in\HH'_{U'_0}$ so that $\phi\circ
h_n=h'_n\circ\phi$ around $y$ because $\Phi$ is a morphism $\HH\to\HH'$ and $(U'_0,U'_0)$ is a completeness
pair of $\overline{\HH'}$. So $\phi\circ h_n=h'_n\circ\phi$ on $X_y$ by
\refc{phi h=h' phi on Xy}. On the other hand, the sequence $\jetgerm{h'_n}{y'}\in\jetgerms{\overline{\HH'}}$
approaches the compact subset $\jetgerms{\overline{\HH'}}_{y'}^{y''}$, where $y'=\phi(y)$ and $y''=\phi\circ
h(y)$. Therefore $\jetgerm{h'_n}{y'}$ has some convergent subsequence; in fact, we can suppose that
$\jetgerm{h'_n}{y'}$ is itself convergent. Its limit is of the form $\jetgerm{h'}{y'}$ for some
$h'\in\overline{\HH'}_{U'_0}$ because $(U'_0,U'_0)$ is a completeness pair of $\overline{\HH'}$.
Moreover $h'_n\to h'$ in $\overline{\HH'}_{U'_0}$ by \refl{jetgerm}. So, from $\phi\circ h_n=h'_n\circ\phi$ on
$X_y$, it follows that $\phi\circ h=h'\circ\phi$ on $X_y$ because $\phi$ is
continuous, completing the proof of \refcl{there is some h'}.

To finish the proof of this lemma, given a sequence $h_n\in\overline{\HH}_{U_0}$ such that
$h_n\to\id_{U_0}$ uniformly, we have to show that there is a subsequence $h_{n_m}$ and a corresponding
sequence
$h''_m\in\overline{\HH'}_{U'_0}$ so that
$h''_m\to\id_{U'_0}$ and $\phi\circ h_{n_m}=h''_m\circ\phi$ on $X_y$ for all $m$. We can 
assume that $h_n(\overline{U_1})\subset U_0$ for all $n$. Then, by \refcl{there is some h'}, there is
a sequence $h'_n\in\overline{\HH'}_{U'_0}$ such that $\phi\circ h_n=h'_n\circ\phi$ on $X_y$. We
have  $h'_n\circ\phi(y)=\phi\circ h_n(y)\to y'$, and thus $\sigma_n=\jetgerm{h'_n}{y'}$ approaches the compact
space $\jetgerms{\overline{\HH'}}_{y'}^{y'}$. Hence some subsequence $\sigma_{n_m}$ is
convergent to some
$\tau\in\jetgerms{\overline{\HH'}}_{y'}^{y'}$. Because $(U'_0,U'_0)$ is a completeness pair of
$\overline{\HH'}$, there is some $f'\in\overline{\HH'}_{U'_0}$ such that $\tau=\jetgerm{f'}{y'}$.
Furthermore $h'_{n_m}\to f'$ on $U'_0$ by \refl{jetgerm}, and thus $f'\circ\phi=\phi$ on $X_y$ because $\phi$ is
continuous, yielding that $f'$ is the identity map on $X'_y$. Therefore there
is some neighborhood of $X'_y$ where the composite
$h'_{n_m}\circ f^{\prime-1}$ is defined for all $m$, and we have
\begin{equation}\label{e:jetgerm h'nm f'-1 to 1y'}
\jetgerm{h'_{n_m}\circ f^{\prime-1}}{y'}=\sigma_{n_m}\cdot\tau^{-1}\lar1_{y'}\;.
\end{equation}
Again, since $(U'_0,U'_0)$ is a completeness pair of $\overline{\HH'}$, there is some
$h''_m\in\overline{\HH'}_{U'_0}$ such that $\sigma_n\cdot\tau^{-1}=\jetgerm{h''_m}{y'}$. Then
$h''_m=h'_{n_m}\circ f^{\prime-1}$ on some neighborhood of $X'_y$ because
this space is connected by property~(H). Therefore
$$
h''_m\circ\phi=h'_{n_m}\circ f^{\prime-1}\circ\phi=h'_{n_m}\circ\phi=\phi\circ h_{n_m}
$$
on $X_y$ because $f'$ is equal to the identity on $X'_y$. Moreover $h''_m\to\id$ on $U'_0$ by \refe{jetgerm h'nm f'-1 to 1y'} and \refl{jetgerm}.
\eprf

\cor{there is calQy}
There is some neighborhood $\PP_2$ of $\id_{U_0}$ in $\overline{\HH}_{U_0}$ such that, for any $y\in
U_1$ and any neighborhood \PP\ of $\id_{U_0}$ in $\overline{\HH}_{U_0}$, there is some neighborhood
$\QQ_y$ of $\id_{U'_0}$ in $\overline{\HH'}_{U'_0}$ such that, if
$h\in\PP_2$ and $h'\in\QQ_y$ satisfy $\phi\circ h=h'\circ\phi$ on $X_y$, then there is some
$f\in\PP$ such that $\phi\circ f=h'\circ\phi$ on $X_y$.
\ecor

\prf
Let $\PP_0$ be the neighborhood of $\id_{U_0}$ in $\overline{\HH}_{U_0}$ given by 
\refc{phi h-1=h'-1 phi}, and choose a neighborhood $\PP_2$ of $\id_{U_0}$ whose closure is contained in
$\PP_0$. 

For each $y\in U_1$, we have to prove that, if some sequences $h_n\in\PP_2$ and
$h'_n\in\overline{\HH'}_{U'_0}$ satisfy $\phi\circ h_n=h'_n\circ\phi$ on $X_y$ and $h'_n\to\id_{U'_0}$,
then there is some subsequence $h'_{n_m}$ and a sequence $f_m\in\overline{\HH}_{U_0}$ such that
$h'_{n_m}\circ\phi=\phi\circ f_m$ on $X_y$ and $f_m\to\id_{U_0}$.

Since $h_n\in\PP_1$, the sequence $\jetgerm{h_n}{y}$ approaches the compact set
$\jetgerms{\overline{\HH}}_y^{\overline{U_0}}$, yielding that some subsequence $h_{n_m}$ is convergent
to some $h$ in $\overline{\HH}_{U_0}$ by \refl{jetgerm}. Since $\phi$ is continuous, we get 
$\phi\circ h=\phi$ on $X_y$, and moreover $h\in\overline{\PP_2}\subset\PP_0$. Therefore 
$U_1\subset\dom h^{-1}$ and $\phi=\phi\circ h^{-1}$ on
$X_y$ by \refc{phi h-1=h'-1 phi}.

Since $h_{n_m}\to h$ in $\overline{\HH}_{U_0}$ and $\overline{U_1}\subset U_0$, it follows that the
composite $h^{-1}\circ h_{n_m}$ is defined on $U_1$ for $m$ large enough. Thus, because
$U_1$ is connected and $(U_0,U_0)$ is a completeness pair of $\overline{\HH}$, it follows that
there is a unique $f_m\in\overline{\HH}_{U_0}$ which is equal to $h^{-1}\circ h_{n_m}$ on $U_1$. From
$h^{-1}\circ h_{n_m}|_{U_1}\to\id_{U_1}$ in $\overline{\HH}_{U_1}$, we get $f_m\to\id_{U_0}$ in
$\overline{\HH}_{U_0}$ by
\refl{jetgerm}. Furthermore
$$
\phi\circ f_m=\phi\circ h^{-1}\circ h_{n_m}=\phi\circ h_{n_m}=h'_{n_m}\circ\phi
$$
on $X_y$ for $m$ large enough.
\eprf

According to \reft{local action}, there are local Lie groups $G$ and $G'$ with local actions on $T$ and $T'$ generating the
pseudogroups $\overline{\HH_0}$ and $\overline{\HH'_0}$, respectively. The identity elements of $G$ and $G'$ will be denoted by
$e$ and $e'$. Since $U_0$ and $U'_0$ are relatively compact, we can suppose that $G$ and $G'$ are small enough around $e$ and $e'$ so
that the local action of all of their elements is defined on the whole of $U_0$ and $U'_0$, respectively. For $g\in G$ and
$g'\in G'$, we shall also denote by $g$ and $g'$ the corresponding elements in
$\overline{\HH}_{U_0}$ and $\overline{\HH'}_{U'_0}$ defined by the local actions. Again, we can assume that $G$ and $G'$
are so small that their elements, considered as maps in
$\overline{\HH}_{U_0}$ and $\overline{\HH'}_{U'_0}$, belong to the neighborhoods
$\PP_i$ and $\QQ_j$ given by Corollaries~\ref{c:phi h-1=h'-1 phi},~\ref{c:phi h1 h2=h'1 h'2 phi} and
\ref{c:there is calQy}, $i=0,1,2$, $j=0,1$. Further assumptions on $G$ and $G'$ will be made when needed.

Fix any $y\in U_1$. 
Let $G_y$ be the set of the elements $g\in G$ such that there is some $g'\in G'$ so that $g'\circ\phi=\phi\circ g$ on
$X_y$. 

\lem{Gy}
$G_y$ is an open local subgroup of $G$.
\elem

\prf
The symmetry of $G_y$ follows from \refc{phi h-1=h'-1 phi}. Also, $G_y$ is a
neighborhood of $e$ in $G$ by \refc{there is calPy}. So it only remains to prove that
$G_y$ is open in $G$

Take any $g\in G_y$, and choose some $g'\in G'$ such that
$g'\circ\phi=\phi\circ g$ on $X_y$. Let $Q$ be a neighborhood of $e'$ in $G'$ so small that, for all
$h'\in Q$, the product $g'h'$ is defined in $G'$. By \refc{there is calPy}, there is some neighborhood $P$
of $e$ in $G$ so small that:
\begin{itemize}

\item for all $h\in P$, there is some $h'\in Q$ such that $h'\circ\phi=\phi\circ h$ on $X_y$; in particular, $P\subset G_y$; and

\item the product $gh$ is well defined in $G$ for all $h\in P$.

\end{itemize}
Then the neighborhood $gP$ of $g$ in $G$ is contained in $G_y$ by \refc{phi h1 h2=h'1 h'2 phi}.
\eprf

Let $G'_y$ be the set of all $g'\in G'$ satisfying $g'\circ\phi=\phi\circ g$ on
$X_y$ for some $g\in G_y$.

\lem{G'y}
$G'_y$ is a local Lie subgroup of $G'$.
\elem

\prf
$G'_y$ clearly contains the identity element $e'$ of $G'$. Moreover $G'_y$ is symmetric by \refc{phi
h-1=h'-1 phi}. 

Now take any $g'_0\in G'_y$ and any $g_0\in G_y$ satisfying $g'_0\circ\phi=\phi\circ g_0$ on $X_y$. Let $P_0$ and $P$
be compact neighborhoods of $g_0$ and $e$ in the open subset $G_y$ of $G$ such that the products
$gh$ and $hg$ are defined in $G$ for all $g\in P_0$ and all $h\in P$. Then, by \refc{there is calQy}, there
are compact neighborhoods $Q_0$ and $Q$ of $g'_0$ and $e'$ in $G'$, respectively, such that:
\begin{itemize}

\item for all $g'\in G'_y\cap Q_0$ and all $h'\in G'_y\cap Q$, there exists some $g\in G_y\cap P_0$ and some
$h\in G_y\cap P$ so that $g'\circ\phi=\phi\circ g$ and $h'\circ\phi=\phi\circ h$ on $X_y$;

\item the products $g'h'$ and $h'g'$ are defined in $G'$ for all $g'\in Q_0$ and all $h'\in Q$.

\end{itemize}
Take $g'\in G'_y\cap Q_0$ and $h'\in G'_y\cap Q$, and choose $g\in G_y\cap P_0$ and $h\in G_y\cap P$
satisfying $g'\circ\phi=\phi\circ g$ and $h'\circ\phi=\phi\circ h$ on $X_y$. Then, by \refc{phi h1 h2=h'1 h'2 phi}, we
get $g'\circ h'\circ\phi=\phi\circ g\circ h$ and $h'\circ g'\circ\phi=\phi\circ h\circ g$ on $X_y$. So both $g'h'$ and $h'g'$ belong to $G'_y$,
and we get that $G'y$ is a local subgroup of $G'$.

Finally, consider a sequence $h'_n\in G'_y\cap Q$ converging to some element $h'\in Q$.
Then there is another sequence $h_n\in P$ such that $h'_n\circ\phi=\phi\circ h_n$ on
$X_y$. Since $P$ is compact, there is a subsequence $h_{n_m}$ converging to some element $h\in P\subset G_y$.
On the other hand, since
$\phi$ is continuous, from $h'_{n_m}\circ\phi=\phi\circ h_{n_m}$ on
$X_y$, it follows that $h'\circ\phi=\phi\circ h$ on $X_y$. So $h'\in G'_y\cap Q$, obtaining that $G'_y\cap Q$ is closed
in $Q$, and thus compact. Hence both $g'_0(G'_y\cap Q)$ and $(G'_y\cap Q)g'_0$ are compact neighborhoods of
$g'_0$ in $G'_y$. Therefore $G'_y$ is a locally compact local subgroup of $G'$, and thus a local Lie
subgroup of $G'$ by \cite[page~227, Th\'eor\`eme~2]{Bourbaki:72Lie}. 
\eprf

\lem{X'y}
$X'_y=\phi(X_y)$ is an open subset of the $G'_y$-orbit of $y'=\phi(y)$ on
$T'$, and thus $X'_y$ is a \cinf~submanifold of $T'$.
\elem

\prf
Since $G_y$ is an open neighborhood of $e$ in $G$, we get that
$\overline{\HH_0}(y)$ is the orbit of the local action of $G_y$ on $T$ that contains $y$. Take any $z\in
X_y$. Because $X_y$ is a connected open subset of $\overline{\HH_0}(y)$, there are $g_1,\dots,g_k\in G_y$
such that $z=g_1\dots g_ky$ and $g_ig_{i+1}\dots g_ky\in X_y$ for $ i=1,\dots,k$. By the definition of $G_y$, there are $g'_1,\dots,g'_k\in G'_y$ such that $g'_i\circ\phi=\phi\circ g_i$ on
$X_y$. Hence $\phi(z)=g'_1\dots g'_ky'$
is in the $G'_y$-orbit of $y'$. Therefore the whole of $X'_y$ is
contained in the $G'_y$-orbit of $y'$ because $z$ is arbitrary.

On the other hand, for all $g'\in G'_y$, we know the existence of some $g\in G_y$ such that
$g'\circ\phi=\phi\circ g$ on $X_y$; in particular, $g'y'=\phi(gy)$. Furthermore, by \refc{there is calQy}, the above
$g$ can be  chosen to approach $e$ in $G_y$ as $g'$ approaches $e'$ in
$G'_y$. So $gy\in X_y$ if $g'$ is close enough to $e'$ in $G'_y$, yielding 
$g'y'=\phi(gy)\in X'_y$. Therefore $X'_y$ contains a nontrivial open subset of the
$G'_y$-orbit of $y'$. It follows that $X'_y$ is open in that orbit as desired.
\eprf

We can assume that $G'$ is an open local subgroup of a local Lie group $\widetilde{G}'$ such that the
product $g'h'$ is defined in $\widetilde{G}'$ for all $g',h'\in G'$, and the local action of $G'$ on $T'$
can be extended to a local action of $\widetilde{G}'$ on $T'$. For the given point $y\in U_1$, two elements
of
$G'_y$ will be said to be equivalent if, as elements of $\overline{\HH'}_{U'_0}$, they have the same
restriction to $X'_y$. The corresponding quotient space will be denoted by $G''_y$.

\lem{G''y}
The local Lie group structure of $G'_y$ canonically induces a local Lie group structure on $G''_y$, and the local action of $G'_y$ on $X'_y$ canonically induces an effective isometric local action of
$G''_y$ on $X'_y$.
\elem

\prf
Let $\widetilde{K}_y\subset\widetilde{G}'$ be the closed normal local Lie subgroup of elements in
$\widetilde{G}'$ that fix all points of $X'_y$. Right translations obviously define a local action of
$\widetilde{K}_y$ on $G'_y$ whose orbit space will be denoted by $G'_y/\widetilde{K}_y$. Since
$\widetilde{K}_y$ is normal in $\widetilde{G}'$, the local group structure of $G'_y$ induces a
local group structure on $G'_y/\widetilde{K}_y$, and the local action of $G'_y$ on $X'_y$ induces
a local action of $G'_y/\widetilde{K}_y$ on $X'_y$.

The pseudogroup on $G'_y$ generated by the local action of $\widetilde{K}_y$ is parallelizable: it preserves the parallelism defined
by any frame of right invariant vector fields. This pseudogroup is closed because
$\widetilde{K}_y$ is closed in $\widetilde{G}'$. Moreover it is obviously complete.
So $G'_y/\widetilde{K}_y$ is a manifold and the quotient map $G'_y\to G'_y/\widetilde{K}_y$ is a
\cinf\ submersion according to \refs{Molino}. Therefore the result easily follows from the following assertion.

\cla{G'y/widetildeKy = G''y}
We have $G'_y/\widetilde{K}_y=G''_y$ as quotient spaces.
\ecla

If the product $g'a'$ is defined in $G'_y$ for some $g'\in G'_y$ and some $a'\in\widetilde{K}_y$, then $g'$
and $g'a'$ have the same restriction to $X'_y$, and thus they are equivalent. Hence, to prove
\refcl{G'y/widetildeKy = G''y}, it only remains to show that two equivalent elements $g',h'\in G'_y$ are in
the same $\widetilde{K}_y$-orbit. But for such $g'$ and $h'$, the inverse $g^{\prime-1}$ is also an element of
$G'_y$, and thus the product $a'=g^{\prime-1}h'$ is defined in $\widetilde{G}'$. Moreover $a'$ fixes every
point of $X'_y$ because $g'$ and $h'$ are equivalent; i.e., $a'\in\widetilde{K}_y$. So $g'$ and $h'=g'a'$ are in
the same $\widetilde{K}_y$-orbit as desired.
\eprf

Let $\frakg$, $\frakg'$, $\frakg'_y$ and $\frakg''_y$ be the Lie algebras of right invariant vector fields on
$G$, $G'$, $G'_y$ and $G''_y$, respectively; we take Lie algebras of right invariant vector fields since we consider left local
actions. By \refl{X'y}, the elements of $\frakg'_y$ are just the elements in $\frakg'$ whose infinitesimal
action is tangent to $X'_y$. If $\frakk_y\subset\frakg'_y$ denotes the normal Lie subalgebra of elements
whose infinitesimal actions vanish on $X'_y$, then there is an induced infinitesimal action of
$\frakg'_y/\frakk_y$ on $X'_y$. 

With the notation of the proof of \refl{G''y}, let $K_y=\widetilde{K}_y\cap G'$, which is an open local
subgroup of $\widetilde{K}_y$, and a closed normal local Lie subgroup of $G'_y$. Then the Lie algebra of
$K_y$ is $\frakk_y$. Moreover, by \refl{G''y}, we have $\frakg''_y\equiv\frakg'_y/\frakk_y$, and the
above infinitesimal action of $\frakg'_y/\frakk_y$ on $X'_y$ can be identified with the infinitesimal action
induced by the local action of $G''_y$ on $X'_y$.

\lem{Fy}
There is a homomorphism of local Lie groups, $F_y:G_y\to G''_y$, such that the restriction
$\phi:X_y\to X'_y$ is $F_y$-equivariant.
\elem 

\prf
For each $g\in G_y$, there is some $g'\in G'_y$ so that $g'\circ\phi=\phi\circ g$ on
$X_y$. By the definition of $G''_y$, any other element of
$G'_y$ satisfying the same property defines the same element in $G''_y$; thus the notation
$F_y(g)\in G''_y$ makes sense for the class represented by $g'$. This defines a map $F_y:G_y\to G''_y$, which is a homomorphism
of local groups by \refc{phi h1 h2=h'1 h'2 phi}. Moreover $\phi:X_y\to X'_y$ is
$F_y$-equivariant, clearly. On the other hand, \refc{there is calPy} implies that $F_y$ is also continuous, and thus
analytic too \cite[Theorem~1, page~225]{Bourbaki:72Lie}. Therefore $F_y$ is a homomorphism of local Lie groups.
\eprf

If $A\in\frakg$ and $A'\in\frakg''_y$, the notation $\widetilde{A}\in\frakX(T)$ and
$\widetilde{B}\in\frakX(X'_y)$ will be used for the corresponding infinitesimal actions.

\cor{phi*}
The map $\phi:X_y\to X'_y$ is \cinf\ and satisfies the following property: if
$A\in\frakg$, $A'=F_{y\ast}(A)\in\frakg''_y$, $z\in X_y$ and $z'=\phi(z)\in X'_y$, then
\begin{equation}\label{e:phi*}
\phi_\ast(\widetilde{A}_z)=\widetilde{A'}_{z'}\;.
\end{equation}
\ecor

\prf
Fix $z\in X_y$ and let $z'=\phi(z)$. According to \refs{locally homogeneous orbits}, there are symmetric open neighborhoods $R$ and $S$ of $e$ and $e''$ in $G_y$ and $G''_y$, and there are open neighborhoods $P$ and $Q$ of $z$ and $z'$ in $X_y$ and $X'_y$, such that the maps $R\to P$ and $S\to Q$,  given by $a\mapsto az$ and $a'\mapsto a'z'$, are well defined \cinf\ surjective submersions. We can assume $\phi(P)\subset Q$ and $F_y(R)\subset S$. The $F_y$-equivariance of $\phi:X_y\to X'_y$
implies the commutativity of the following diagram:
\begin{equation}\label{e:Fy phi}
\begin{CD}
R@>{F_y}>>S\\
@VVV@VVV\\
P@>{\phi}>>Q\;.
\end{CD}
\end{equation}
So $\phi:P\to Q$ is \cinf, and thus $\phi:X_y\to X'_y$ is \cinf\ since $z$ is arbitrary.
Finally, equality~\refe{phi*} is a direct consequence of the commutativity of~\refe{Fy phi}.
\eprf

Now, we can suppose that $\exp:B\to G$ and $\exp:B'\to G'$ are diffeomorphisms for some balls $B$ and $B'$ in
$\frakg$ and $\frakg'$, respectively, centered at the origin. Then, for $-1\le r\le1$, multiplication by $r$ defines maps
$B\to B$ and $B'\to B'$, which correspond via the exponential map to maps $G\to G$ and $G'\to G'$, which will be
respectively denoted by $g\mapsto g^r$ and $g'\mapsto g^{\prime r}$. 
%
%If $S$ is any subset of $G$ or $G'$, then let $S^r$ denote the subset of elements $a^r$ with $a\in S$,
%$-1\le r\le1$. 
%

\cor{g'r phi=phi gr}
If $g'\circ\phi=\phi\circ g$ on $X_y$ for $g\in G$ and $g'\in G'$, then
$g^{\prime r}\circ\phi=\phi\circ g^r$ on $X_y$ for $-1\le r\le1$. In particular, we have $g^r\in G_y$ and
$g^{\prime r}\in G'_y$ for all $g\in G_y$ and all $g'\in G'_y$, and thus $G_y$ and $G'_y$ are connected.
\ecor

\prf
The condition $g'\circ\phi=\phi\circ g$ on $X_y$ means that $g\in G_y$, $g'\in G'_y$ and $F_y(g)=\bar g'$, where
$\bar g'\in G''_y$ denotes the class represented by $g'$. Take $A\in B$ and $A'\in B'$ such that $g=\exp A$ and $\bar g'=\exp A'$. Then the local actions of $g^r$ and $\bar g^{\prime r}$ on
$X_y$ and $X'_y$ are the corresponding uniparametric local groups of diffeomorphisms $\exp r\widetilde{A}$ and $\exp r\widetilde{A'}$.
Therefore it is enough to show that
$\phi_{\ast}(\widetilde{A})=\widetilde{A'}$ on $X'_y$. 
But we have 
$$
\exp F_{y\ast}(A)=F_y\exp(A)=F_y(g)=\bar g'=\exp(A')\;.
$$
So $F_{y\ast}(A)=A'$, and the result follows by \refc{phi*}.
\eprf

Now, fix an arbitrary point $x\in\dom\phi$ and let $x'=\phi(x)$. Take neighborhoods $U_1$ and $U'_1$ of $x$ and $x'$ satisfying properties~(A)--(I). Define $\rho:G\to[0,R]$ in the following way. Given $g\in G$, if there is some
$g'\in G'$ with $\phi\circ g=g'\circ\phi$ around $x$, then let $\rho(g)$ be the supremum of $r\in(0,R]$ such that
$\phi\circ g=g'\circ\phi$ on $B(x,r)$ for some $g'\in G'$;  otherwise let $\rho(g)=0$. The following is a key
result.

\lem{rho ge C}
There is some $C>0$ and some open local subgroup $H\subset G$ such that $\rho(g)\ge C$
for all $g\in H$.
\elem

\prf
Let $C_1=d(x,T\sm U_1)>0$.

\cla{rho(g1 g2)}
We have
$$
\rho(g_1g_2)\ge\min\{C_1,\rho(g_1),\rho(g_2)\}
$$
for all $g_1,g_2\in G$.
\ecla

To prove this assertion, we can assume $\rho(g_i)>0$, $i=1,2$. Set 
$$
0<s<\min\{C_1,\rho(g_1),\rho(g_2)\}\;.
$$ 
Then, by \refc{phi h=h' phi on Xy}, there are $g'_1,g'_2\in G'$ such that
$g'_i\circ\phi=\phi\circ g_i$ on a tube in $U_1$ around $X_x$ of radius $s$. Therefore $(g'_1g'_2)\circ\phi=\phi\circ(g_1g_2)$
on a tube in $U_1$ around $X_x$ of radius
$s$ by \refc{phi h1 h2=h'1 h'2 phi}, and \refcl{rho(g1 g2)} follows.

\cla{rho(gr)}
We have
$$
\rho(g^r)\ge\min\{\rho(g),C_1\}
$$
for $-1\le r\le1$ and $g\in G$.
\ecla

We can assume $\rho(g)>0$ to prove this claim, and set $0<s<\min\{\rho(g),C_1\}$. Then,
by \refc{phi h=h' phi on Xy}, there is some $g'\in G'$
such that $g'\circ\phi=\phi\circ g'$ on the tube in $U_1$ around $X_x$ of radius $s$. So
$g^{\prime r}\circ\phi=\phi\circ g^r$ on such a tube for $-1\le r\le1$ by \refc{g'r phi=phi gr}, and \refcl{rho(gr)}
follows.

Let $\Lambda\subset G$ be the dense local subgroup whose local action on $T$ generates $\HH_0$. Since
$\Phi$ is a morphism $\HH\to\HH'$ and $(U'_0,U'_0)$ is a completeness pair of $\HH'$, for all
$g\in\Lambda$, there is some $h'\in\HH'_{U'_0}$ so that $h'\circ\phi=\phi\circ g$ around $x$. But this does not directly imply $\rho(g)>0$ because $h'$ may not be defined by the action of some element of $G'$; recall that the local action of $G'$ on $T'$ generates $\overline{\HH'_0}$, which may be different from
$\overline{\HH'}$. So the following assertion is not completely obvious, where the notation $g^{(0,1]}=\{g^r\ |\ r\in(0,1]\}$ is used for any $g\in G$.

\cla{rho>0 on Lambda}
For all $g\in\Lambda$ and any neighborhood $U$ of $g^{(0,1]}$ in $G$, there is some $a\in\Lambda\cap U$
with $\rho(a)>0$.
\ecla

Fix any $g\in\Lambda$ and any neighborhood $U$ of $g^{(0,1]}$ in $G$. Take some $h'\in\HH'_{U'_0}$ such
that $h'\circ\phi=\phi\circ g$ on some neighborhood $O$ of $x$ in $U_1$. It is easy to see that there are
sequences, $g_n\in\Lambda$ and $h'_n\in\HH'_{U'_0}$, so that:
\begin{itemize}

\item[(a)] $h'_n\circ\phi=\phi\circ g_n$ on some neighborhood $O_n$ of $x$ in $O$;

\item[(b)] $g_n\to g$;

\item[(c)] $g_n(O_n)\subset g_m(O_m)$ if $m<n$;

\item[(d)] $g_m^{-1}g_n\in U$ if $m<n$.

\end{itemize}
Recall the notation $x'=\phi(x)$, and let $y=g(x)$ and $y'=\phi(y)$. Let also $x'_n=h'_n(x')=\phi\circ g_n(x)$.
Since $\phi$ is continuous, from~(a) and~(b) we get that $\jetgerm{h'_n}{x'}$ approaches the compact
subspace $\jetgerms{\overline{\HH'}}_{x'}^{y'}$ of $J^1(T')$. So some subsequence of
$\jetgerm{h'_n}{x'}$ is convergent in $\jetgerms{\overline{\HH'}}_{x'}^{y'}$. Therefore there are $m$ and $n$
large enough, say $m<n$, so that 
$$
\jetgerm{h^{\prime-1}_m\circ h'_n}{x'}=\jetgerm{h'_m}{x'_n}^{-1}\cdot\jetgerm{h'_n}{x'}
$$
is as close to $1_{x'}$ as desired. In particular, we can assume $h^{\prime-1}_m\circ h'_n$ is defined by the
local action of some $a'_{m,n}\in G'$ by \refl{jetgerm} and \reft{local action}.
On the other hand, from~(c) we get a commutative diagram
$$
\begin{CD}
O_n@>{g_n}>>g_m(O_m)@>{g_n^{-1}}>>O_m\\
@V{\phi}VV@V{\phi}VV@VV{\phi}V\\
\phi(O_n)@>{h'_n}>>\phi\circ g_m(O_m)@>{h^{\prime-1}_m}>>\phi(O_m)\;.
\end{CD}
$$
So $a'_{m,n}\circ\phi=\phi\circ g_m^{-1}g_n$ on $O_n$, yielding $\rho(g_m^{-1}g_n)>0$, and \refcl{rho>0 on Lambda}
follows by~(d).

The proof of \refl{rho ge C} can be completed as follows. Since $\Lambda$ is dense in $G$, there are
elements $g_1,\dots,g_k$ in $\Lambda$ which are the image by the exponential map of the elements of some base
of \frakg. It is easy to see that there is a neighborhood $U_i$ of each $g_i^{(0,1]}$
such that any choice of elements $a_i\in U_i$, $i=1,\dots,k$, is the image by the exponential map of some
base of \frakg. So all elements in some open local subgroup $H\subset G$ can be
written as products of the form $a_1^{r_1}\dots a_k^{r_k}$ with $-1\le r_i\le1$. Moreover, by
\refcl{rho>0 on Lambda} such elements $a_i$ can be chosen so that $\rho(a_i)>0$. Then the result easily
follows from \refcls{rho(g1 g2)}{rho(gr)}.
\eprf

Let us finish the proof of \refp{main}. For $\phi$ and $x$ as above, by \refl{rho ge C}, there is an open
local subgroup $H\subset G$ and some open connected neighborhood $U$ of $x$ in $U_1$ such that, for all
$h\in H$, there is some $h'\in G'$ so that $h'\circ\phi=\phi\circ h$ on $U$. Indeed, we can assume $X_y\subset U$ if
$y\in U$ by 
\refc{phi h=h' phi on Xy}. So \refp{main}-(ii) follows with such a $U$ and taking as \OO\
the neighborhood of $\id_{U}$ in $\overline{\HH}_{U}$ defined by the local action of $H$ on $U$. Because
$U\subset U_0$, this $U$ also satisfies \refp{main}-(i); i.e., $(U,U)$ is a completeness
pair of $\overline{\HH}$. We can choose such a $U$ satisfying \refp{main}-(iii) by
\refls{E(y)=E(y,r)}{h'i}.

Only \refp{main}-(iv) remains to be proved; that is, it remains to show that $\phi$ is \coinf\ on
$U$. This will take some more work, but the main ideas were already used; it is only needed a slight sharpening of the arguments by using the existence of $H$ and $U$.

Let $H'$ be the set of elements $h'\in G'$ such that $h'\circ\phi=\phi\circ h$ on $U$ for some $h\in
H$. Observe that $H'$ contains $e'$.

\lem{there is P}
For any neighborhood $Q$ of $e'$ in $H'$, there is some
neighborhood $P$ of $e$ in $H$ such that, for all $g\in P$, there 
exists some $g'\in Q$ such that $\phi\circ g=g'\circ\phi$ on $U$.
\elem

\prf
Given a sequence $g_n\in H$ such that
$g_n\to e$, we have to show that there is a subsequence $g_{n_m}$ and a corresponding
sequence $h'_m\in H'$ so that $h'_m\to e'$ and $\phi\circ g_{n_m}=h'_m\circ\phi$ on $U$ for all $m$. We know that
there is a sequence $g'_n\in H'$ such that $\phi\circ g_n=g'_n\circ\phi$ on $U$. Then $g'_n\circ\phi=\phi\circ g_n\to
\phi$ on $U$, and thus $\sigma_n=\jetgerm{g'_n}{x'}$ approaches the compact space
$\jetgerms{\overline{\HH'}}_{x'}^{x'}$ (recall that $x'=\phi(x)$). Hence some subsequence $\sigma_{n_m}$ is
convergent to some
$\tau\in\jetgerms{\overline{\HH'}}_{x'}^{x'}$. Because $(U'_0,U'_0)$ is a completeness pair of
$\overline{\HH'}$, there is some $f'\in\overline{\HH'}_{U'_0}$ such that $\tau=\jetgerm{f'}{x'}$.
Furthermore $g'_{n_m}\to f'$ on $U'_0$ by \refl{jetgerm}. So $f'\circ\phi=\phi$ on $U$, and thus $f'$ is the
identity on $\phi(U)$.  Hence there
is some neighborhood of $\phi(U)$ where the composite
$g'_{n_m}\circ f^{\prime-1}$ is defined for all $m$, and we have
\begin{equation}\label{e:jetgerm g'nm f'-1 to 1y'}
\jetgerm{g'_{n_m}\circ f^{\prime-1}}{x'}=\sigma_{n_m}\cdot\tau^{-1}\lar1_{x'}\;.
\end{equation}
Again, since $(U'_0,U'_0)$ is a completeness pair of $\overline{\HH'}$ and $U$ is connected, there is some
$h'_m\in\overline{\HH'}_{U'_0}$ such that $\sigma_n\cdot\tau^{-1}=\jetgerm{h'_m}{x'}$. Then
$h'_m=g'_{n_m}\circ f^{\prime-1}$ on some neighborhood of $\phi(U)$. Therefore
$$
h'_m\circ\phi=g'_{n_m}\circ f^{\prime-1}\circ\phi=g'_{n_m}\circ\phi=\phi\circ g_{n_m}
$$
on $U$ because $f'$ equals the identity on $\phi(U)$. Moreover $h'_m\to\id_{U'_0}$ by \refe{jetgerm
g'nm f'-1 to 1y'} and \refl{jetgerm}; thus the maps $h'_m$ can be considered as elements of $H'$ for $m$
large enough.
\eprf

\lem{there is Q}
For any neighborhood $P$ of $e$ in $G$, there is some neighborhood
$Q$ of $e'$ in $G'$ such that for all $g'\in H'\cap Q$ there is some
$f\in H\cap P$ such that $\phi\circ f=h'\circ\phi$ on $U$.
\elem

\prf
We have to prove that, if some sequences $g_n\in H$ and
$g'_n\in H'$ satisfy $\phi\circ g_n=g'_n\circ\phi$ on $U$ and $g'_n\to e'$,
then there is some subsequence $g'_{n_m}$ and a sequence $f_m\in H$ such that
$g'_{n_m}\circ\phi=\phi\circ f_m$ on $U$ and $f_m\to e$.

The sequence $\jetgerm{g_n}{x}$ approaches the compact set
$\jetgerms{\overline{\HH}}_x^{\overline{U_0}}$. So some subsequence $g_{n_m}$ is convergent
to some $h$ in $\overline{\HH}_{U_0}$ by \refl{jetgerm}. Since $\phi$ is continuous, we get 
$\phi\circ h=\phi$ on $U$, and moreover $h\in\overline{\PP_2}\subset\PP_0$. Therefore 
$U_1\subset\dom h^{-1}$ and $\phi=\phi\circ h^{-1}$ on
$U$ by \refc{phi h-1=h'-1 phi}.

Since $g_{n_m}\to h$ in $\overline{\HH}_{U_0}$ and $\overline{U_1}\subset U_0$, the
composite $h^{-1}\circ h_{n_m}$ is defined on $U_1$ for $m$ large enough. Thus, because
$U_1$ is connected and $(U_0,U_0)$ is a completeness pair of $\overline{\HH}$, there is a unique $f_m\in\overline{\HH}_{U_0}$ which equals $h^{-1}\circ h_{n_m}$ on $U_1$. From the convergence of $h^{-1}\circ h_{n_m}$ to the identity map on $U_1$, we get $f_m\to\id_{U_0}$ in $\overline{\HH}_{U_0}$ by \refl{jetgerm}. Furthermore
$$
\phi\circ f_m=\phi\circ h^{-1}\circ g_{n_m}=\phi\circ g_{n_m}=g'_{n_m}\circ\phi
$$
on $U$ for $m$ large enough. Finally, since $f_m\to\id_{U_0}$, the maps $f_m$ can be considered as
elements of $G$, and thus of $H$, for $m$ large enough.
\eprf

\lem{H'}
$H'$ is a local Lie subgroup of $G'$.
\elem

\prf
$H'$ contains the identity element, and is symmetric by \refc{phi h-1=h'-1 phi}. 

Now take any $g'_0\in H'$ and any $g_0\in H$ satisfying $g'_0\circ\phi=\phi\circ g_0$ on $U$. Let $P_0$ and $P$
be compact neighborhoods of $g_0$ and $e$ in $H$, respectively, such that the products
$gh$ and $hg$ are defined in $G$ for all $g\in P_0$ and all $h\in P$. Then, by \refl{there is Q}, there
are compact neighborhoods $Q_0$ and $Q$ of $g'_0$ and $e'$ in $G'$, respectively, such that:
\begin{itemize}

\item for all $g'\in H'\cap Q_0$ and all $h'\in H'\cap Q$, there exists some $g\in P_0$ and some
$h\in P$ so that $g'\circ\phi=\phi\circ g$ and $h'\circ\phi=\phi\circ h$ on $U$;

\item the products $g'h'$ and $h'g'$ are defined in $G'$ for all $g'\in Q_0$ and all $h'\in Q$.

\end{itemize}
Take $g'\in H'\cap Q_0$ and $h'\in H'\cap Q$, and choose $g\in P_0$ and $h\in P$
satisfying $g'\circ\phi=\phi\circ g$ and $h'\circ\phi=\phi\circ h$ on $U$. Then, by \refc{phi h1 h2=h'1 h'2 phi}, we
get $(g'h')\circ\phi=\phi\circ (gh)$ and $(h'g')\circ\phi=\phi\circ (hg)$ on $U$. So both $g'h'$ and $h'g'$ belong to $H'$,
and we get that $H'$ is a local subgroup of $G'$.

Finally, consider a sequence $h'_n\in H'\cap Q$ converging to some element $h'\in Q$.
Then there is another sequence $h_n\in P$ such that $h'_n\circ\phi=\phi\circ h_n$ on
$U$. Since $P$ is compact, there is a subsequence $h_{n_m}$ converging to some element $h\in P$.
On the other hand, since
$\phi$ is continuous, from $h'_{n_m}\circ\phi=\phi\circ h_{n_m}$ on
$U$, it follows that $h'\circ\phi=\phi\circ h$ on $U$. So $h'\in H'\cap Q$, obtaining that $H'\cap Q$ is closed
in $Q$, and thus compact. Hence both $g'_0(H'\cap Q)$ and $(H'\cap Q)g'_0$ are compact neighborhoods of
$g'_0$ in $H'$. Therefore $H'$ is a locally compact local subgroup of $G'$, and thus a local Lie
subgroup of $G'$ by \cite[page~227, Th\'eor\`eme~2]{Bourbaki:72Lie}. 
\eprf

\lem{H'-orbits}
For each $y\in U$, $X'_y$ is open in the $H'$-orbit of $y'=\phi(y)$ on $T'$.
\elem

\prf
Since $H$ is an open neighborhood of $e$ in $G$, we get that
$\overline{\HH_0}(y)$ is the orbit of the local action of $H$ on $T$ that contains $y$. Take any $z\in
X_y$. Then, because $X_y$ is a connected open subset of $\overline{\HH_0}(y)$, there are $g_1,\dots,g_k\in H$
such that $z=g_1\dots g_ky$ and $g_ig_{i+1}\dots g_ky\in X_y$ for $i=1,\dots,k$. By definition of $H$, there are $g'_1,\dots,g'_k\in H'$ such that $g'_i\circ\phi=\phi\circ g_i$ on
$U$. Hence $\phi(z)=g'_1\dots g'_ky'$ is in the $H'$-orbit of $y'$. Therefore the whole $X'_y$ is
contained in the $H'$-orbit of $y'$ because $z$ is arbitrary.
Now the result follows from \refl{X'y}.
\eprf

Two elements of $H'$ will be said to be equivalent if they have the same action on $\phi(U)$. The
corresponding quotient space will be denoted by $H''$.

\lem{H''}
The local Lie group structure of $H'$ canonically induces a local Lie group structure on $H''$, and the
local action of $H'$ on $\phi(U)$ canonically induces a effective isometric local action of
$H''$ on some $G'$-submanifold
$U''$ of $U'_1$ that contains $\phi(U)$.
\elem

\prf
Let $\widetilde{L}\subset\widetilde{G}'$ be the closed normal local Lie subgroup of elements that fix the
points of $\phi(U)$. Let $U''$ be the space of points in $U'_1$ that are fixed by $\widetilde{L}$. Since
the local action of $\widetilde{L}$ is isometric, $U''$ is a \cinf\ submanifold of
$U_1$. This follows because the local action of $\widetilde{K}$ on $T'$ can be linearized around each $y'\in U''$ via the exponential map, and the fixed points of this linearized local action is a linear subspace of $\TT_{y'}T'$; thus normal coordinates of $T'$ around points of $U''$ restrict to local coordinates on $U''$; and the changes of such local coordinates on $U''$ are \cinf, obviously. Furthermore $U''$ is a
$G'$-subset of $U'_1$: if $y'\in U''$, $g'\in G'$ and $a'\in\widetilde{L}$, there is some
$b'\in\widetilde{L}$ such that $a'g'=g'b'$ because $\widetilde{L}$ is normal in $\widetilde{G}'$, yielding
$a'g'y'=g'b'y'=g'y'$, and thus $g'y'\in U''$.

Now, as in the proof of \refl{G''y}, the local Lie group structure
of $H'$ induces a local Lie group structure on $H''$, and the local action of
$H'$ on $U''$ induces a local action of $H''$ on $U''$.
\eprf

With the notation of the proofs of \refls{G''y}{H''}, observe that
$$
H'\subset\bigcap_{y\in U}G'_y\;,\quad\widetilde{L}\subset\bigcap_{y\in U}\widetilde{K}_y\;.
$$
So there is a canonical homomorphism $H''\to G''_y$ for each $y\in U$. The identity element of $H''$ will be
also denoted by $e''$.

\lem{F}
There is a homomorphism of local Lie groups, $F:H\to H''$, such that the diagram
\begin{equation}\label{e:F}
\begin{CD}
H@>F>>H''\\
@VVV@VVV\\
G_y@>{F_y}>>G''_y
\end{CD}
\end{equation}
is commutative for all $y\in U$.
\elem

\prf
For each $h\in H$, there is some $h'\in H'$ so that $g'\circ\phi=\phi\circ g$ on
$U$. By the definition of $H''$, any other element in
$H'$ satisfying the same property represents the same element in $H''$; thus the notation
$F(h)\in H''$ makes sense for the element represented by $h'$. This defines a map $F:H\to H''$, which is a homomorphism of local groups by \refc{phi h1 h2=h'1 h'2 phi}. On the other hand, \refl{there is P} implies that $F$ is also
continuous, and thus it is a homomorphism of local Lie groups \cite[Theorem~1, page~225]{Bourbaki:72Lie}. The
commutativity of~\refe{F} is easy to check.
\eprf

\cor{F}
The restriction $\phi:U\to U''$ is $F$-equivariant.
\ecor

\prf
This follows from \refl{Fy} and the commutativity of \refe{F}.
\eprf 

Let $\mu:\Omega\to U$ and $\mu'':\Omega''\to U''$
denote the local actions of $H$ on $U$ and of $H''$ on $U''$, where $\Omega\subset H\times U$ and $\Omega''\subset H'\times U''$ are open neighborhoods of $\{e\}\times U$ and $\{e''\}\times U''$. The equivariance of $\phi$ means that the following diagram is commutative:
\begin{equation}\label{e:mu mu''}
\begin{CD}
\Omega@>{F\times\phi}>>\Omega''\\
@V{\mu}VV@VV{\mu''}V\\
U@>{\phi}>>U''\;.
\end{CD}
\end{equation}

Let $\GG$ and $\GG''$ be the foliations on $\Omega$ and $\Omega''$  whose leaves are the the intersections of $\Omega$ and $\Omega''$ with the corresponding product slices $H\times\{y\}$ and $H''\times\{y''\}$, for $y\in U$ and $y''\in U''$. Also, the $\overline{\HH}$-orbits $X_y$ in $U$ are the leaves of a singular Riemannian foliation \FF\ on $U$,
and the $\overline{\HH'}$-orbits define a singular Riemannian foliation $\FF''$ on $U''$ since $U''$ is a $G'$-submanifold of $U'_1$ (\refl{H''}). With respect to these singular Riemannian foliations, \refe{mu mu''} consists of foliated maps. Moreover the top map of \refe{mu mu''} is \coinf, obviously, and the vertical maps $\mu$ and $\mu''$
of \refe{mu mu''} are \cinf\ surjective submersions. So the composite
$$
\begin{CD}
\GG@>{\mu}>> \FF @>{\phi}>> \FF''
\end{CD}
$$
is \coinf\ by the commutativity of \refe{mu mu''}. Therefore the fact that $\phi$ is \coinf\ follows from 
\refl{coinf} since $\mu$ satisfies the following key property.

\lem{there is some W}
For each $V\in\frakX(\FF)$, there is some $W\in\frakX(\GG)$ so that $\mu_\ast(W)=V$.
\elem

\prf
Since $H$ is an open local group of $G$,  the Lie algebra of right invariant vector fields
on $H$ can be also identified with \frakg, and for each
$A\in\frakg$, let $\widetilde{A}\in\frakX(\FF)$ be the corresponding infinitesimal action. Since
infinitesimal actions of elements of \frakg\ generate $\frakX(\FF)$ as $\cinf(U)$-module, it is enough to
check the result when $V=\widetilde{A}$ for some
$A\in\frakg$. In this case, it is easy to check that the result holds with $W=(A,0)\in\frakX(\GG)$, where we use the
canonical injection of $\frakX(H)\oplus\frakX(U)$ into $\frakX(H\times U)$.
\eprf

\sec{SP}{The strong plaquewise topology}

In this section, we give a version for continuous foliated maps of the strong and weak compact-open topologies.  Roughly speaking, two foliated maps will be close with respect to this topology when, on foliated charts, they have the same representations with transverse coordinates, and close representations with plaquewise coordinates.

Let $X$ and $Y$ be foliated spaces with foliated structures $\FF$ and $\GG$. Fix the following data:
\begin{itemize} 

\item Any foliated map $f:\FF\to\GG$.

\item Any locally finite collection $\UU=\{U_i\}$ of simple open sets of $X$. 

\item A family $\VV=\{V_i\}$ of simple open sets of $Y$, with the same index set. Let
$p_i:V_i\to T_i$ be the canonical projection to the local quotient.

\item A family  $\KK=\{K_i\}$ of compact subsets of $X$, with the same index set, such that $K_i\subset
U_i$ and $f(K_i)\subset V_i$ for all $i$.

\end{itemize}
Let $\NN(f,\UU,\VV,\KK)$ be the set of foliated maps
$g:\FF\to\GG$ such that $g(K_i)\subset V_i$ and $p_i\circ g(x)=p_i\circ f(x)$ for each $i$ and every $x\in K_i$. Such sets $\NN(f,\UU,\VV,\KK)$ form a base of a topology on $C(\FF,\GG)$, called the {\em strong plaquewise topology\/}, and the corresponding space
is denoted by $C_{SP}(\FF,\GG)$. A weak version of this topology can be defined by taking finite families; it can be called the {\em weak plaquewise topology\/}, and the corresponding space can be denoted by $C_{WP}(\FF,\GG)$. Both of these topologies are equal if $X$ is compact; in this case, we can use the term {\em plaquewise topology\/} and the notation $C_P(\FF,\GG)$. 

From now on in this section, assume that $X$ is locally compact and Polish.

\th{integrable homotopy}
With the above notation and conditions, if two continuous foliated maps
$f,g:\FF\to\GG$ are close enough in $C_{SP}(\FF,\GG)$, then there is an integrable
homotopy between them. If moreover $f$ and $g$ are proper, then there is a proper integrable
homotopy between them.
\eth

The proof of \reft{integrable homotopy} uses the following lemma.

\lem{integrable homotopy}
Let $f,g:\FF\to\GG$ be continuous foliated maps. Suppose that, for each
$x\in M$, there exists some simple open set of $Y$ such that some of its plaques contain
$f(x)$ and $g(x)$. Then there is an integrable homotopy between $f$ and $g$. If moreover $f$ and $g$ are proper, then there is a proper integrable homotopy between them.
\elem

\prf
It easily follows from the hypotheses that we can get the following:
\begin{itemize}

\item A locally finite family $\{U_i\}$ of open subsets of $X$;

\item A covering of $X$ by compact subset, $\{K_i\}$, with the same index set and such that $K_i\subset U_i$ for all $i$;

\item For each $i$, a foliated chart $\theta_i:V _i\to T_i\times B$ of $\GG$ such that $f(x)$ and $g(x)$ lie in the same plaque of $V_i$ for each $x\in U_i$. Here, $B$ is an open ball of $\R^n$, where $n=\dim\GG$. Let $p_i:V_i\to T_i$ and $q_i:V_i\to B$ denote the composites of $\theta_i$ with the factor projections of $T_i\times B$.

\end{itemize}
Since $\{U_i\}$ is locally finite, its index set is countable; say, either $\Z_+$, or $\{1,\dots,k\}$ for some $k\in\Z_+$. Set $f_0=f$.

\cla{integrable homotopy}
For each $i$, there is a continuous foliated map $f_i:\FF\to\GG$ such that:
\begin{itemize}

\item $f_i(U_j)\subset V_j$ for all $j$;

\item $f_i(x)$ and $f_{i-1}(x)$ lie in the same plaque of $V_i$ for all $i$ and $x\in U_i$;

\item $f_i(x)=g(x)$ for all $i$ and all $x\in K_1\cup\dots\cup K_i$;

\item $f_i=f_{i-1}$ on some neighborhood of $M\sm U_i$ for all $i$;

\item there is an integrable homotopy $H_i:\FF\times I\to\GG$ between $f_{i-1}$ and $f_i$ for all $i$; and

\item if $f_{i-1}$ is proper, then $f_i$ and $H_i$ are proper.

\end{itemize}
\ecla

This assertion is proved by induction on $i$.  We
have already set $f_0=f$. Now suppose that $f_{i-1}$ is given for some $i>0$. Let
$\lambda_i:X\to I$ be a continuous map such that $\supp\lambda_i\subset U_i$ and 
$\lambda_i(x)=1$ for all $x\in K_i$. Then it is easy to show that the statement of
\refcl{integrable homotopy} holds with the following definitions of $f_i$ and $H_i$. 
Let $H_i(x,t)=f_{i-1}(x)$ for $x\in X\sm U_i$, 
$$
H_i(x,t)=
\theta_i^{-1}(p_i\circ g(x),(1-t\,\lambda_i(x))\cdot q_i\circ f_{i-1}(x)+
t\,\lambda_i(x)\cdot q_i\circ g(x))
$$
for $x\in U_i$, and $f_i(x)=H_i(x,1)$ for any $x\in M$.

The result follows from \refcl{integrable homotopy} in the following way. Take
a partition 
$$
0=t_1<\dots<t_k<t_{k+1}=1
$$ 
of $I$ when $i$ runs in $\{1,\dots,k\}$, and take any
sequence $t_i\uparrow1$ with $t_1=0$ when $i$ runs in $\Z_+$. In any case, let
$\rho_i:[t_i,t_{i+1}]\to I$ be any orientation preserving homeomorphism for each $i$.
Then it is easy to check that an integrable homotopy $H:\FF\times I\to\GG$ between $f$ and $g$
can be defined by 
$$
H(x,t)=
\begin{cases}
H_i(x,\rho_i(t))&\text{if $t\in[t_i,t_{i+1}]$ for some $i$,}\\
g(x)&\text{if $t=1$.}
\end{cases}
$$
Moreover $H$ is proper if so $f$ and $g$ are proper.
\eprf

\prf[Proof of \reft{integrable homotopy}]
Consider a basic neighborhood of $f$ in $C_{SP}(\FF,\GG)$ of the form 
$\NN(f,\UU,\VV,\KK)$ defined above. 
We can suppose that $\KK$ covers $X$.  
Then, by \refl{integrable homotopy}, there is an integrable homotopy
between $f$ and any map in $\NN(f,\UU,\VV,\KK)$.
\eprf

\cor{SP-close foliated maps}
If two foliated maps $f,g:\FF\to\GG$ are close enough in
$C_{SP}(\FF,\GG)$, then $\Hol(f)=\Hol(g)$.
\ecor

\prf
This follows from \refp{integrable homotopy} and \reft{integrable homotopy}.
\eprf

Now, assume that $Y$ is locally compact and Polish too, and let $\Proper(\FF,\GG)$ denote the set of proper continuous foliated maps $\FF\to\GG$.

\th{SP proper}
With the above notation and conditions, $\Proper(\FF,\GG)$ is open in $C_{SP}(\FF,\GG)$.
\eth

\prf
Given $f\in\Proper(\FF,\GG)$, take a neighborhood $\NN=\NN(f,\UU,\VV,\KK)$ of $f$ in $C_{SP}(\FF,\GG)$ as above. 

\cla{SP proper}
We can choose \UU, \VV\ and \KK\ such that \KK\ covers $X$ and \VV\ is locally finite.
\ecla

To prove this assertion, let $\Lambda$ be the set of indices $i$.
We can assume that \KK\ covers $X$ because $X$ is locally compact and Polish. On the other hand, since $Y$ is locally compact and Polish, it has a locally finite open covering \WW. Then, for each $i\in\Lambda$, the set
$$
\WW_i=\{W\in\WW\ |\ W\cap f(K_i)\neq\emptyset\}
$$
is finite, and there is an expression
$$
K_i=\bigcup_{W\in\WW_i}K_{i,W}\;,
$$
where each $K_{i,W}$ is compact and so that $f(K_{i,W})\subset W$. For each $i\in\Lambda$ and $W\in\WW_i$, there is a finite covering $\VV'_{i,W}$ of $f(K_{i,W})$ by simple open sets which are uniform in $U_i$. Then, for each $i$ and $W\in\WW_i$, there is another expression
$$
K_{i,W}=\bigcup_{V'\in\VV'_{i,W}}K_{i,W,V'}\;,
$$
where each $K_{i,W,V'}$ is compact and so that $f(K_{i,W,V'})\subset V'$. Now let
$$
\Lambda'=\{(i,W,V')\ |\ i\in\Lambda,\ W\in\WW_i,\ V'\in\VV'_{i,W}\}\;,
$$
and consider the indixed sets 
$$
\UU'=\{U'_{i,W,V'}\}\;,\quad\VV'=\{V'_{i,W,V'}\}\;,\quad\KK'=\{K'_{i,W,V'}\}\;,
$$ 
with $(i,W,V')\in\Lambda'$, where 
$$
U'_{i,W,V'}=U_i\;,\quad V'_{i,W,V'}=V'\;,\quad K'_{i,W,V'}=K_{i,W}\;.
$$
Then $\UU'$, $\VV'$ and $\KK'$ satisfy the conditions of \refcl{SP proper}, and $\NN(f,\UU',\VV',\KK')$ is defined because each $V'\in\VV'_{i,W}$ is uniform in $V_i$.

When the properties of \refcl{SP proper} are satisfied, for any compact subset $R\subset Y$, there is a finite subfamily $\Lambda_R\subset\Lambda$ such that 
$$
i\in\Lambda\sm\Lambda_R\;\Longrightarrow R\cap V_i=\emptyset\;.
$$
On the other hand, by the conditions of \refcl{SP proper} and the definition of \NN,  for any $g\in\NN$ and $i\in\Lambda$, we have:
\begin{itemize}

\item $g(X)\subset\bigcup_{i\in\Lambda}V_i$; and

\item there is some finite subset $\Lambda_{g,i}\subset\Lambda$ such that 
$g^{-1}(V_i)\subset\bigcup_{j\in\Lambda_{g,i}}K_j$.

\end{itemize}
Therefore
$$
g^{-1}(R)\subset\bigcup_{i\in\Lambda_R}g^{-1}(V_i)
\subset\bigcup_{i\in\Lambda_R}\bigcup_{j\in\Lambda_i}K_j\;,
$$
and thus $g^{-1}(R)$ is compact. So $\NN\subset\Proper(\FF,\GG)$, and the result follows.
\eprf

\sec{coinf approx}{\coinf\ approximations of foliated maps}

\ssec{general}{The general case}

With the notation of the above section, suppose that \FF\ and \GG\ are \coinf\ (\refs{foliated sp}), and assume that $X$ is locally compact and Polish. 
%We will consider \coinf\ foliated maps $\FF\to\GG$ (\refs{foliated map}).

\th{coinf}
With the above notation and conditions, suppose that some $f\in C(\FF,\GG)$ is of class $C^{0,\infty}$ on some 
neighborhood of a closed subset $E\subset X$. Then any neighborhood of $f$ in
$C_{SP}(\FF,\GG)$ contains some $g\in C^{0,\infty}(\FF,\GG)$ such that $g=f$ on some
neighborhood of $E$. In particular, $C^{0,\infty}(\FF,\GG)$ is dense in $C_{SP}(\FF,\GG)$.
\eth

The first step to prove \reft{coinf} is the following local result. 
Let $T$ and $T'$ be topological spaces, and let $L$ and $L'$ be connected open subsets of euclidean spaces; assume that $T$ is locally compact and Polish. Let $\FF$ and $\GG$ be the \cinf\ foliations on $X=T\times L$  and $Y=T'\times L'$ with leaves
$\{x\}\times L$ and $\{x'\}\times L'$ for $x\in T$ and $x'\in T'$, respectively.
For any $u\in C(\FF,\GG)$, write $u(x,y)=(\bar u(x),\tilde u(x,y))$ for all
$(x,y)\in X$. Obviously, the foliated map $u$ is of class \coinf\ if and only if the mapping $y\mapsto\tilde u(x,y)$ is \cinf\ for each $x\in T$ and its partial derivatives of arbitrary order depend continuously on $x$.

\lem{coinf}
Let $Q\subset X$ be a closed subset and $W\subset X$ an open subset. If some $u\in
C(\FF,\GG)$ is of class \coinf\ on some neighborhood of $Q\sm W$, then any
neighborhood of $u$ in $C_{SP}(\FF,\GG)$ contains some map $v$
which is \coinf\ on some neighborhood of $Q$, and which equals $u$ on $X\sm W$.
\elem

\prf
Fix any neighborhood $\NN$ of $u$ in $C_{SP}(\FF,\GG)$. Since $X$ is locally compact and Polish, if $\NN$ is small enough, we can assume that $\bar v=\bar u$ for all $v\in\NN$. So the strong compact-open topology and strong plaquewise topology have the same restrictions to $\NN$, and the corresponding space will be denoted by $\NN_S$. Observe also that the second factor projection, $\operatorname{pr}_2:Y\to L'$, induces a
homeomorphism
$$
\operatorname{pr}_{2*}:\NN_S\to C_S(X,L')\;,\quad v\mapsto\operatorname{pr}_2\circ v=\tilde v\;.
$$
Therefore, by the relative approximation theorem \cite[pp.~48--49]{Hirsch}, there is some
$h\in\operatorname{pr}_{2*}(\NN)$ which is \coinf\ on some neighborhood of $Q$, and which equals $\tilde u$ on $M\sm W$. Then the result follows with the map $v\in\NN$ satisfying $\tilde v=h$.
\eprf

\prf[Proof of \reft{coinf}]
This is similar to the proof of the usual approximation theorem \cite[p.~49]{Hirsch}.
Fix the following:
\begin{itemize}

\item A locally finite open covering $\{U_i\}$ of $X$ by relatively compact
domains of foliated charts of $\FF$.

\item A covering $\{K_i\}$ of $X$ by compact subsets, with the same index set, satisfying $K_i\subset U_i$ for all $i$.

\item An open subset $V_i\subset M$ for each $i$ so that $K_i\subset V_i$ and
$\overline{V_i}\subset U_i$.

\item For each $i$, a simple open set $U'_i$ of $Y$ such that $f(\overline{U_i})\subset U'_i$. Let $p_i:U'_i\to T'_i$ denote the canonical projection to the local quotient.

\item An open neighborhood $O$ of $E$ so that $f$ is \coinf\ around $\overline{O}$.

\item Any neighborhood $\MM$ of $f$ in $C_{SP}(\FF,\GG)$ satisfying
$g(\overline{U_i})\subset U'_i$ for all $i$ and all $g\in\MM$; this
property holds when $\MM$ is small enough.

\end{itemize}
Since the open covering $\{U_i\}$ of $X$ is locally
finite, and $X$ is locally compact and Polish, the above index set is countable; say, either \N, or $\{0,\dots,k\}$ for some
$k\in\N$. 

\cla{coinf}
For each $i$, there is some $g_i\in\MM$ which is \coinf\ on some
neighborhood of $K_0\cup\dots\cup K_i$, which equals $f$ on some neighborhood of
$\overline{O}$, and which equals $g_{i-1}$ on some neighborhood of $X\sm U_i$ if $i>0$.
\ecla

The result follows from this assertion because a foliation map $g\in\MM$ is well
defined by
$g(x)=g_{i_x}(x)$ for each $x\in M$, where 
$$
i_x=\max\{i\ |\ x\in U_i\}\;.
$$ 
Such a $g$ is \coinf\ because the sets $K_i$ cover $M$, and $g$ is obviously equal to
$f$ on $O$.

To prove \refcl{coinf}, each $g_i$ is defined by induction on $i$. To
simplify the construction of this sequence, we can assume $U_0=K_0=\emptyset$, and take
$g_0=f$. Now let $i>0$, and assume that $g_{i-1}$ is defined. Let
$\FF_i=\FF|_{U_i}$ and
$\GG_i=\GG|_{U'_i}$, and let $u_{i-1}$ denote the restriction
$g_{i-1}:\FF_i\to\GG_i$. Let $\YY_i$ denote the
subspace of $C_{SP}(\FF_i,\GG_i)$ whose elements are the foliated maps
$v\in C(\FF_i,\GG_i)$ satisfying 
$$
p_i\circ v=p_i\circ u_{i-1}\;,\quad 
v|_{U_i\sm V_i}=u_{i-1}|_{U_i\sm V_i}\;.
$$
Let $H_i:\YY_i\to C_{SP}(\FF,\GG)$ be the extension map given by 
$$
H_i(v)=
\begin{cases}
v&\text{on $U_i$}\\
g_{i-1}&\text{on $X\sm U_i$}\;.
\end{cases}
$$
It is easy to prove that $H_i$ is continuous, and thus there is some neighborhood $\NN_i$ of
$u_{i-1}$ in $\YY_i$ such that $H_i(\NN_i)\subset\MM$. Let $\NN'_i$ be a neighborhood of
$u_{i-1}$ in $C_{SP}(\FF_i,\GG_i)$ with $\NN_i=\YY_i\cap\NN'_i$.
Then, by \refl{coinf}, there
is some $v_i\in\NN'_i$ which is 
\coinf\ on some neighborhood of $K_i$, and which equals $u_{i-1}$ on some neighborhood of $U_i\sm W_i$, where $W_i=V_i\cap(U_i\sm\overline{O})$. It follows that $v_i\in\YY_i$, and 
\refcl{coinf} holds with $g_i=H_i(v_i)$.
\eprf

\cor{coinf integrable homotopy}
If there is an integrable homotopy between two \coinf\ foliated maps $\FF\to\GG$, then there also exists a \coinf\ integrable homotopy
between them.
\ecor

\prf
Let $H$ be an integrable homotopy between two maps $f,g\in\coinf(\FF,\GG)$. 
Then let $\overline{H}:\FF\times\R\to\GG$ be the foliated map defined by
$$
\overline{H}(x,t)=
\begin{cases}
H(x,t)&\text{if $t\in I$}\\
f(x)&\text{if $t\le0$}\\
g(x)&\text{if $t\ge1$.}\\
\end{cases}
$$ 
This $\overline{H}$ is \coinf\ on $X\times(\R\setminus I)$. Then, by \reft{coinf}, there
is a \coinf\ foliated map $F:\FF\times\R\to\GG$ which equals $\overline{H}$
on $X\times(-\infty,-1]$ and $X\times[2,\infty)$. Thus the foliated map
$F':\FF\times I\to\GG$, defined by 
$F'(x,t)=F(x,3t-1)$, is a \coinf\ integrable homotopy between $f$ and $g$.
\eprf

\cor{coinf}
Any map in $C(\FF,\GG)$ is integrably homotopic to a map in $C^{0,\infty}(\FF,\GG)$.
\ecor

\prf
This is a direct consequence of \reft{integrable homotopy} and \refc{coinf integrable homotopy}.
\eprf

\ssec{proper}{The case of proper maps}
Now, assume that $Y$ is locally compact and Polish too. When $\Proper(\FF,\GG)$ is endowed with the restriction of the strong plaquewise topology, it will be denoted by $\Proper_{SP}(\FF,\GG)$. Let 
$$
\Proper^{0,\infty}(\FF,\GG)=\Proper(\FF,\GG)\cap\coinf(\FF,\GG)\;.
$$

\th{proper coinf}
With the above notation and conditions, in \reft{coinf}, if $f$ is proper, then $g$ can be chosen to be proper too. In particular, $\Proper^{0,\infty}(\FF,\GG)$ is dense in $\Proper_{SP}(\FF,\GG)$.
\eth

\prf
This follows from \refts{SP proper}{coinf}.
\eprf

\cor{proper coinf integrable homotopy}
If there is a proper integrable homotopy between two proper \coinf\ foliated maps $\FF\to\GG$, then there also exists a proper \coinf\ integrable homotopy
between them.
\ecor

\prf
This can be proved like \refc{coinf integrable homotopy}, by using \reft{proper coinf} instead of \reft{coinf}.
\eprf

\cor{proper coinf}
There is a proper integrable homotopy between any proper continuous foliated map $\FF\to\GG$ and some proper \coinf\ foliated map.
\ecor

\prf
This is a direct consequence of \reft{integrable homotopy} and \refc{proper coinf integrable homotopy}.
\eprf

\ssec{cinfo}{The case of \cinfo\ foliated maps}

If \FF\ and \GG\ are \cinfo\ foliations, the notation $C^{\infty,0}_{SP}(\FF,\GG)$ is used when the set $\cinfo(\FF,\GG)$ is endowed with the restriction of the strong plaquewise topology.
If moreover \FF\ and \GG\ are \cinf\ foliations, the regularity of the approximation in \reft{coinf} can be improved as follows.

\th{cinf approx of cinfo}
With the above notation and conditions, suppose that some $f\in\cinfo(\FF,\GG)$ is of class \cinf\ on some  neighborhood of a closed subset $E\subset M$. Then any neighborhood of $f$ in
$C_{SP}(\FF,\GG)$ contains some $g\in\cinf(\FF,\GG)$ such that $g=f$ on some
neighborhood of $E$. In particular, $\cinf(\FF,\GG)$ is dense in $C^{\infty,0}_{SP}(\FF,\GG)$.
\eth

\reft{cinf approx of cinfo} follows by adapting the proof of \reft{coinf} and using the following version of \refl{coinf}.

\lem{coinf}
With the notation of \refl{coinf}, suppose that $T$ and $T'$ are open sets in Euclidean spaces.
Let $Q\subset X$ be a closed subset and $W\subset X$ an open subset. If some $u\in
\cinfo(\FF,\GG)$ is of class \cinf\ on some neighborhood of $Q\sm W$, then any
neighborhood of $u$ in $C_{SP}(\FF,\GG)$ contains some map $v$
which is \cinf\ on some neighborhood of $Q$, and which equals $u$ on $X\sm W$.
\elem

\prf
The condition that $u$ is \cinfo\  means that $\bar u$ is \cinf. On the other hand, the proof of the approximation theorem \cite[Chapter~2, Theorem~2.6]{Hirsch} can be easily adapted to get a relative approximation theorem. Then, with the notation of the proof of \refl{coinf}, we can assume that $h$ is \cinf, and thus $v$ is \cinf.
\eprf

\cor{cinf approx of cinfo integrable homotopy}
With the above notation and conditions, let $f$ and $g$ be \cinfo\ foliated maps $\FF\to\GG$. If there is a
integrable homotopy between $f$ and $g$, then there also exists a \cinf\ integrable homotopy
between $f$ and $g$.
\ecor

\prf
This follows by using \reft{cinf approx of cinfo} instead of \reft{coinf} in the proof \refc{coinf integrable homotopy}.
\eprf

\cor{cinfo is integrably homotopic to cinf}
Any map in $\cinfo(\FF,\GG)$ is integrably homotopic to a map in $\cinf(\FF,\GG)$.
\ecor

\prf
This is a direct consequence of \reft{integrable homotopy} and \refc{cinf approx of cinfo integrable homotopy}.
\eprf

\ssec{proper cinfo}{The case of proper \cinfo\ foliated maps}
In \refss{cinfo}, suppose that $Y$ is also locally compact and Polish.
Let $\Proper^{\infty,0}_{SP}(\FF,\GG)$ denote the set $\Proper^{\infty,0}(\FF,\GG)$ endowed with the restriction of the strong plaquewise topology, and let 
$$
\Proper^\infty(\FF,\GG)=\Proper(\FF,\GG)\cap\cinf(\FF,\GG)\;.
$$

\th{proper cinf approx of cinfo}
With the above notation and conditions, in \reft{cinf approx of cinfo}, if $f$ is proper, then $g$ can be chosen to be proper too. In particular, $\Proper^\infty(\FF,\GG)$ is dense in $\Proper^{\infty,0}_{SP}(\FF,\GG)$.
\eth

\prf
This follows from \refts{SP proper}{cinf approx of cinfo}.
\eprf

\cor{proper cinf approx of cinfo integrable homotopy}
If there is a proper
integrable homotopy between two maps in $\Proper^\infty(\FF,\GG)$, then there also exists a proper \cinf\ integrable homotopy
between them.
\ecor

\prf
This follows by using \reft{proper cinf approx of cinfo} instead of \reft{cinf approx of cinfo} in the proof \refc{cinf approx of cinfo integrable homotopy}.
\eprf

\cor{proper cinfo is integrably homotopic to cinf}
There is a proper integrable homotopy between any proper \cinfo\ foliated map $\FF\to\GG$ and some proper \cinf\ foliated map.
\ecor

\prf
This is a direct consequence of \reft{integrable homotopy} and \refc{proper cinf approx of cinfo integrable homotopy}.
\eprf

\ssec{Riem folns}{The case of Riemannian foliations with dense leaves}

Assume that \FF\ and \GG\ are transversely complete Riemannian foliations, and that the leaves of \GG\ are dense. Then $C(\FF,\GG)=C^{\infty,0}(\FF,\GG)$ by \reft{main}-(iii). So \reft{cinf approx of cinfo} and \refcs{cinf approx of cinfo integrable homotopy}{cinfo is integrably homotopic to cinf} have the following consequences.

\th{cinf approx for Riem}
With the above notation and conditions, suppose that some $f\in C(\FF,\GG)$ is of class \cinf\ on some  neighborhood of a closed subset $E\subset M$. Then any neighborhood of $f$ in
$C_{SP}(\FF,\GG)$ contains some $g\in\cinf(\FF,\GG)$ such that $g=f$ on some
neighborhood of $E$. In particular, $\cinf(\FF,\GG)$ is dense in $C_{SP}(\FF,\GG)$.
\eth

\cor{cinf approx of integrable homotopy for Riem}
If there is a
integrable homotopy between two continuous foliated maps $\FF\to\GG$, then there also exists a \cinf\ integrable homotopy
between them.
\ecor

\cor{integrably homotopic to cinf for Riem}
Any continuous foliated map $\FF\to\GG$ is integrably homotopic to a \cinf\ foliated map.
\ecor

\ssec{proper Riem folns}{The case of proper foliated maps between Riemannian foliations with dense leaves}

\reft{proper cinf approx of cinfo}, and \refcs{proper cinf approx of cinfo integrable homotopy}{proper cinfo is integrably homotopic to cinf} have the following consequences.

\th{proper cinf approx for Riem}
In \reft{cinf approx for Riem}, if $f$ is proper, then $g$ can be chosen to be proper too. In particular, $\Proper^\infty(\FF,\GG)$ is dense in $\Proper_{SP}(\FF,\GG)$.
\eth

\cor{proper cinf approx of integrable homotopy for Riem}
If there is a proper
integrable homotopy between two maps in $\Proper(\FF,\GG)$, then there also exists a proper \cinf\ integrable homotopy
between them.
\ecor

\cor{proper is integrably homotopic to cinf for Riem}
There is a proper integrable homotopy between each proper continuous foliated map $\FF\to\GG$ and some proper \cinf\ foliated map.
\ecor

\sec{SH}{The strong horizontal topology}

In this section, we introduce a second version of the strong and weak compact-open topologies. Now, we consider the set of foliated maps between two transversely complete Riemannian foliations, and the topology is closely related with the horizontal metric. 

Let \FF\ and \GG\ be transversely complete Riemannian foliations on manifolds $M$ and $N$, and fix any transversely complete bundle-like metric on $N$. Consider the following data:
\begin{itemize} 

\item Any $f\in C(\FF,\GG)$.

\item Any locally finite family $\QQ=\{Q_a\}$ of saturated closed subsets of
$M$ with compact projection to $M/\overline{\FF}$. 

\item A family $\EE=\{\epsilon_a\}$ of positive real numbers, with the same index set as \QQ.

\end{itemize}
Let $\MM(f,\QQ,\EE)$ be the set of continuous foliated maps
$g:\FF\to\GG$ such that there is some sequence $(f_0,f_1,\dots)$ in $C(\FF,\GG)$ satisfying the following properties:
\begin{itemize}

\item $f_0=f$;

\item $f_k=g$ on each $Q_a$ for all but finitely many $k\in\N$; and

\item for each $x\in M$ and $k\in\N$, there is some horizontal geodesic arc $c_{x,k}$ between $f_k(x)$ and $f_{k+1}(x)$ so that
$$
x\in Q_a\;\Longrightarrow\;\sum_{k=0}^\infty\length(c_{x,k})<\epsilon_a\;.
$$

\end{itemize}
Suppose that $\epsilon_a\ge\epsilon'_a+\epsilon''_a$ with $\epsilon'_a,\epsilon''_a>0$, and let $\EE'=\{\epsilon'_a\}$ and $\EE''=\{\epsilon''_a\}$. Then
\begin{equation}\label{e:calM}
g\in\MM(f,\QQ,\EE')\;\Longrightarrow\;\MM(g,\QQ,\EE'')\subset\MM(f,\QQ,\EE)\;.
\end{equation}
It follows from~\refe{calM} that the family of such sets $\MM(f,\QQ,\EE)$ form a base of a topology, called the
{\em strong horizontal topology\/}, and the corresponding
space is denoted by $C_{SH}(\FF,\GG)$. A weak version of this topology can be defined by considering finite families; it can be called the {\em weak horizontal topology\/}, and the notation $C_{WH}(\FF,\GG)$ can be used for the corresponding space. Both of these topologies are equal when $M/\overline{\FF}$ is compact; in this case, the term {\em horizontal topology\/} and the notation $C_H(\FF,\GG)$ can be used.

\prop{SH foliated homotopy}
With the above notation and conditions, if two continuous foliated maps $f,g:\FF\to\GG$ are close enough in $C_{SH}(\FF,\GG)$, then there exists a foliated homotopy $G:\FF\times I_{\text{\rm pt}}\to\GG$ between them. If moreover $f$ and $g$ are proper, then $G$ can also be chosen to be proper.
\eprop

\prf
Consider an open set $\MM=\MM(f,\QQ,\EE)$ as above, which is a neighborhood of $f$ in $C_{SH}(\FF,\GG)$. For each $g\in\MM$, take the maps $f_k$ given by the definition of \MM. Suppose that \QQ\ covers $M$. Then, according to \refs{Riem folns}, there is a \cinf\ regular submanifold $\Omega\subset N\times N$ and a foliated homotopy $H:(\GG\times\GG)|_\Omega\times I_{\text{\rm pt}}\to(\GG\times\GG)|_\Omega$ such that:
\begin{itemize}

\item if the numbers $\epsilon_a$ are small enough, then $(f_k,f_{k+1})(M)\subset\Omega$; and

\item $H_0(x,y)=(x,x)$ for all $(x,y)\in\Omega$, and $H_1=\id_\Omega$.

\end{itemize}
Therefore, when the numbers $\epsilon_a$ are small enough, a foliated homotopy $G_k$ between each $f_k$ and $f_{k+1}$ is given by the composite
$$
\begin{CD}
\FF\times I_{\text{\rm pt}} @>{(f_k,f_{k+1})\times\id_I}>> (\GG\times\GG)|_\Omega\times I_{\text{\rm pt}} @>H>> (\GG\times\GG)|_\Omega @>{\pr_2}>> \GG\;,
\end{CD}
$$
where $\pr_2:\GG\times\GG\to\GG$ is the second factor projection. Observe that $G_k(x,t)=f_k(x)$ if $f_k(x)=f_{k+1}(x)$. 

Now, fix any increasing sequence $0=t_0<t_1<\dots\uparrow1$, and let $\rho_k:[t_k,t_{k+1}]\to I$ be any orientation preserving homeomorphism for each $k\in\N$. Then a foliated homotopy $G:\FF\times I_{\text{\rm pt}}\to\GG$ between $f$ and $g$ can be defined by 
$$
G(x,t)=
\begin{cases}
G_k(x,\rho_k(t)) & \text{if $t_k\le t\le t_{k+1}$ for some $k\in\N$}\\
g(x) & \text{if $t=1$.}
\end{cases}
$$

Assume that moreover $\epsilon=\max\EE<\infty$. Observe that the curve $t\mapsto H(x,t)$ is horizontal for each $x\in M$. Hence, if $\pr_1:N\times I\to N$ denotes the first factor projection, we have
$$
f\circ\pr_1(G^{-1}(R))\subset\Pen_H(R,\epsilon)
$$
for any $R\subset N$, yielding 
$$
G^{-1}(R)\subset f^{-1}(\Pen_H(R,\epsilon))\times I\;.
$$
Therefore, by \refp{Pen H}, $G$ is proper if $f$ is proper.
\eprf

\prop{SH proper}
With the above notation and conditions, $\Proper(\FF,\GG)$ is open in $C_{SH}(\FF,\GG)$.
\eprop

\prf
Let $f\in\Proper(\FF,\GG)$. Choose \QQ\ and \EE\ as above such that \QQ\ covers $M$ and $\epsilon=\max\EE<\infty$. Let $\MM=\MM(f,\QQ,\EE)$, which is a neighborhood of $f$ in $C_{SH}(\FF,\GG)$. For any $g\in\MM$ and $R\subset N$, we have
$$
f(g^{-1}(R))\subset\Pen_H(R,\epsilon)
$$
because \QQ\ covers $M$, yielding
$$
g^{-1}(R)\subset f^{-1}(\Pen_H(R,\epsilon))\;. 
$$
So $\MM\subset\Proper(\FF,\GG)$ by \refp{Pen H}, and the result follows.
\eprf

\sec{transv approx}{\cinfo\ approximations of foliated maps}

Approximation of foliated maps by \cinfo\ ones is impossible in general (Example~\ref{ex:there is no cinf approx}). But, in the case of
transversely complete Riemannian foliations, such a \cinfo\ approximation is a consequence of \reft{main} and the structure of neighborhoods of orbit closures. The strong horizontal topology is appropriate to get this approximation.

\th{cinfo}
With the notation and conditions of \refs{SH}, suppose that some $f\in C(\FF,\GG)$ is \cinfo\ on some 
neighborhood of a closed saturated subset $E\subset M$. Then any neighborhood of $f$ in
$C_{SH}(\FF,\GG)$ contains some $g\in\cinfo(\FF,\GG)$ such that $g=f$ on some
neighborhood of $E$. In particular, $\cinfo(\FF,\GG)$ is dense in $C_{SH}(\FF,\GG)$.
\eth

\prf
Let $\mu=\dim N/\overline{\GG}\in\N$. If $\mu=0$, then $N/\overline{\GG}$ is totally disconnected \cite[Chapter~3, Proposition~1.3]{Pears:75}, and thus the leaf closures of \GG\ are open in $N$. It follows that $f$ itself is \cinfo\ by \reft{main}-(iii). Therefore we can assume $\mu>0$.

Let $\MM=\MM(f,\QQ,\EE)$ be a neighborhood of $f$ of the type described in \refs{SH}, with $\QQ=\{Q_a\}$ and $\EE=\{\epsilon_a\}$.

Let $p=\dim\GG$ and $q=\codim\GG$. Consider the nested sequence
$$
\emptyset=N_{-1}\subset N_0\subset\dots\subset N_q=N\;,
$$
where each $N_\ell$ is the union of all leaf closures of \GG\ with dimension $\le p+\ell$.  According to \refs{Riem folns}, every $N_\ell$ is closed in $N$, each $N_\ell\sm N_{\ell-1}$ is open and dense in $N_\ell$, $N_\ell\sm N_{\ell-1}$ is a $\overline{\GG}$-saturated \cinf\ submanifold of $N$, and the restriction of $\overline{\GG}$ to $N_\ell\sm N_{\ell-1}$ is a regular Riemannian foliation. Let $\EE'=\{\epsilon'_a\}$ with $\epsilon'_a=\epsilon_a/(q+1)$, and set $f_{-1}=f$.

\cla{f ell}
For each $\ell\in\{0,\dots,q\}$, there is some $f_\ell\in\MM(f_{\ell-1},\QQ,\EE')$ which is \cinfo\ on some neighborhood of $f_\ell^{-1}(N_\ell)$, and such that $f_\ell=f_{\ell-1}$ on some neighborhood of $E$.
\ecla

By~\refe{calM}, the result follows from \refcl{f ell} with $g=f_q$. 

\refcl{f ell} is proved by induction on $\ell\in\{-1,0,\dots,q\}$. We have already defined $f_{-1}=f$, which is not required to satisfy any condition. Now, suppose that $f_{\ell-1}$ is already constructed for some $\ell\in\{0,\dots,q\}$. Let $U$ be an open neighborhood of $E\cup f^{-1}(N_{\ell-1})$ where $f_{\ell-1}$ is \cinfo. We can assume that $U$ is $\overline{\FF}$-saturated according to \refs{Riem folns}. Let $U_0$ be another $\overline{\FF}$-saturated open neighborhood of $E\cup f^{-1}(N_{\ell-1})$ whose closure is contained in $U$. Consider the notation of \refs{Riem folns} applied to \GG\ and $N$.

\cla{V m, ...}
For each $m\in\Z_+$, there are saturated open subsets $V_m,W_m\subset M\sm\overline{U_0}$, a leaf closure $F_m$ of \GG\ in $N_\ell\sm N_{\ell-1}$, and some $\delta_m>0$ such that:

\begin{itemize}

\item[(i)] $\{V_m\}$ is locally finite and of order $\le\mu$;

\item[(ii)] $\{W_m\}$ covers $f_\ell^{-1}(N_\ell\sm N_{\ell-1})\sm U$;

\item[(iii)] $\overline{W_m}\subset V_m$;

\item[(iv)] $\widetilde{\Pen}(F_m,\delta_m)\subset\widetilde{\Omega}$;

\item[(v)] $\exp:\widetilde{\FF}_{F_m,\delta_m}\to\FF_{F_m,\delta_m}$ is a foliated diffeomorphism;

\item[(vi)] $f_{\ell-1}(V_m)\subset\Pen(F_m,2^{-\mu}\delta_m)$; and

\item[(vii)] $\delta_m\le\bar\epsilon_m/\mu$, where
$$
\bar\epsilon_m=\min\{\epsilon'_a\ |\ V_m\cap Q_a\neq\emptyset\}\;.
$$

\end{itemize}
\ecla

Since $\overline{\GG}$ is a regular Riemannian foliation on $N_\ell\sm N_{\ell-1}$, there is a covering $\{O_\alpha\}$ of $N_\ell\sm N_{\ell-1}$ by saturated open subsets, and a collection $\{\gamma_\alpha\}$ of positive real numbers such that $\widetilde{\Pen}(F,\gamma_\alpha)\subset\widetilde{\Omega}$ and $\exp:\widetilde{\FF}_{F,\gamma_\alpha}\to\FF_{F,\gamma_\alpha}$ is a foliated diffeomorphism for each index $\alpha$ and every leaf closure $F$ of \GG\ in $O_\alpha$. Let $\{V'_k\}$ ($k\in\Z_+$) be a locally finite family of saturated open subsets of $M\sm U_0$ that covers the closed set $f_{\ell-1}^{-1}(N_\ell)\sm U$, and such that $V'_k\cap(N_\ell\sm N_{\ell-1})\subset f_{\ell-1}^{-1}(O_{\alpha_k})$ for some mapping $k\mapsto\alpha_k$. Then
$$
\bar\epsilon'_k=\min\{\epsilon'_a\ |\ V'_k\cap Q_a\neq\emptyset\}>0
$$
for all $k$, and let $\delta'_k=\min\{\gamma_{\alpha_k},\bar\epsilon'_k/\mu\}$. When $k\in\Z_+$ and $F$ runs in the family of leaf closures of \GG\ in $O_{\alpha_k}$, the family $\{f_{\ell-1}^{-1}(\Pen(F,2^{-\mu}\delta'_k))\}$ is another covering of $f_{\ell-1}^{-1}(N_\ell)\sm U$ by saturated open subsets of $M\sm U_0$. From \cite[Chapter~3, Theorem~4.3]{Pears:75}, it follows that there is a common locally finite refinement $\{V_m\}$ ($m\in\Z_+$) of the coverings $\{V'_k\}$ and $\{f_{\ell-1}^{-1}(\Pen(F,2^{-\mu}\delta'_k))\}$ by saturated open subsets of $M\sm U_0$, and whose order is $\le\mu$. Thus there are mappings $m\mapsto k_m$ and $m\mapsto F_m$, where $F_m$ is a leaf closure of \GG\ in $O_{\alpha_{k_m}}$, so that 
\begin{equation}\label{e:V m, ...}
V_m\subset V'_{k_m}\;,\qquad f_{\ell-1}(V_m)\subset\Pen(F_m,2^{-\mu}\delta'_{k_m})\;.
\end{equation}
Thus the sets $V_m$ satisfy \refcl{V m, ...}-(i), and \refcl{V m, ...}-(vi) follows from the second part of~\refe{V m, ...} with $\delta_m=\delta'_{k_m}$. 
By the first part of~\refe{V m, ...}, we get $\bar\epsilon'_{k_m}\le\bar\epsilon_m$ with the definition of $\bar\epsilon_m$ given in \refcl{V m, ...}-(vii). So 
$$
\delta_m\le\min\{\gamma_{\alpha_{k_m}},\bar\epsilon_m/\mu\}\;,
$$
yielding the conditions~(iv),~(v) and~(vii) of \refcl{V m, ...}. A shrinking of the covering $\{V_m\}$ gives the family $\{W_m\}$ satisfying conditions~(ii) and~(iii) of \refcl{V m, ...}, finishing the proof of this assertion.

According to \refs{Riem folns}, from \refcl{V m, ...}, we get the following for each $m\in\Z_+$:
\begin{itemize}

\item There is a \cinf\ function $\lambda_m:M\to I$, constant on every leaf of \FF, such that $\lambda_m\equiv1$ on some saturated open neighborhood $O_m$ of $M\sm V_m$ and $\lambda_m\equiv0$ on $\overline{W_m}$.

\item There is a \cinf\ foliated homotopy $H_m:\GG_{F_m,\delta_m}\times I_{\text{\rm pt}}\to\GG_{F_m,\delta_m}$ such that:
\begin{itemize}

\item $\im H_{m,0}\subset F_m$;

\item $H_{m,1}=\id_{\GG_{F_m,\delta_m}}$; and

\item the mapping $t\mapsto H_m(x',t)$ is a horizontal geodesic segment of length $<d(x',F_m)$ for each $x'\in\Pen(F_m,\delta_m)$, where $d$ denotes the distance function of $N$.

\end{itemize}

\end{itemize}

Let $\mu_0=0$ and, for $k\in\Z_+$, let $\mu_k$ be the order of the family $\{V_1,\dots,V_k\}$; thus $\mu_0\le\mu_1\le\cdots$ and $\mu_k\le\mu$ for all $k\in\N$. 

\cla{f ell,k}
For each $k\in\N$, there is some $f_{\ell,k}\in C(\FF,\GG)$ such that:
\begin{itemize}

\item[(i)] $f_{\ell,k}=f_{\ell-1}$ on $M\sm V_k$;

\item[(ii)] $f_{\ell,k}$ is \cinfo\ on $W_1\cup\dots\cup W_k\cup U$ for $k>0$;

\item[(iii)] $f_{\ell,k}(V_m)\subset\Pen(F_m,2^{\mu_k-\mu}\delta_m)$ for all $k$ and $m$; and,

\item[(iv)] for each $x\in V_m\cap V_k$, there is a horizontal geodesic segment $c_{x,k}$ between $f_{\ell,k-1}(x)$ and $f_{\ell,k}(x)$ whose length is $<2^{\mu_{k-1}-\mu}\delta_m$.

\end{itemize}
\ecla

First, set $f_{\ell,0}=f_{\ell-1}$, which satisfies the conditions of \refcl{f ell,k} because $\mu_0=0$. Now assume that $f_{\ell,k-1}$ is defined for some $k>0$ and satisfies the conditions of \refcl{f ell,k}. Then let
$f_{\ell,k}$ be the combination of the restriction $f_{\ell,k-1}:\FF|_{O_k}\to\GG$ and the composite
$$
\begin{CD}
\FF|_{V_k} @>{(f_{\ell,k-1},\lambda_k)}>> \GG_{F_k,\delta_k}\times I_{\text{\rm pt}} @>{H_k}>> \GG_{F_k,\delta_k} \hookrightarrow\GG\;.
\end{CD}
$$
Observe that $f_{\ell,k}$ is well defined because $f_{\ell,k-1}$ satisfies \refcl{f ell,k}-(iii). Moreover $f_{\ell,k}$ is a continuous foliated map $\FF\to\GG$ because $f_{\ell,k-1}$ and $H_k$ are continuous foliated maps, and $\lambda_k$ is constant on the leaves. 

We have $f_{\ell,k}=f_{\ell-1}$ on $M\sm V_k$, yielding \refcl{f ell,k}-(i) for $f_{\ell,k}$ because $f_{\ell,k-1}$ satisfies this condition. 

On the one hand, $f_{\ell,k}$ is \cinfo\ on any open set where $f_{\ell,k-1}$ is \cinfo\ because $H_k$ and $\lambda_k$ are \cinf. On the other hand, $f_{\ell,k}(W_k)\subset F_k$, and thus $f_{\ell,k}$ is \cinfo\ on $W_k$ by \reft{main}-(iii). Therefore $f_{\ell,k}$ satisfies \refcl{f ell,k}-(ii) since so does $f_{\ell,k-1}$. 

Take any $m\in\Z_+$ and any $x\in V_m$. If $x\in V_m\sm V_k$, then $f_{\ell,k}(x)=f_{\ell,k-1}(x)$, yielding
$$
d(f_{\ell,k}(x),F_m)=d(f_{\ell,k-1}(x),F_m)<2^{\mu_{k-1}-\mu}\delta_m
\le2^{\mu_k-\mu}\delta_m\;.
$$
If $x\in V_m\cap V_k$, then $\mu_k=\mu_{k-1}+1$, and the mapping 
$$
t\mapsto H_k(f_{\ell,k-1}(x),1-t+t\,\lambda_k(x))\;,\quad0\le t\le1\;,
$$
is a horizontal geodesic segment $c_{x,k}$ between $f_{\ell,k-1}(x)$ and $f_{\ell,k}(x)$ whose length is 
$<2^{\mu_{k-1}-\mu}\delta_m$ since $f_{\ell,k-1}$ satisfies \refcl{f ell,k}-(iii). Then \refcl{f ell,k}-(iv) follows, and moreover
\begin{align*}
d(f_{\ell,k}(x),F_m)&\le d(f_{\ell,k}(x),f_{\ell,k-1}(x))+d(f_{\ell,k-1}(x),F_m)\\
&<2^{\mu_{k-1}-\mu}\delta_m+2^{\mu_{k-1}-\mu}\delta_m\\
&=2^{\mu_k-\mu}\delta_m
\end{align*}
because \refcl{f ell,k}-(iii) holds for $f_{\ell,k-1}$. Therefore $f_{\ell,k}$ also satisfies \refcl{f ell,k}-(iii), completing the proof of \refcl{f ell,k}.

Now, let $f_\ell:M\to N$ be the map defined by $f_\ell(x)=f_{\ell,k_x}(x)$, where
$$
k_x=
\begin{cases}
\max\{k\in\Z_+\ |\ x\in V_k\} & \text{if $x\in\bigcup_mV_m$}\\
0 & \text{if $x\not\in\bigcup_mV_m$.} 
\end{cases}
$$
Observe that  $f_\ell=f_{\ell,k}$ on some saturated neighborhood of each point if $k$ is large enough because $\{V_m\}$ is of finite order. So the following properties hold: 
\begin{itemize}

\item $f_\ell=f_{\ell-1}$ on $U_0$ by \refcl{f ell,k}-(i) since $U_0\cap\bigcup_mV_m=\emptyset$;

\item $f_\ell$ is a continuous foliated map $\FF\to\GG$ because each $f_{\ell,k}$ is a continuous foliated map $\FF\to\GG$; and

\item $f_\ell$ is \cinfo\ on the neighborhood $U\cup\bigcup_mW_m$ of $f_\ell^{-1}(N_\ell)$ by \refcl{f ell,k}-(ii).

\end{itemize}
Take any $x\in M$. Consider the horizontal geodesic segments $c_{x,k}$ given by \refcl{f ell,k}-(iv) if $x\in V_k$, and let $c_{x,k}$ be the constant geodesic segment at $f_{\ell,k-1}(x)$ if $x\in M\sm V_k$. . Then
$$
\sum_{k=1}^\infty\length(c_{x,k})=0
$$
if $x\not\in\bigcup_mV_m$, whilst, by \refcl{V m, ...}-(vii),
$$
\sum_{k=1}^\infty\length(c_{x,k})<\sum_{x\in V_k}2^{\mu_{k-1}-\mu}\delta_m
\le\mu\delta_m\le\bar\epsilon_m
$$
if $x\in V_m$ for some $m$. So
$$
x\in Q_a\;\Longrightarrow\;\sum_{k=1}^\infty\length(c_{x,k})<\epsilon'_a\;,
$$
which means that $f_\ell\in\MM(f_{\ell-1},\QQ,\EE')$, and \refcl{f ell} follows.
\eprf

\sec{SA}{The strong adapted topology}

With the notation and conditions of the above section, the intersection of the strong plaquewise topology and the strong horizontal topology on $C(\FF,\GG)$ is called the {\em strong adapted topology\/}, and the corresponding
space is denoted by $C_{SA}(\FF,\GG)$. The {\em weak adapted topology\/} can be also defined as the intersection of the weak plaquewise topology and the weak horizontal topology, and the corresponding space can be denoted by $C_{WA}(\FF,\GG)$. Both of these topologies are equal when $M$ is compact; in this case, the term {\em adapted topology\/} and the notation $C_A(\FF,\GG)$ can be used. If \FF\ has compact leaf closures, then the strong adapted topology equals the strong compact-open topology on $C(\FF,\GG)$. Given any $f\in C(\FF,\GG)$, let $\QQ$ and $\EE$ be like in \refs{SH}, and let $\NN$ be an open neighborhood of $f$ in $C_{SP}(\FF,\GG)$. Then, with the notation of \refs{SH}, let 
$$
\MM(f,\QQ,\EE,\NN)=\bigcup_{g\in\NN}\MM(g,\QQ,\EE)\;.
$$
The following result is easily verified.

\lem{SA} 
The sets $\MM(f,\QQ,\EE,\NN)$ form a base of the strong adapted topology. 
\elem

\th{SA foliated homotopy}
With the above notation and conditions, if two foliated maps are close enough in $C_{SA}(\FF,\GG)$, then there exists a foliated homotopy between them.
\eth

\prf
This follows from \reft{integrable homotopy}, \refp{SH foliated homotopy} and \refl{SA}.
\eprf

\cor{SA-close foliated maps}
If two foliated maps $f,g:\FF\to\GG$ are close enough in
$C_{SA}(\FF,\GG)$, then $\Hol(f)$ is homotopic to $\Hol(g)$.
\ecor

\prf
This is a consequence of \refp{foliated homotopy} and \reft{SA foliated homotopy}.
\eprf

\th{SA proper}
With the above notation and conditions, $\Proper(\FF,\GG)$ is open in $C_{SA}(\FF,\GG)$.
\eth

\prf
This is a consequence of \reft{SP proper} and \refp{SH proper}.
\eprf

\sec{cinf approx}{\cinf\ approximations of foliated maps}

\ssec{general cinf approx}{The general case}
Combining the \coinf\ and \cinfo\ approximations, we get \cinf\ approximation of foliated maps for transversely complete Riemannian foliations.

\th{cinf Riem foln}
With the notation and conditions of \refs{SA}, suppose that some $f\in C(\FF,\GG)$ is of class \cinf\ on some  neighborhood of a closed saturated subset $E\subset M$. Then any neighborhood of $f$ in
$C_{SA}(\FF,\GG)$ contains some $g\in\cinf(\FF,\GG)$ such that $g=f$ on some
neighborhood of $E$. In particular, $\cinf(\FF,\GG)$ is dense in $C_{SA}(\FF,\GG)$.
\eth

\prf
This is a consequence of \refts{coinf}{cinfo}.
\eprf

\cor{cinf foliated homotopy}
If there is a foliated homotopy between two \cinf\ foliated maps $\FF\to\GG$, then there also exists a \cinf\ foliated homotopy between them.
\ecor

\prf
This follows with the arguments of \refc{coinf integrable homotopy}, by using foliated homotopies instead of integrable homotopies, and using \reft{cinf Riem foln} instead of \reft{coinf}.
\eprf

\cor{homotopic cinf Riem foln}
Any foliated map $\FF\to\GG$ is foliatedly homotopic to a \cinf\ foliated map.
\ecor

\prf
This is a direct consequence of \reft{SA foliated homotopy} and \reft{cinf Riem foln}.
\eprf

\ssec{proper cinf approx}{The case of proper foliated maps}
The notation $\Proper_{SA}(\FF,\GG)$ will be used when the set $\Proper(\FF,\GG)$ is endowed with the restriction of the strong adapted topology.

\th{proper cinf approx Riem}
In \reft{cinf Riem foln}, if $f$ is proper, then $g$ can be chosen to be proper too. In particular, $\Proper^\infty(\FF,\GG)$ is dense in $\Proper^{\infty}_{SA}(\FF,\GG)$.
\eth

\prf
This follows from \refts{SA proper}{cinf Riem foln}.
\eprf

\cor{proper cinf approx of foliated homotopy}
If there is a proper foliated homotopy between two maps in $\Proper^\infty(\FF,\GG)$, then there also exists a proper \cinf\ foliated homotopy between them.
\ecor

\prf
This can be proved with the arguments of \refc{coinf integrable homotopy}, by using proper foliated homotopies instead of integrable homotopies, and using \reft{proper cinf approx Riem} instead of \reft{coinf}.
\eprf

\cor{proper cinfo is integrably homotopic to cinf}
There is a proper foliated homotopy between any proper continuous foliated map $\FF\to\GG$ and some proper \cinf\ foliated map.
\ecor

\prf
This is a direct consequence of \refts{SA foliated homotopy}{proper cinf approx Riem}.
\eprf

\sec{spectral seq}{The spectral sequence of a \cinf\ foliation}

Let \FF\ be a \cinf\ foliation of dimension $p$ and codimension $q$ on a \cinf\ manifold $M$. Consider the decreasing filtration of the de~Rham differential algebra $(\Omega=\Omega(M),d)$ by the differential ideals
$$
\Omega=F^0\Omega\supset F^1\Omega\supset\dots\supset F^q\Omega\supset F^{q+1}\Omega=0\;,
$$
where and $r$-form is in $F^k\Omega$ if it vanishes when $r-k+1$ vectors are tangent to the leaves; intuitively, this means that its ``transverse degree'' is $\ge k$. The induced spectral sequence $(E_i=E_i(\FF),d_i)$ is a differentiable invariant of \FF\ \cite{Masa:85Lie}, \cite{Sarkaria:78finite}, \cite{KamberTondeur:83Fm}, \cite{Alv:89finiteness} (see also \cite{McCleary} for the general theory of spectral sequences). 

A compactly supported version $(E_{c,i}=E_{c,i}(\FF),d_i)$ can be also defined by restricting the above filtration to the subcomplex $\Omega_c\subset\Omega$ of compactly supported differential forms. 

Let \GG\ be another \cinf\ foliation on a manifold $N$. For each $f\in\cinf(\FF,\GG)$, the corresponding homomorphism $f^*:\Omega(N)\to\Omega(M)$ preserves the filtrations and induces a spectral sequence homomorphism $E_i(f):(E_i(\GG),d_i)\to(E_i(\FF),d_i)$. 

If $f$ is a proper map, then $f^*$ restricts to a homomorphism $\Omega_c(N)\to\Omega_c(M)$, inducing a spectral sequence homomorphism $E_{c,i}(f):(E_{c,i}(\GG),d_i)\to(E_{c,i}(\FF),d_i)$.

\prop{Ei}
For \cinf\ foliated maps $f,g:\FF\to\GG$, we have the following: 
\begin{itemize}

\item[(i)] If there is a \cinf\ integrable homotopy between $f$ and $g$, then $E_i(f)=E_i(g)$ for $i\ge1$. 

\item[(ii)] If there is a \cinf\ foliated homotopy between $f$ and $g$, then $E_i(f)=E_i(g)$ for $i\ge2$. 

\end{itemize}
\eprop

\prf
This follows easily from the following observations with a general spectral sequence argument. The operator on differential forms induced by a \cinf\ integrable homotopy, as defined e.g. in \cite{BottTu}, preserves the filtration, whilst, for a \cinf\ foliated homotopy, the corresponding operator decreases the filtration degree at most by one.
\eprf

The following result has a similar proof.

\prop{Ec,i}
For \cinf\ proper foliated maps $f,g:\FF\to\GG$, we have the following: 
\begin{itemize}

\item[(i)] If there is a \cinf\ proper integrable homotopy between $f$ and $g$, then $E_{c,i}(f)=E_{c,i}(g)$ for $i\ge1$. 

\item[(ii)] If there is a \cinf\ proper foliated homotopy between $f$ and $g$, then $E_{c,i}(f)=E_{c,i}(g)$ for $i\ge2$. 

\end{itemize}
\eprop

\sec{invariance}{Invariance of the spectral sequence}

Let \FF\ and \GG\ be transversely complete Riemannian foliations. 

\ssec{general spect seq}{The general case}

According to \reft{cinf Riem foln}, any $f\in C(\FF,\GG)$ is foliatedly homotopic to some $g\in\cinf(\FF,\GG)$. Moreover, by \refc{cinf foliated homotopy} and \refp{Ei}-(ii), the homomorphism $E_i(g)$ is independent of the choice of $g$ for
$i\ge2$, and can be denoted by $E_i(f)$. 

\th{invariance}
With the above notation and conditions, for $f,g\in C(\FF,\GG)$, we have the following:
\begin{itemize}

\item[(i)] If $f$ and $g$ are foliatedly homotopic, then $E_i(f)=E_i(g)$ for $i\ge2$.

\item[(ii)] If $f$ is a foliated homotopy equivalence, then $E_i(f)$ is an isomorphism for $i\ge2$.

\end{itemize}
\eth

\prf
This is a consequence of \refc{cinf foliated homotopy} and \refp{Ei}-(ii).
\eprf

\ssec{proper general spect seq}{The case of proper foliated maps}
According to \refc{proper cinfo is integrably homotopic to cinf}, there is a proper foliated homotopy between each $f\in\Proper(\FF,\GG)$ and some $g\in\cinf(\FF,\GG)$. Moreover, by \refc{proper cinf approx of foliated homotopy} and \refp{Ec,i}-(ii), the homomorphism $E_{c,i}(g)$ is independent of the choice of $g$ for $i\ge2$, and can be denoted by $E_{c,i}(f)$.

\th{proper invariance}
With the above notation and conditions, for $f,g\in\Proper(\FF,\GG)$, we have the following:
\begin{itemize}

\item[(i)] If there is a proper foliated homotopy between $f$ and $g$, then $E_{c,i}(f)=E_{c,i}(g)$ for $i\ge2$.

\item[(ii)] If $f$ is a proper foliated homotopy equivalence, then $E_{c,i}(f)$ is an isomorphism for $i\ge2$.

\end{itemize}
\eth

\prf
This is a consequence of \refc{proper cinf approx of foliated homotopy} and \refp{Ec,i}-(ii).
\eprf

\ssec{dense spect seq}{The case with dense leaves}

Now assume that the leaves of \GG\ are dense. According to \refc{integrably homotopic to cinf for Riem}, any $f\in C(\FF,\GG)$ is integrably homotopic to some $g\in\cinf(\FF,\GG)$. Moreover, by \refc{cinf approx of integrable homotopy for Riem} and \refp{Ei}-(i), the homomorphism $E_1(g)$ is also independent of the choice of $g$ and can be denoted by $E_1(f)$. 

\th{invariance dense}
With the above notation and conditions, for $f,g\in C(\FF,\GG)$, we have the following:
\begin{itemize}

\item[(i)] If $f$ and $g$ are integrably homotopic, then $E_1(f)=E_1(g)$.

\item[(ii)] If the leaves of \FF\ are dense too and $f$ is an integrable homotopy equivalence, then $E_1(f)$ is an isomorphism.

\end{itemize}
\eth

\prf
This is a consequence of \refc{cinf approx of integrable homotopy for Riem} and \refp{Ei}-(i).
\eprf

\ssec{proper dense spect seq}{The case of proper foliated maps and dense leaves}
By \refc{proper is integrably homotopic to cinf for Riem}, there is a proper integrable homotopy between any $f\in\Proper(\FF,\GG)$ and some $g\in\Proper^\infty(\FF,\GG)$. Moreover, by \refc{proper cinf approx of integrable homotopy for Riem} and \refp{Ec,i}-(i), the homomorphism $E_{c,1}(g)$ is also independent of the choice of $g$ and can be denoted by $E_{c,1}(f)$.

\th{proper invariance dense}
With the above notation and conditions, for $f,g\in\Proper(\FF,\GG)$, we have the following:
\begin{itemize}

\item[(i)] If there is a proper integrable homotopy between $f$ and $g$, then $E_{c,1}(f)=E_{c,1}(g)$.

\item[(ii)] If the leaves of \FF\ are dense too and $f$ is a proper integrable homotopy equivalence, then $E_{c,1}(f)$ is an isomorphism.

\end{itemize}
\eth

\prf
This is a consequence of \refc{proper cinf approx of integrable homotopy for Riem} and \refp{Ec,i}-(i).
\eprf

\sec{S for morphisms}{Strong compact-open topology for morphisms}

\begin{notation}
For spaces $X$ and $Y$, the notation $C_{\text{\rm c-o}}(X,Y)$ or $C_S(X,Y)$ will be used when $C(X,Y)$ is endowed with the compact-open or strong compact-open topology, respectively.
\end{notation}

Let  \HH\ and $\HH'$ be pseudogroups of local transformations of spaces $T$ and $T'$, respectively. 
Define the {\em strong compact-open topology\/} on $C(\HH,\HH')$ to be the initial topology induced by the orbit map $C(\HH,\HH')\to C_S(\HH\backslash T,\HH'\backslash T')$, $\Phi\mapsto\Phi_{\text{orb}}$; the notation $C_S(\HH,\HH')$ will be used for the corresponding space. The {\em compact-open topology\/} on $C(\HH,\HH')$ is defined similarly by using $C_{\text{\rm c-o}}(\HH\backslash T,\HH'\backslash T')$, and the notation $C_{\text{\rm c-o}}(\HH,\HH')$ will be used in this case. These topologies can be described as follows. Suppose that the following data are given:
\begin{itemize}

\item A locally finite family $\QQ=\{Q_a\}$ of \HH-invariant closed subsets of $T$ with compact projection to $\HH\backslash T$.

\item A family $\UU'=\{U'_a\}$ of $\HH'$-invariant open subsets of $T'$, with the same index set.

\end{itemize}
Let $\MM(\QQ,\UU')$ be the family of morphisms $\Phi:\HH\to\HH'$ such that $\Phi(Q_a)\subset U'_a$ for all $a$. All such sets $\MM(\QQ,\UU')$ form a base of open sets of $C_S(\HH,\HH')$. If we consider only finite families, we get a base of $C_{\text{\rm c-o}}(\HH,\HH')$.

Observe that $C_{\text{\rm c-o}}(\HH,\HH')$ and $C_S(\HH,\HH')$ may not be Hausdorff; for instance, they have the trivial topology when $\HH'$ has dense orbits.

From now on, assume that \HH\ and $\HH'$ are complete Riemannian pseudogroups. Then the orbit closure map $C_S(\HH,\HH')\to C_S(\overline{\HH}\backslash T,\overline{\HH'}\backslash T')$ is an identification. A similar property holds for the compact-open topologies. In this case, another description of these topologies can be given as follows. Fix an $\HH'$-invariant metric on $T'$, and let $d$ be the corresponding distance function. Consider the following data:
\begin{itemize} 

\item Any $\Phi\in C(\HH,\HH')$.

\item Any locally finite family $\QQ=\{Q_a\}$ of saturated closed subsets of
$T$ with compact projection to $\overline{\HH}\backslash T$. 

\item A family $\EE=\{\epsilon_a\}$ of positive real numbers, with the same index set as \QQ.

\end{itemize}
Let $\MM(\Phi,\QQ,\EE)$ be the set of morphisms
$\Psi:\HH\to\HH'$ such that $d(\Phi(F),\Psi(F))<\epsilon_a$ for all index $a$ and every orbit closure $F$ in $Q_a$. Notice that this condition means that, between $\Phi(F)$ and $\Psi(F)$,  there is some geodesic segment on $\HH'$ of length $<\epsilon_a$ (in the sense of \refs{generalization}). Then the sets $\MM(\Phi,\QQ,\EE)$ form a base of open sets of $C_S(\HH,\HH')$. A base of $C_{\text{\rm c-o}}(\HH,\HH')$ is given by considering only finite families. 

\prop{Hol is cont}
Let \FF\ and \GG\ be transversely complete Riemannian foliations, where \GG\ is endowed with a bundle-like metric. Then the holonomy functor defines continuous maps
$$
C_{SA}(\FF,\GG)\to C_S(\Hol(\FF),\Hol(\GG))\;,\quad 
C_{WA}(\FF,\GG)\to C_{\text{\rm c-o}}(\Hol(\FF),\Hol(\GG))\;.
$$
\eprop

\prf
This follows easily from the observation that each horizontal geodesic segment of the bundle-like metric of \GG\ project to geodesic segments of $\Hol(\GG)$ with the same length.
\eprf

\th{S homotopy}
With the above notation and conditions, if two morphisms are close enough in $C_S(\HH,\HH')$, then there exists a homotopy between them.
\eth

\prf
According to \refs{around the orbit closures}, there is an open neighborhood $\Omega'$ of the diagonal in $T'\times T'$, which is invariant by $\HH'\times\HH'$, and there is a homotopy $\Psi:(\HH'\times\HH')|_{\Omega'}\times I\to(\HH'\times\HH')|_{\Omega'}$ such that $\Psi_0$ is generated by the map $(x,y)\mapsto(x,x)$, and $\Psi_1$ is the identity morphism at $(\HH'\times\HH')|_{\Omega'}$. If two morphisms $\Phi_1$ and $\Phi_2$ are close enough in $C_S(\HH,\HH')$, then $\Phi_1(F)\times\Phi_2(F)\subset\Omega'$ for all \HH-orbit closure $F$. Then a homotopy between $\Phi_1$ and $\Phi_2$ is given by the composite
$$
\begin{CD}
\HH\times I @>{(\Phi_1,\Phi_2)\times\id_I}>> (\HH'\times\HH')|_{\Omega'}\times I @>\Psi>> (\HH'\times\HH')|_{\Omega'} @>>> \HH'\;,
\end{CD}
$$
where the right hand side morphism is generated by the restriction $\Omega'\to T'$ of second factor projection $T'\times T'\to T'$. 
\eprf

\cor{S homotopy}
If the orbits of $\HH'$ are dense, then all morphisms $\HH\to\HH'$ are homotopic to each other.
\ecor

\sec{cinf for morphisms}{\cinf\ approximations of morphisms}

The following is a pseudogroup version of \reft{cinf Riem foln}. 

\th{cinf}
Let \HH\ and $\HH'$ be complete Riemannian pseudogroups. Suppose that some $\Phi\in C(\HH,\HH')$ is  \cinf\ on some neighborhood of a closed saturated subset $E\subset T$. Then any neighborhood of $\Phi$ in $C_S(\HH,\HH')$ contains some $\Psi\in\cinf(\HH,\HH')$ such that $\Phi=\Psi$ on some
neighborhood of $E$. In particular, $\cinf(\HH,\HH')$ is dense in $C_S(\HH,\HH')$.
\eth

\prf
This result can be proved with an easy adaptation to pseudogroups of the arguments of \reft{cinfo}. Instead of observations from  \refs{Riem folns}, we have to use the corresponding observations from Sections~\ref{s:Molino} and~\ref{s:around the orbit closures}. Moreover we have to use the morphism versions of some operations with maps (see \refs{morphisms}); e.g., the combination of morphisms is used in the needed version of \refcl{f ell,k}.
\eprf

Like in the foliated setting, \reft{cinf} has the following consequences.

\cor{cinf homotopy}
If there is a homotopy between two \cinf\ morphisms $\HH\to\HH'$, then there also exists a \cinf\ homotopy between them.
\ecor

\cor{cinf}
Any morphism $\HH\to\HH'$ is homotopic to a \cinf\ morphism.
\ecor

\sec{cohoms of psgrs}{Cohomologies of \cinf\ pseudogroups}

Let \HH\ be a \cinf\ pseudogroup acting on a \cinf\ manifold $T$. The {\em invariant complex\/} of \HH\ is the subcomplex $\Omega_{\text{\rm inv}}(\HH)\subset\Omega(T)$ consisting of differential forms that are \HH-{\em invariant\/}; i.e., those forms $\alpha$ satisfying $h^*(\alpha|_{\im h})=\alpha|_{\dom h}$ for all $h\in\HH$. The corresponding cohomology is called the {\em invariant cohomology\/} and denoted by $H_{\text{\rm inv}}(\HH)$.

A \cinf\ morphism $\Phi:\HH\to\HH'$ induces a homomorphism of complexes $\Phi^*:\Omega_{\text{\rm inv}}(\HH')\to\Omega_{\text{\rm inv}}(\HH)$, which is well defined by the condition $(\Phi^*\alpha)|_{\dom\phi}=\phi^*\alpha$ for $\alpha\in\Omega_{\text{\rm inv}}(\HH')$ and $\phi\in\Phi$.  This assignment $\Phi\mapsto\Phi^*$ is functorial; thus isomorphic pseudogroups have isomorphic invariant complexes. We get an induced homomorphism $\Phi^*:H_{\text{\rm inv}}(\HH')\to H_{\text{\rm inv}}(\HH)$.

Let \FF\ be a \cinf\ foliation on a \cinf\ manifold $M$, and $\PP:M\to\Hol(\FF)$ the corresponding canonical morphism. Then $\PP^*:\Omega_{\text{\rm inv}}(\Hol(\FF))\to\Omega(M)$ restricts to an isomorphism $\PP^*:\Omega_{\text{\rm inv}}(\Hol(\FF))\to E_1^{\bullet,0}(\FF)$, yielding the well known isomorphism $H_{\text{\rm inv}}(\Hol(\FF))\cong E_2^{\bullet,0}(\FF)$.

\prop{H inv}
If there is a \cinf\ homotopy between two \cinf\ morphisms $\Phi_0,\Phi_1:\HH\to\HH'$, then $\Phi_0^*=\Phi_1^*:H_{\text{\rm inv}}(\HH')\to H_{\text{\rm inv}}(\HH)$.
\eprop

\prf
This is a direct adaptation of the proof of the corresponding result for manifolds \cite{BottTu} by using morphisms instead of maps.
\eprf

In the sense of \refs{generalization}, we can also consider the {\em de~Rham cohomology\/} $H_{\text{\rm dR}}(\HH)$  of a \cinf\ pseudogroup \HH\ acting on some \cinf\ manifold $T$, which is precisely defined as follows. Let $\bigwedge\TT\HH^*$ be the pseudogroup on $\bigwedge\TT T^*$ generated by the local transformations of the form $\bigwedge\TT h^*:\bigwedge\TT(\im h)^*\to\bigwedge\TT(\dom h)^*$ with $h\in\HH$; notice that $h^*\alpha=\bigwedge\TT h^*\circ\alpha\circ h$ for all $\alpha\in\Omega(\im h)$. For each degree $r$, the restriction of $\bigwedge\TT\HH^*$ to $\bigwedge^r\TT T^*$ will be denoted by $\bigwedge^r\TT\HH^*$. The fiber bundle projection $\pi:\bigwedge\TT T^*\to T$ generates a morphism $\Pi:\bigwedge\TT\HH^*\to\HH$ since $\pi\circ\bigwedge\TT h^*=h^{-1}\circ\pi$ for all $h\in\HH$. A morphism $\Theta:\HH\to\bigwedge\TT\HH^*$ is called a {\em differential form\/} on \HH\ if it is a {\em section\/} of $\Pi$ in the sense that $\Pi\circ\Theta=\id_{\HH}$; it is said that $\Theta$ is of {\em degree\/} $r$ if $\im\Theta\subset\bigwedge^r\TT T^*$. Any differential form $\Theta$ on \HH\ is generated by differential forms on open subsets of $T$: if the domain of some $\theta\in\Theta$ is small enough, then $h=\pi\circ\theta\in\HH$ and $\theta\circ h^{-1}\in\Omega(\im h)\cap\Theta$. The {\em exterior derivative\/} of a \cinf\ differential form $\Theta$ on \HH\ is the \cinf\ differential form $d\Theta$ on \HH\ generated by the exterior derivatives $d\theta$ of those $\theta\in\Theta$ that are differential forms on open subsets of $T$. The vector space $\Omega(\HH)$ of \cinf\ differential forms on \HH\ becomes a complex endowed with the exterior derivative $d$ and the grading defined as above. Then $H_{\text{\rm dR}}(\HH)$ is the cohomology of $(\Omega(\HH),d)$.  There is a canonical isomorphism  $\Omega_{\text{\rm inv}}(\HH)\cong\Omega_{\text{\rm dR}}(\HH)$, which assigns to each $\alpha\in\Omega_{\text{\rm inv}}(\HH)$ the morphism generated by $\alpha$. Thus $H_{\text{\rm inv}}(\HH)\cong H_{\text{\rm dR}}(\HH)$, obtaining a known simpler expression for the de~Rham cohomology of \HH. 

Another complex associated to \HH\ was introduced by A.~Haefliger as follows \cite{Haefliger:80minimal}. If some $\alpha\in\Omega_c(T)$ is supported in the image of some $h\in\HH$, let $h^*\alpha\in\Omega_c(T)$ denote the extension by zero of $h^*(\alpha|_{\im h})$. The {\em Haefliger complex\/} $(\Omega_{\text{\rm Ha}}(\HH),d)$ of \HH\ is the quotient complex of $(\Omega_c(T),d)$ over the subcomplex generated by forms of the type $\alpha-h^*\alpha$, where $h\in\HH$ and $\alpha\in\Omega_c(T)$ with $\supp\alpha\subset\im h$. The corresponding cohomology is called the {\em Haefliger cohomology\/} and denoted by $H_{\text{\rm Ha}}(\HH)$. The Haefliger cohomology is not Hausdorff in general with the topology induced by the \cinf\ topology on $\Omega_c(T)$ \cite{Haefliger:80minimal}; so it makes sense to consider its maximal Hausdorff quotient, which is called the {\em reduced Haefliger cohomology\/} and denoted by $\overline{H}_{\text{\rm Ha}}(\HH)$. 

Exterior product and integration defines a pairing
$$
\Omega_{\text{\rm inv}}^r(\HH)\otimes\Omega_{\text{\rm Ha}}^{n-r}(\HH)\to\R
$$
for each degree $r$, where $n=\dim T$. It induces a pairing
$$
H_{\text{\rm inv}}^r(\HH)\otimes\overline{H}_{\text{\rm Ha}}^{n-r}(\HH)\to\R\;,
$$
which degenerates in general; nevertheless, the pseudogroup version of the arguments of \cite{Alv:89Duality} and \cite{Masa:92Duality} show that this pairing is non-degenerated for complete Riemannian pseudogroups.

Any \'etal\'e morphism $\Phi:\HH\to\HH'$  functorially induces a continuous homomorphism $\Phi_*:\Omega_{\text{\rm Ha}}(\HH)\to\Omega_{\text{\rm Ha}}(\HH')$  \cite[Section~1.2]{Haefliger:80minimal}; thus isomorphic pseudogroups have isomorphic Haefliger complexes.

For any \cinf\ foliation \FF, $\Omega_{\text{\rm Ha}}(\Hol(\FF))$ and $H_{\text{\rm Ha}}(\Hol(\FF))$ are called the {\em transverse complex\/} and {\em transverse cohomology\/} of \FF, and there is a homomorphism $\int_{\FF}:\Omega_c(M)\to\Omega_{\text{\rm Ha}}(\Hol(\FF))$, called {\em integration along the leaves\/}, which induces an isomorphism $\int_{\FF}:E_{c,1}^{\bullet,p}(\FF)\to\Omega_{\text{\rm Ha}}(\Hol(\FF))$ ($p=\dim\FF$) \cite[Section~3]{Haefliger:80minimal}, yielding $E_{c,2}^{\bullet,p}(\FF)\cong H_{\text{\rm Ha}}(\Hol(\FF))$. 

A morphism $\Phi:\HH\to\HH'$ is called a {\em \cinf\ submersion\/} when its maps are \cinf\ submersions. In this case, $\Phi$ induces a continuous open homomorphism $\Phi_*:\Omega_{\text{\rm Ha}}(\HH)\to\Omega_{\text{\rm Ha}}(\HH')$ in the following way. For each $\alpha\in\Omega_c(T)$, there are some $\alpha_1,\dots,\alpha_k\in\Omega_c(T)$ and some $\phi_1,\dots,\phi_k\in\Phi$ such that $\alpha=\alpha_1+\dots+\alpha_k$ and $\supp\alpha_i\subset\dom\phi_i$ for all $i\in\{1,\dots,k\}$. Since each $\phi_i$ is a submersion, the integration along the fibers $\int_{\phi_i}\alpha_i\in\Omega_c(\im\phi_i)$ is defined, whose extension by zero to $T'$ is also denoted by $\int_{\phi_i}\alpha_i\in\Omega_c(T')$. With the obvious generalization of the arguments of \cite[Theorem~3.1]{Haefliger:80minimal}, it follows that the class of $\int_{\phi_1}\alpha_1+\dots+\int_{\phi_k}\alpha_k$ in $\Omega_{\text{\rm Ha}}(\HH')$ is independent of the choices of $\alpha_1,\dots,\alpha_k$ and $\phi_1,\dots,\phi_k$. Then, if $\zeta\in\Omega_{\text{\rm Ha}}(\HH)$ is the element represented by $\alpha$, define $\Phi_*(\zeta)\in\Omega_{\text{\rm Ha}}(\HH')$ to be the element represented by $\int_{\phi_1}\alpha_1+\dots+\int_{\phi_k}\alpha_k$. The arguments of \cite[Theorem~3.1]{Haefliger:80minimal} also show that $\Phi_*$ is continuous and open. The induced continuous homomorphism $H_{\text{\rm Ha}}(\HH)\to H_{\text{\rm Ha}}(\HH')$ is also denoted by $\Phi_*$, or by $H_{\text{\rm Ha}}(\Phi)$. This defines a covariant functor from $\PsGr$ to the category of continuous homomorphisms between topological vector spaces.

As a particular example, for any \cinf\ foliation \FF\ on a \cinf\ manifold $M$, Haefliger's integration along the fibers, $\int_\FF:\Omega_c(M)\to\Omega_{\text{\rm Ha}}(\Hol(\FF))$, is the homomorphism $\PP_*$ induced by the canonical morphism $\PP:M\to\Hol(\FF)$.

\sec{invariance for psgrs}{Invariance of the invariant cohomology}

The following is a version for pseudogroups of \reft{invariance}. 

\th{invariance psgr}
Let $\Phi,\Psi:\HH\to\HH'$ be morphisms between complete Riemannian pseudogroups. We have the following:
\begin{itemize}

\item[(i)] If $\Phi$ and $\Psi$ are homotopic, then $\Phi^*=\Psi^*:H_{\text{\rm inv}}(\HH')\to H_{\text{\rm inv}}(\HH)$.

\item[(ii)] If $\Phi$ is a homotopy equivalence, then $\Phi^*:H_{\text{\rm inv}}(\HH')\to H_{\text{\rm inv}}(\HH)$ is an isomorphism.

\end{itemize}
\eth

\prf
This follows from \refcs{cinf homotopy}{cinf}, and \refp{H inv}.
\eprf

\section{Examples}\label{s:examples}

\begin{exmp}
Theorem~\ref{t:main}-(i) supplies a large class of complete morphisms. Another source of complete morphisms is the following: any pseudogroup generated by a group action is complete, and any equivariant map generates a complete morphism. 
\end{exmp}

\begin{exmp}\label{ex:non-complete psdgrs}
For $\lambda>1$, the mapping $x\mapsto\lambda x$ generates a complete pseudogroup \HH, whose restriction to $U=(-1,1)$ is not complete. This \HH\ is equivalent to $\HH|_U$ since $U$ cuts every \HH-orbit. So completeness is not invariant by pseudogroup equivalences. As pointed out in \cite[Section~1.3]{Haefliger:88Leaf}, the pseudogroup of all local isometries of a Riemannian manifold is not complete in general; for instance, it is not complete in the case of a sphere with a metric which is flat on some part with non-empty interior, and has positive scalar curvature elsewhere. The identity morphism of any of the above non-complete pseudogroups is not complete. Thus, for Riemannian pseudogroups, the relation between completeness and geodesic completeness seems to be weak. 
\end{exmp}

\begin{exmp}\label{ex:product}
A pseudogroup \HH\ acting on a space $T$ is called {\em globally complete\/} when, for all $h\in\HH$ and any $x\in\dom h$, there is some $\tilde h\in\HH$ such that $\dom\tilde h=T$ and $\gamma(\tilde h,x)=\gamma(h,x)$; i.e., $(T,T)$ is a completeness pair for \HH.

Let $\HH'$ be another pseudogroup acting on a space $T'$, and let $\Phi:\HH\to\HH'$ be a morphism. It is said that $\Phi$ is {\em globally complete\/} when there is some $\phi\in\Phi$ such that $\dom\phi=T$, and moreover, for all $h\in\HH$ and every $z\in\dom h$, there is some $\tilde
h\in\HH$ and some $h'\in\HH'$ so that 
$\dom\tilde h=T$, $\gamma(\tilde h,x)=\gamma(h,x)$, $\im\phi\subset\dom h'$, and
$h'\circ\phi=\phi\circ\tilde h$. In this case, \HH\ is globally complete, $\Phi$ is generated by $\phi$, and $(\phi,T;\phi,T)$ is a completeness quadruple of $\Phi$.

Now, consider a family of pseudogroups $\HH_i$, where each $\HH_i$ acts on a space $T_i$. As a straightforward generalization of the finite product of pseudogroups (\refs{morphisms}), we can define the {\em product pseudogroups\/} $\prod_i\HH_i$ acting on $\prod_iT_i$. It is easy to check that $\prod_i\HH_i$ is complete if and only if each $\HH_i$ is complete, and all but finitely many pseudogroups $\HH_i$ are globally complete.

Consider another family of pseudogroups $\HH'_i$, with the same index set, and a family of morphisms $\Phi_i:\HH_i\to\HH'_i$. As a straightforward generalization of the finite product of morphisms (\refs{morphisms}), we can define the {\em product morphism\/} $\prod_i\Phi_i:\prod_i\HH_i\to\prod_i\HH'_i$. This $\prod_i\Phi_i$ is complete if and only if each $\Phi_i$ is complete, and all but finitely many of the morphisms $\Phi_i$ are globally complete.

This allows the construction of concrete examples of non-complete morphisms between complete pseudogroups, which shows that the Riemannian condition cannot be removed in \reft{main}-(i).  For instance, let \HH\ be the globally complete pseudogroup acting on the line \R\ generated by all homeomorphisms $\R\to\R$, and let $\HH'$ be the globally complete pseudogroup acting on circle $S^1$ generated by all homeomorphisms $S^1\to S^1$. The universal covering map $\R\to S^1$ generates a complete morphism $\Phi:\HH\to\HH'$, which is not globally complete.  Then, with $T_i=\R$, $T'_i=S^1$, $\HH_i=\HH$ and $\HH'_i=\HH'$ for all $i\in\N$, the pseudogroups $\prod_i\HH_i$ and $\prod_i\HH'_i$ are complete, but the morphism $\prod_i\Phi_i:\prod_i\HH_i\to\prod_i\HH'_i$ is not complete.
\end{exmp}

\begin{exmp}\label{ex:R2}
This is another example of a non-complete morphism between complete pseudogroups, which furthermore is \cinf. Let \HH\ be the \cinf\ pseudogroup acting on $\R^2$ generated by the group of diffeomorphisms $h:\R^2\to\R^2$ such that there is some $\epsilon>0$ and some $c\in\R$ so that $h(x,y)=(x+c,y)$ for all $(x,y)\in\R\times(-\epsilon,\epsilon)$. Let $\HH'$ be the \cinf\ pseudogroup acting on $\R^2$ generated by the group of diffeomorphisms $\R^2\to\R^2$ that fix the origin and have the same germ at the origin as a rotation. Let $\phi:\R^2\to\R^2$ be the \cinf\ map defined by $\phi(x,y)=(y\cos x,y\sin x)$. This $\phi$ generates a \cinf\ morphism $\Phi:\HH\to\HH'$. It is easy to prove that \HH\ and $\HH'$ are complete, whilest $\Phi$ is not complete: its completeness condition fails around the origin.
\end{exmp}

\begin{exmp}\label{ex:Arnold}
Consider Arnold's example of a diffeomorphism $h:S^1\to S^1$ which is topologically conjugated but not $C^1$ conjugated to a rotation $h':S^1\to S^1$ \cite{Arnold:61}. By suspension, we get examples of homeomorphic \cinf\ foliations with non-isomorphic basic cohomologies \cite{KacimiNicolau:93topological}; so these foliations cannot be diffeomorphic. Thus, if \HH\ and $\HH'$ denote the pseudogroups generated by $h$ and $h'$, respectively, the isomorphism $\Phi:\HH\to\HH'$ generated by any fixed homeomorphism $\phi:S^1\to S^1$ satisfying $\phi\circ h=h'\circ\phi$ is not \cinf. This does not contradict \reft{main}-(iii) and \reft{cinf} because \HH\ is not Riemannian. However \HH\ is equicontinuous \cite[Appendix~E, Section~5]{Molino:88Riemannian}, \cite{AlvCandel:equicont}, showing that \reft{main}-(iii) cannot be generalized to equicontinuous pseudogroups.
\end{exmp}

\begin{exmp}\label{ex:there is no cinf approx}
Let \HH\ be the pseudogroup on \R\ generated by the map $h$ defined by $h(x)=\lambda x$ for some some $\lambda>1$. The map $\phi$ defined by $\phi(x)=|x|$ satisfies $\phi\circ h=h\circ\phi$. So $\phi$ generates a morphism $\Phi:\HH\to\HH$. It is easy to see that $\Phi$ cannot be approximated by any \cinf\ morphism in $C_S(\HH,\HH)$. By suspension, we get a foliated map between \cinf\ foliations with compact space of leaves that cannot be strongly approximated by \cinf\ foliated maps. This shows that Theorems~\ref{t:cinfo},~\ref{t:cinf Riem foln} and~\ref{t:cinf} cannot be generalized to arbitrary \cinf\ foliations.
\end{exmp}
 
\begin{exmp}
Let \HH\ and $\HH'$ be pseudogroups acting on spaces $T$ and $T'$, respectively. A morphism $\Phi:\HH\to\HH'$ is said to be {\em simply complete\/} if all $x\in T$ and $x'\in T'$ have respective open neighborhoods $U$ and $U'$ such that, for any $\phi\in\Phi$ and every $y\in U\cap\dom\phi$ with $\phi(y)\in U'$, there is some $\tilde\phi\in\Phi$ such that $U\subset\dom\tilde\phi$ and $\gamma(\tilde\phi,y)=\gamma(\phi,y)$; in this case, it will be said that $(U,U')$ is a {\em simple completeness pair\/} of $\Phi$. Observe that \HH\ is complete if and only if the identity morphism $\id_{\HH}$ simply complete; thus this is another version of completeness for morphisms different from \refd{complete morphism}. It is easy to see that any morphism generated by a globally defined map is simply complete. The following result shows that the reciprocal statement is also true under certain topological conditions.

\prop{simply complete}
If $\Phi$ is simply complete, $T$ is simply connected and $T'$ is compact, then $\Phi$ is generated by a single map $T\to T'$.
\eprop

\prf
Since $\Phi$ is simply complete, for each $x\in T$, there is a family of simple completeness pairs of $\Phi$, $\{(U_i,U'_i)\}$, such that each $U_i$ contains $x$ and $\{U'_i\}$ covers $T'$. Because $T'$ is compact, $\{U'_i\}$ admits a finite subcovering $\{U'_{i_1},\dots,U'_{i_n}\}$. Then $U=U_{i_1}\cap\dots\cap U_{i_n}$ is an open neighborhood of $x$ and satisfies the following condition:  for any $\phi\in\Phi$ and every $y\in U\cap\dom\phi$, there is some $\tilde\phi\in\Phi$ such that $U\subset\dom\tilde\phi$ and $\gamma(\tilde\phi,y)=\gamma(\phi,y)$. This means that the source map $\gamma(\Phi)\to T$ is a covering map, which admits a global section $\theta:T\to\gamma(\HH)$ because $T$ is simply connected. Then $\Phi$ is generated by the composite of $\theta$ with the target projection $\gamma(\Phi)\to T'$.
\eprf

A morphism $\Phi:\HH\to\HH'$ is called {\em quasi-analytic\/} when, for all $\phi,\psi\in\Phi$ and any connected open subset $U\subset\dom\phi\cap\dom\psi$, we have $\phi=\psi$ on $U$ if $\phi$ and $\psi$ have the same germ at some point of $U$. The following result is an elementary relation between completeness and simple completeness.

\prop{simply complete and complete}
If $T$ and $T'$ are locally connected, and $\Phi$ is quasi-analytic and simply complete, then $\Phi$ is complete.
\eprop
\end{exmp}

\begin{exmp}\label{ex:Haefliger}
In~\cite{Haefliger:84classifiants}, Haefliger has introduced certain morphisms between topological groupoids in the following way.  Let $\Gamma$ be a topological groupoid with space of units $T$. First, the concept of right action of $\Gamma$ on a space $E$ relative to a map $E\to T$ is defined as an obvious generalization of the case of topological groups. Left actions of $\Gamma$ are defined similarly, as well as $\Gamma$-principal bundles. Let $\Gamma'$ be another topological groupoid with space of units $T'$. Let $E$ be a space endowed with the following data: continuous maps $p:E\to T$ and $p':E\to T'$, a left action of $\Gamma$ on $E$ relative to $p$, and a right action of $\Gamma'$ on $E$ relative to $p'$. Suppose that these two actions commute. If $E$ becomes a $\Gamma'$-principal bundle with projection $p$, then $E$ is called a {\em morphism\/} of $\Gamma$ to $\Gamma'$. Consider also the induced right $\Gamma$-action and left $\Gamma'$-action on $E$ defined by inverting the elements of $\Gamma$ and $\Gamma'$. If moreover $E$ becomes a $\Gamma$-principal bundle with projection $p'$, then $E$ is called an {\em equivalence\/}. These morphisms form a category with the operation defined by taking fiber products in an obvious way. The isomorphisms of this category are the equivalences.

Recall that a topological groupoid is called {\em \'etal\'e\/} when its source and target projections are local homeomorphisms. A topological groupoid is \'etal\'e just when it is isomorphic to the groupoid of germs of some pseudogroup with the \'etal\'e topology. Thus Haefliger morphisms between \'etal\'e groupoids can be also considered as another type of morphisms between pseudogroups, which also form a category denoted by $\PsGr_{\text{\rm Ha}}$. Let $E$ be a Haefliger morphism of a pseudogroup $\HH$ acting on $T$ to a pseudogroup $\HH'$ acting on $T'$, with maps $p:E\to T$ and $p':E\to T'$. Then the composites $p'\circ\sigma$, for local sections $\sigma$ of $p$, generate a morphism $\Phi:\HH\to\HH'$ of our type. If $\Phi$ consists of open maps, then there is a unique Haefliger morphism $E$ inducing $\Phi$. The assignment $E\mapsto\Phi$ defines a functor $\PsGr_{\text{\rm Ha}}\to\PsGr$.

Haefliger's definition of equivalence between topological groupoids  is a version of the concept of Morita equivalence. Thus Heafliger morphisms could be more appropriate to establish relations with the $C^*$-algebras of topological groupoids and their $K$-theory.
\end{exmp}

\section{Open problems}\label{s:problems}

\begin{prob}\label{pb:Hurewicz}
Develop the algebraic topology, differential topology and differential geometry of pseudogroups in the sense of \refs{generalization}. To begin with, we may ask whether results like the Hurewicz isomorphism
theorem \cite{Spanier:81} or the Van~Kampen theorem \cite{Hatcher:02} can be generalized to pseudogroups. Or what about any version of the Hopf-Rinow theorem for Riemannian pseudogroups?
\end{prob}

\begin{prob}
Problem~\ref{pb:Hurewicz} can be also considered for Haefliger morphisms between pseudogroups, or between topological groupoids (see Example~\ref{ex:Haefliger}). In this way, another definition of algebraic invariants can be given for pseudogroups and topological groupoids. The functor $\PsGr_{\text{\rm Ha}}\to\PsGr$ induces natural homomorphisms between Haefliger invariants and our invariants. Are these homomorphisms surjective? Which kind of information is contained in their kernels?

Haefliger has indicated us the following concrete example. Let $n$ be a positive integer, and let $\HH$ be the pseudogroup on $\R^2$ generated by a rotation of order $n$. Then $\HH$ is contractible with our definition, whilst, in $\PsGr_{\text{\rm Ha}}$, there are $n$ Haefliger morphisms of the circle to $\HH$ which are not homotopic to each other.
\end{prob}

\begin{prob}\label{pb:lift}
Let \FF\ and \GG\ be foliated structures, $f:\FF\to\GG$ a continuous foliated map, and $\Psi:\Hol(\FF)\times I\to\Hol(\GG)$ a homotopy such that $\Psi_0\equiv\Hol(f)$. What conditions are needed for the existence of a lift $H:\FF\times I_{\text{\rm pt}}\to\GG$ of $\Psi$ so that $H_0=f$? The needed conditions must be very restrictive (see \cite{Meigniez:02Submersions} for the solution of a related problem). It is easy to produce counterexamples to the existence of such a lift by using vanishing cycles.
\end{prob}

\begin{prob}\label{pb:measurable morphisms}
In the case of pseudogroups generated by the left translations on Lie groups, Theorem~\ref{t:main}-(iii) asserts that any continuous homomorphism between Lie groups is \cinf. But indeed it is well known that any measurable homomorphism between Lie groups is \cinf. Then, according to the generalization of the theorem of Myers-Steenrod for complete Riemannian pseudogroups given in \cite{Salem:88}, there should be a version of Theorem~\ref{t:main}-(iii) for ``measurable morphisms''.
\end{prob}

\begin{prob}
\reft{main}-(iii) cannot be generalized to equicontinuous pseudogroups (Example~\ref{ex:Arnold}), but what about its assertions~(i) and~(ii)? Observe that \reft{main}-(ii) would make sense for equicontinuous pseudogroups because they also have a closure under mild topological conditions \cite{AlvCandel:equicont}, \cite{Kellum}, \cite{Tarquini04:conformes}. Moreover equicontinuity is a key property in the topological description of Riemannian foliations with dense leaves given in \cite{AlvCandel:top}.
\end{prob}

\begin{prob}
Is there any version of completeness for Haefliger morphisms? And a version of \reft{main}?
\end{prob}

\begin{prob}
Is it possible to adapt the convolution technique \cite{Hirsch} to improve differentiability of foliated maps between transversely complete Riemannian foliations? This is clearly possible in the case of transversely parallelizable foliations. 
\end{prob}

\begin{prob}
The spectral sequence can be given by the differential sheaf of local basic differential forms \cite{Masa:02A-S}, which leads to the following question. For transversely complete Riemannian foliations, is it possible to prove the homotopy invariance of the spectral sequence directly from \reft{main}-(iii) with a sheaf theoretic argument?
\end{prob}

\begin{prob}
We may ask whether a ``proper'' morphism $\Phi:\HH\to\HH'$ induces a homomorphism $\Phi^*:H_{\text{\rm Ha}}(\HH')\to H_{\text{\rm Ha}}(\HH)$ in a functorial way. The commutativity of the diagram
$$
\begin{CD}
E_{c,1}^{\bullet,p}(\GG) @>{E_{c,1}^{\bullet,p}(f)}>> E_{c,1}^{\bullet,p}(\FF) \\
@V{\int_{\GG}}VV  @VV{\int_{\FF}}V \\
\Omega_{\text{\rm Ha}}(\Hol(\GG)) @>{\Hol(f)^*}>> \Omega_{\text{\rm Ha}}(\Hol(\FF))
\end{CD}
$$
must be required for any proper \cinf\ foliated map $f:\FF\to\GG$ between \cinf\ foliations of dimension $p$. If this can be made, then there must be a version of \reft{proper invariance} for the Haefliger cohomology.

Of course, first, the right definition of {\em proper morphism\/} must be given! We may say that a morphism $\Phi$ is {\em proper\/} if $\Phi_{\text{\rm orb}}$ is proper, but this does not seem to be enough. With a good definition of {\em proper morphism\/}, the canonical injective functor $\Top\to\PsGr$ and the holonomy functor $\Fol\to\PsGr$ must assign proper morphisms to proper maps. Thus this condition may have some relation with Haefliger's definition of {\em compact generation\/} for pseudogroups \cite{Haefliger:85Pseudogroups}, \cite{Haefliger:02compactly}.

\end{prob}

\begin{prob}
A {\em pseudomonoid\/} can be defined like a pseudogroup with arbitrary continuous maps and without using inversion. {\em Orbits\/}, {\em morphisms\/} and {\em equivalences\/} can be generalized to pseudomonoids. But, by the lack of inversion, it also makes sense to consider {\em $\alpha$-\/} and {\em $\omega$-orbits\/}, referring to the backwards and forwards direction. We similarly have {\em $\alpha$-\/} and {\em $\omega$-morphisms\/}, and {\em $\alpha$-\/} and {\em $\omega$-equivalences\/}. For instance, for any foliated space, all maps between local quotients induced by inclusions of simple open sets form a pseudomonoid. Its $\alpha$-class gives the holonomy pseudogroup of all open subsets. This may be useful to deal with invariants like the {\em transverse LS~category\/} $\cat_\pitchfork\FF$ of a foliated structure \FF, whose definition involves non-saturated open sets \cite{ColmanMacias:01transvLS}. For instance, to know how far it is from being a transverse invariant.
\end{prob}

\begin{prob}
The {\em saturated transverse LS~category\/} $\cat_\pitchfork^S\FF$ of a foliated structure \FF\ is defined like $\cat_\pitchfork\FF$, but using only saturated open sets. It is not a transverse invariant either, but it is closer to being so. The ``transverse invariant LS~category'' of $\FF$ is the LS~category $\cat\Hol(\FF)$ of its holonomy pseudogroup, defined in the sense of \refs{generalization}. We easily get 
$$
\cat\Hol(\FF)\le\cat_\pitchfork^S\FF\;,
$$ 
which may be a strict inequality by the possible non-existence of lifts of homotopies on holonomy pseudogroups (Problem~\ref{pb:lift}). We propose the study of $\cat\Hol(\FF)$. For instance, it should be equal to the ``transverse category'' of the classifying foliated space of \FF\ \cite{Haefliger:71Homotopy}, \cite{BuffetLor}, \cite{Haefliger:84classifiants}. Moreover the known results for  $\cat_\pitchfork\FF$ or  $\cat_\pitchfork^S\FF$ should have a version for $\cat\Hol(\FF)$ (see \cite{ColmanMacias:01transvLS}, \cite{ColmanHurder:Lcompact}, \cite{Colman:04transvLSRiem}, \cite{Hurder}, \cite{HurderWalczak:CompactLS}).
\end{prob}

\begin{prob}
Give a good definition of \coinf\ singular foliated spaces.
\end{prob}

\begin{prob}[Based on an idea of Jos\'e Luis Arraut]
The concept of {\em holonomy groupoid\/} can be extended to the singular case as follows.
Let \FF\ be a singular foliated structure on a space $X$, let $c:I\to\FF$ be a foliated curve between points $x$ and $y$, and let $\Sigma_x$ and $\Sigma_y$ be local transversals of \FF\ through $x$ and $y$, respectively. There is some open neighborhood $\Delta$ of $x$ in $T_x$ and some continuous foliated map $H:\Delta_{\text{\rm pt}}\times I\to\FF$ such that:
\begin{itemize}

\item $H(x,\cdot)=c$;

\item each map $H(\cdot,t)$ is an embedding whose image is a local transversal of \FF\ through $c(t)$;

\item $H(\cdot,0)$ is the inclusion map $\Delta\hookrightarrow X$; and

\item $H(\Delta\times\{1\})$ is an open subset of $\Sigma_y$.

\end{itemize}
Then $h=H(\cdot,1):\Delta\to\Sigma_y$ is an open embedding, but, in the singular case, the germ $\gamma(h,x)$ may depend on the choice of $H$. To avoid this dependence, we introduce an equivalence relation on the set of such germs as follows. First, observe that $h$ is a foliated map $\FF_\Delta\to\FF_{\Sigma_y}$ (see \refs{singular}). Then, given another continuous foliated map $H':\Delta'_{\text{\rm pt}}\times I\to\FF$ satisfying the same properties as $H$, for $h'=H'(\cdot,1):\Delta'\to\Sigma_y$, let as say that $\gamma(h,x)$ is {\em equivalent\/} to $\gamma(h',x)$ if there is some open neighborhood $\Delta_0$ of $x$ in $\Delta\cap\Delta'$ and some integrable homotopy $G:\FF_{\Delta_0}\times I\to\FF_{\Sigma_y}$ between $h$ and $h'$ such that each $G(\cdot,t)$ is an open embedding. This defines an equivalence relation, and it is easy to check that the equivalence class of $\gamma(h,x)$ depends only on $c$ (once $\Sigma_x$ and $\Sigma_y$ are give), and can be called the {\em holonomy\/} defined by $c$. On the set of foliated curves $I\to\FF$, define an equivalence relation by declaring that to curves are equivalent when they have the same end points and define the same holonomy for any choice of local transversals. The quotient set $\mathfrak{G}=\mathfrak{G}(\FF)$ becomes a topological groupoid with the operation induced by the path product and the topology induced by the compact-open topology; this $\mathfrak{G}$ can be called the {\em holonomy groupoid\/}, and can be identified to the usual holonomy groupoid when \FF\ is regular. 

Like in the regular case, another version of the holonomy groupoid can be also given by a singular foliated cocycle (\refs{singular}), by adapting the steps of the above construction of $\mathfrak{G}$.

It is a natural question whether the important role played by the holonomy groupoid of foliated spaces can be extended to the singular case. As a first approach, we may ask for a singular version of the Reeb's local stability theorem.
\end{prob}

\begin{prob}
Is it possible to generalize our approximation and invariance results to the case of singular Riemannian foliations? Such generalizations should be consequences of a singular version of \reft{main}.
\end{prob}

\end{document}